\documentclass[12pt]{amsart}
\usepackage{amscd,amssymb,amsthm,verbatim}
\newcommand{\p}{p} 
\newcommand{\pp}{p_2} 
\def\bp{\star} 

\newcommand{\const}{\operatorname{const.}}

\newcommand{\diam}{\operatorname{diam}}

\newcommand{\dvol}{\operatorname{dvol}}

\newcommand{\Hess}{\operatorname{Hess}}

\newcommand{\Id}{\operatorname{Id}}

\newcommand{\Lip}{\operatorname{Lip}}

\newcommand{\N}{{\mathbb N}}

\newcommand{\Q}{{\mathbb Q}}
\newcommand{\R}{{\mathbb R}}

\newcommand{\Ric}{\operatorname{Ric}}
\newcommand{\Riem}{\operatorname{Riem}}

\newcommand{\supp}{\operatorname{supp}}

\newcommand{\Tr}{\operatorname{Tr}}

\newcommand{\vol}{\operatorname{vol}}
\newcommand{\Z}{{\mathbb Z}}

\newcommand{\DC}{{\cal D}{\cal C}}
\numberwithin{equation}{section}
\theoremstyle{plain}
\newtheorem{definition}[equation]{Definition}

\newtheorem{lemma}[equation]{Lemma}
\newtheorem{theorem}[equation]{Theorem}
\newtheorem{proposition}[equation]{Proposition}
\newtheorem{corollary}[equation]{Corollary}

\newtheorem{partics}[equation]{Particular cases}
\theoremstyle{remark}

\newtheorem{remark}[equation]{Remark}
\newtheorem{example}[equation]{Example}

\usepackage{graphicx}
\input{amssym.def}
\input{amssym}
\input{epsf}
\usepackage{a4wide}

\def\R{\mathbb R}
\def\N{\mathbb N}

\def\Q{\mathbb Q}
\def\Z{\mathbb Z}

\def\var{\varepsilon}

\def\ov{\overline}
\def\cal{\mathcal}

\def\tilde{\widetilde}

\def\Hess{\mathop{{\rm Hess}\,}}

\def\Id{{\rm Id}\,}
\def\dps{\displaystyle}

\def\2dr#1#2{\left. \frac{d^2}{d{#1}^2} \right |_{#2}}
\def\d2#1{\frac{d^2}{d{#1}^2}}

\def\LSI{{\rm LSI}}

\DeclareMathOperator*{\grad}{grad}

\def\med{\medskip}
\def\sm{\smallskip}

\def\begeq{\begin{equation}}
\def\endeq{\end{equation}}
\def\begar{\begin{eqnarray}}
\def\endar{\end{eqnarray}}
\def\begar*{\begin{eqnarray*}}
\def\endar*{\end{eqnarray*}}
\def\begal{\begin{align}}
\def\endal{\end{align}}
\def\begal*{\begin{align*}}
\def\endal*{\end{align*}}

\theoremstyle{definition}

\theoremstyle{remark}

\newtheorem*{Thm*}{Theorem}
\newtheorem*{Lem*}{Lemma}
\newtheorem*{Conj*}{Conjecture}
\newtheorem*{Cor*}{Corollary}
\newtheorem*{Def*}{Definition}
\newtheorem*{Prop*}{Proposition}
\newtheorem*{Exo*}{Exercise}
\newtheorem*{Exs*}{Examples}
\newtheorem*{Ex*}{Example}
\newtheorem*{Rk*}{Remark}
\newtheorem*{Rks*}{Remarks}

\begin{document}

\title[Ricci curvature via optimal transport]
{Ricci curvature for metric-measure spaces via optimal transport}

\author{John Lott}
\address{Department of Mathematics\\
University of Michigan\\
Ann Arbor, MI  48109-1109\\
USA} \email{lott@umich.edu}
\author{C\'edric Villani}
\address{UMPA (UMR CNRS 5669)\\ ENS Lyon\\
46 all\'ee d'Italie, 69364 Lyon Cedex 07\\
FRANCE} \email{cvillani@umpa.ens-lyon.fr}

\thanks{The research of the first author was 
supported by NSF grant DMS-0306242 and the Clay Mathematics
Institute}
\date{June 23, 2006}

\begin{abstract}
We define a notion of a measured length space $X$ having
nonnegative $N$-Ricci curvature, for $N \in [1, \infty)$, or
having $\infty$-Ricci curvature bounded below by $K$, for
$K \in \R$. The definitions are in terms of the displacement
convexity of certain functions on the associated Wasserstein
metric space $P_2(X)$ of probability measures. We show that these 
properties are preserved under measured Gromov--Hausdorff limits.
We give geometric and analytic consequences.
\end{abstract}

\maketitle

This paper has dual goals.  One goal is to extend results about
optimal transport from the setting of smooth Riemannian manifolds to the 
setting of length spaces. A second goal is
to use optimal transport to give a notion for a measured
length space to have Ricci curvature bounded below.
We refer to \cite{Burago-Burago-Ivanov (2001)} and 
\cite{Villani (2003)} for background material on 
length spaces and optimal transport, respectively.
Further bibliographic notes on optimal transport
are in Appendix \ref{bibnotes}.
In the present introduction we motivate the questions that we address
and we state the main results.

To start on the geometric side, there are various reasons to
try to extend notions of curvature from smooth Riemannian manifolds
to more general spaces.  A fairly general setting is that of 
{\em length spaces},
meaning metric spaces $(X, d)$ 
in which the distance between two points equals
the infimum of the lengths of curves joining the points. In the
rest of this introduction
we assume that $X$ is a compact length space. 
Alexandrov gave a good notion of a length space having 
{\em ``curvature bounded below by $K$''}, with $K$ a real number,
in terms of the geodesic triangles in $X$. In the case of a Riemannian
manifold $M$ with the induced length structure, one recovers the
Riemannian notion of having sectional curvature bounded below by $K$.
Length spaces with Alexandrov curvature bounded below by $K$ behave
nicely with respect to the Gromov--Hausdorff topology on compact metric spaces
(modulo isometries); they form a closed subset.

In view of Alexandrov's work, it is natural to ask whether there are
metric space versions of other types of Riemannian curvature, such 
as Ricci curvature.  This question takes substance from {\em Gromov's
precompactness theorem} for Riemannian manifolds with Ricci
curvature bounded below by $K$, dimension bounded above by $N$ 
and diameter bounded above by $D$
\cite[Theorem 5.3]{Gromov (1999)}. 
The precompactness indicates that there could be a
notion of a length space having ``Ricci curvature bounded below by
$K$'', special cases of which would be Gromov--Hausdorff limits of 
manifolds with lower Ricci curvature bounds.

Gromov--Hausdorff limits of manifolds with Ricci curvature bounded
below have been studied by various authors, notably Cheeger
and Colding \cite{Cheeger-Colding (1996),Cheeger-Colding 
(1997),Cheeger-Colding II (2000),Cheeger-Colding (2000)}. One feature
of their work, along with the earlier work of 
Fukaya \cite{Fukaya (1987)}, is that it turns out to be useful to
add an auxiliary Borel probability measure $\nu$ 
and consider {\em metric-measure spaces} $(X,d,\nu)$. (A compact 
Riemannian manifold $M$ has a canonical measure $\nu$ given by the
normalized Riemannian density $\frac{\dvol_M}{\vol(M)}$.) There is a measured
Gromov--Hausdorff topology on such triples $(X,d,\nu)$ (modulo isometries)
and one again has precompactness for 
Riemannian manifolds with Ricci
curvature bounded below by $K$, dimension bounded above by $N$
and diameter bounded above by $D$.
Hence the question is whether there is a good notion of a measured
length space $(X, d, \nu)$ having {\em ``Ricci curvature bounded below by
$K$''}. Whatever definition one takes,
one would like the set of such triples to be closed in the
measured Gromov--Hausdorff topology. One would also like to
derive some nontrivial consequences from the definition, and of course
in the case of Riemannian manifolds one would like to recover 
classical notions. We refer to
\cite[Appendix 2]{Cheeger-Colding (1997)} for further discussion of 
the problem of giving a ``synthetic'' treatment of Ricci curvature.

Our approach is in terms of
a metric space $(P(X), W_2)$ that is canonically
associated to the original metric space
$(X, d)$. Here $P(X)$ is the space of Borel probability measures
on $X$ and $W_2$ is the so-called {\em Wasserstein distance} of order $2$.
The square of the Wasserstein distance $W_2(\mu_0, \mu_1)$ between
$\mu_0, \mu_1 \in P(X)$ is defined to be
the infimal cost to transport the total mass
from the measure $\mu_0$ to the measure $\mu_1$, where the
cost to transport a unit of mass between points $x_0, x_1 \in X$ is
taken to be $d(x_0, x_1)^2$. A transportation scheme with infimal cost is
called an {\em optimal transport}. The topology on $P(X)$ coming from the
metric $W_2$ turns out to be the weak-$*$ topology. We will write
$P_2(X)$ for the metric space $(P(X), W_2)$, which we call
the {\em Wasserstein space}. If $(X, d)$ is a length space then
$P_2(X)$ turns out to also be a length space. Its geodesics will be
called {\em Wasserstein geodesics}. If $M$ is a Riemannian
manifold then we write $P^{ac}_2(M)$ for the elements of 
$P_2(M)$ that are absolutely continuous with respect to the
Riemannian density.

In the past fifteen years, optimal transport of measures has been extensively 
studied in the case $X = \R^n$, with motivation coming from the study
of certain partial differential equations. A notion which has proved
useful is that of {\em displacement convexity},
i.e. convexity along Wasserstein geodesics, which
was introduced by McCann in order to show the existence
and uniqueness of minimizers for certain relevant functions on $P_2^{ac}(\R^n)$
\cite{McCann (1997)}.

In the past few years, some regularity results for optimal transport on $\R^n$ 
have been extended to Riemannian manifolds
\cite{Cordero-Erausquin-McCann-Schmuckenschlager (2001),McCann (2001)}.
This made it possible to study displacement convexity in a Riemannian
setting. Otto and Villani \cite{Otto-Villani (2000)}
carried out Hessian computations for certain functions on $P_2(M)$
using a formal infinite-dimensional Riemannian structure on $P_2(M)$ defined by
Otto \cite{Otto (2001)}.
These formal computations indicated a relationship between 
the Hessian of an ``entropy'' function on $P_2(M)$ and
the Ricci curvature of $M$. Later, a rigorous
displacement convexity result for a class of functions on $P^{ac}_2(M)$, 
when $M$ has nonnegative Ricci curvature, was proven by
Cordero-Erausquin, McCann and Schmuckenschl\"ager
\cite{Cordero-Erausquin-McCann-Schmuckenschlager (2001)}.
This work was extended by von Renesse and Sturm \cite{Sturm-von Renesse}.

Again in the case of Riemannian manifolds, a further circle of ideas 
relates displacement convexity to log Sobolev inequalities,
Poincar\'e inequalities, Talagrand inequalities and concentration of
measure \cite{Bobkov-Gentil-Ledoux (2001),Bobkov-Gotze (1999),Ledoux 
(2001),Otto-Villani (2000)}.

In this paper we use optimal transport and displacement convexity
in order to {\em define} a notion of a measured length space $(X,d,\nu)$ 
having Ricci curvature bounded below.  
If $N$ is a finite parameter (playing the role of a dimension)
then we will define a notion of $(X,d,\nu)$ having nonnegative 
$N$-Ricci curvature. 
We will also define a notion of $(X,d,\nu)$ having $\infty$-Ricci curvature 
bounded below by $K \in \R$.
(The need to input the 
possibly-infinite parameter $N$ can be seen from the Bishop--Gromov
inequality for complete
$n$-dimensional Riemannian manifolds with nonnegative Ricci
curvature. It states that $r^{-n} \vol(B_r(m))$ is nonincreasing in $r$,
where $B_r(m)$ is the $r$-ball centered at $m$
\cite[Lemma 5.3.bis]{Gromov (1999)}. When we go from manifolds to 
length spaces there is no {\it a priori} value for the parameter $n$.
This indicates the need to specify a dimension parameter in the
definition of Ricci curvature bounds.)

We now give the main results of the paper, sometimes in a simplified form.
For consistency, we assume in the body of the paper that the relevant length
space $X$ is compact. The necessary modifications to deal with complete pointed
locally compact length spaces are given in Appendix \ref{noncompact}.

Let $U \: : \: [0, \infty) \rightarrow \R$ be a continuous
convex function with $U(0) \: = \: 0$.
Given a reference probability measure $\nu \in P(X)$,
define the function
$U_\nu \: : \: P_2(X) \rightarrow \R \cup \{\infty\}$ by
\begeq
U_\nu(\mu) = 
\int_X U(\rho(x))\,d\nu(x) + U'(\infty)\,\mu_s(X), 
\endeq
where 
\begin{equation} 
\mu = \rho \nu + \mu_s 
\end{equation}
is the Lebesgue decomposition of $\mu$ with respect to $\nu$ 
into an absolutely continuous part $\rho\nu$ and a singular part $\mu_s$, and
\begin{equation}
U^\prime(\infty) \: = \: \lim_{r \rightarrow \infty} \frac{U(r)}{r}.
\end{equation}

If $N\in [1,\infty)$ then 
we define $\DC_N$ to be the set of such functions $U$
so that the function
\begin{equation} 
\psi(\lambda) = \lambda^N \:U(\lambda^{-N}) 
\end{equation}
is convex on $(0, \infty)$.
We further define $\DC_\infty$ 
to be the set of such functions $U$ so that the function
\begin{equation}
\psi(\lambda) = e^\lambda \: U(e^{-\lambda})
\end{equation} 
is convex on $(- \infty, \infty)$.
A relevant example of an element of $\DC_N$ is given by
\begin{equation}
U_N(r) \: = \: 
\begin{cases}
Nr (1 - r^{-1/N})& \text{ if $1 < N < \infty$}, \\
r \log{r} & \text{ if $N = \infty$}.
\end{cases}
\end{equation}

\begin{definition}
Given $N \in [1, \infty]$, we say that a compact measured length space
$(X,d,\nu)$ has nonnegative $N$-Ricci curvature if for all
$\mu_0, \mu_1 \in P_2(X)$ with
$\supp(\mu_0) \subset \supp(\nu)$ and $\supp(\mu_1) \subset \supp(\nu)$, 
there is {\em some} Wasserstein geodesic
$\{\mu_t\}_{t \in [0,1]}$ from $\mu_0$ to $\mu_1$ so that
for all $U \in \DC_N$ and all $t \in [0,1]$,
\begin{equation} \label{ineq1}
U_\nu(\mu_t) \: \le \: t \: U_\nu(\mu_1) \: + \: (1-t) \: 
U_\nu(\mu_0).
\end{equation}
Given $K \in \R$, we say that $(X,d,\nu)$ has $\infty$-Ricci curvature
bounded below by $K$ if for all
$\mu_0, \mu_1 \in P_2(X)$ with
$\supp(\mu_0) \subset \supp(\nu)$ and $\supp(\mu_1) \subset \supp(\nu)$,
there is {\em some} Wasserstein geodesic
$\{\mu_t\}_{t \in [0,1]}$ from $\mu_0$ to $\mu_1$ so that
for all $U \in \DC_\infty$ and all $t \in [0,1]$,
\begin{equation} \label{ineq2}
U_\nu(\mu_t) \: \le \: t \: U_\nu(\mu_1) \: + \: (1-t) \: 
U_\nu(\mu_0)
 \: - \:
\frac12 \: \lambda(U) \: t(1-t) W_2(\mu_0, \mu_1)^2,
\end{equation}
where $\lambda \: : \: \DC_\infty \rightarrow \R \cup \{-\infty\}$
is defined in~\eqref{lambdainfty} below.
\end{definition}

Note that the inequalities (\ref{ineq1}) and (\ref{ineq2}) are only assumed to
hold along {\em some} Wasserstein geodesic from $\mu_0$ to $\mu_1$, 
and not necessarily along all such geodesics.  This is what we call
{\em weak displacement convexity}.

Naturally, one wants to know that in the case of a Riemannian manifold,
our definitions are equivalent to classical ones.  Let $M$ be
a smooth compact connected $n$-dimensional manifold with Riemannian
metric $g$. We let $(M, g)$ denote the corresponding metric space. 
Given $\Psi \in C^\infty(M)$ with $\int_M e^{-\Psi} \: \dvol_M \: = \: 1$,
put $d\nu \: = \: e^{-\Psi} \: \dvol_M$.

\begin{definition}
For $N \in [1, \infty]$, let the $N$-Ricci tensor $\Ric_N$ of $(M, g, \nu)$
be defined by
\begin{equation}
\Ric_N \: = 
\begin{cases}
\Ric \: + \: \Hess(\Psi)  & \text{ if $N = \infty$}, \\
\Ric \: + \: \Hess(\Psi) \: - \: \frac{1}{N-n} \: d \Psi \otimes
d \Psi & \text{ if $n \: < \: N \: < \: \infty$}, \\
\Ric \: + \: \Hess(\Psi) \: - \: \infty \: (d \Psi \otimes
d \Psi) & \text{ if $N = n$}, \\
- \infty & \text{ if $N < n$,}
\end{cases}
\end{equation}
where by convention $\infty \cdot 0 \: = \: 0$.
\end{definition}

\begin{theorem} \label{Riemequiv}
(a) For $N \in [1, \infty)$, the measured length space 
$(M, g, \nu)$ has nonnegative $N$-Ricci curvature if and
only if $\Ric_N \: \ge \: 0$. \\
(b) $(M, g, \nu)$ has $\infty$-Ricci curvature bounded below
by $K$ if and only if $\Ric_\infty \: \ge \: K g$.
\end{theorem}

In the special case when $\Psi$ is constant, and so
$\nu \: = \: \frac{\dvol_M}{\vol(M)}$,
Theorem \ref{Riemequiv} shows that we recover the usual notion
of a Ricci curvature bound from our length space definition
as soon as $N\geq n$.

The next theorem, which is the 
main result of the paper, says that our notion of $N$-Ricci
curvature has good behavior under measured Gromov--Hausdorff limits.

\begin{theorem} \label{main}
Let $\{(X_i, d_i, \nu_i)\}_{i=1}^\infty$ be a sequence of
compact measured length spaces with $\lim_{i \rightarrow \infty} 
(X_i,d_i,\nu_i) \: = \: (X,d,\nu)$ in the
measured Gromov--Hausdorff topology. \\
(a) For any $N \in [1, \infty)$, if 
each $(X_i, d_i, \nu_i)$ has nonnegative $N$-Ricci curvature then
$(X, d, \nu)$ has nonnegative $N$-Ricci curvature.\\
(b) If each $(X_i, d_i, \nu_i)$ has $\infty$-Ricci curvature bounded
below by $K$ then $(X, d, \nu)$ has $\infty$-Ricci curvature bounded
below by $K$.
\end{theorem}

Theorems \ref{Riemequiv} and \ref{main} imply that measured Gromov--Hausdorff
limits $(X,d,\nu)$ of smooth manifolds
$\left( M, g, \frac{\dvol_M}{\vol(M)} \right)$ with lower Ricci
curvature bounds fall under our
considerations. Additionally, we obtain
the following new characterization of such limits $(X,d,\nu)$ 
which happen to be {\em smooth}, meaning that $(X, d)$ is a smooth
$n$-dimensional
Riemannian manifold $(B, g_B)$ and $d\nu \: = \: e^{-  \Psi} \: \dvol_B$ for
some $\Psi \in C^\infty(B)$ :

\begin{corollary} \label{maincor}
(a) If $(B, g_B, \nu)$ is a measured 
Gromov--Hausdorff limit of Riemannian manifolds with nonnegative
Ricci curvature and dimension at most $N$ then $\Ric_N(B) \: \ge \:0$. \\
(b) If $(B, g_B, \nu)$ is a measured 
Gromov--Hausdorff limit of Riemannian manifolds with Ricci curvature
bounded below by $K \in \R$ then $\Ric_\infty(B) \: \ge \: K \: g_B$.
\end{corollary}

There is a partial converse to Corollary \ref{maincor}
(see Corollary~\ref{limitspacecor}(ii,ii')).

Finally, if a measured length space has lower Ricci curvature bounds then there
are analytic consequences, such as a \emph{log Sobolev inequality}.
To state it, we define the {\em gradient norm} of a Lipschitz function $f$
on $X$ by the formula
\begin{equation}
|\nabla f| (x) \: = \:
\limsup_{y\to x} \frac{|f(y)-f(x)|}{d(x,y)}.
\end{equation}

\begin{theorem} \label{introls}
Suppose that a compact measured length space
$(X,d,\nu)$ has $\infty$-Ricci curvature bounded below by
$K \in \R$. Suppose that
$f \in \Lip(X)$ satisfies $\int_X f^2 \: d\nu \: = \: 1$.\\
(a) If $K > 0$ then
\begin{equation} \label{introeqn}
\int_X f^2 \: \log(f^2) \: d\nu \: \le \: \frac{2}{K} \:
\int_X |\nabla f|^2 \: d\nu.
\end{equation}
\noindent
(b) If $K \: \le \: 0$ then
\begin{equation}
\int_X f^2 \: \log(f^2) \: d\nu \: \le \: 2 \: \diam(X) \:
\sqrt{\int_X |\nabla f|^2 \: d\nu} \: - \: \frac12 \: 
K \: \diam(X)^2.
\end{equation}
\end{theorem}

In the case of Riemannian manifolds, one recovers from 
(\ref{introeqn}) the log Sobolev inequality of Bakry and
\'Emery \cite{Bakry-Emery (1985)}. 

A consequence of (\ref{introeqn}) is a
Poincar\'e inequality.

\begin{corollary} \label{introcor}
Suppose that a compact measured length space
$(X,d,\nu)$ has $\infty$-Ricci curvature bounded below by
$K > 0$. Then for all
$h \in \Lip(X)$ with $\int_X h\,d\nu = 0$, we have
\begin{equation} 
\int_X h^2\,d\nu \leq 
\frac1{K} \int_X |\nabla h|^2\,d\nu. 
\end{equation}
\end{corollary}

In the case of Riemannian manifolds, Corollary \ref{introcor}
follows from the Lichnerowicz inequality for the smallest
positive eigenvalue of the Laplacian \cite{Lichnerowicz (1958)}.

We now give the structure of the paper.  More detailed descriptions
appear at the beginnings of the sections.

Section \ref{notation} gives basic definitions about length spaces
and optimal transport.  Section \ref{appsecgeom} 
shows that the Wasserstein space
of a length space is also a length space, and that Wasserstein geodesics
arise from displacement interpolations. Section \ref{secfunct} defines
weak displacement convexity and its variations. This is used to prove
functional inequalities called the HWI inequalities.

Section \ref{secGH} proves, modulo the technical results of Appendices
\ref{applsc} and \ref{secmoll}, 
that weak displacement convexity is preserved by
measured Gromov--Hausdorff limits.  The notion of $N$-Ricci
curvature is defined in Section \ref{secricci}, which contains the proof of
Theorem \ref{main}, along with a Bishop--Gromov-type inequality. Section
\ref{secineq} proves log Sobolev, Talagrand and Poincar\'e inequalities
for measured length spaces, such as Theorem \ref{introls} and
Corollary \ref{introcor},
along with 
a weak Bonnet--Myers theorem.  Section \ref{Riemcase} looks at the case
of smooth Riemannian manifolds and proves, in particular,
Theorem \ref{Riemequiv} and Corollary \ref{maincor}.

There are six appendices that contain either technical results or
auxiliary results.  Appendix \ref{secgeom}, which is a sequel to
Section \ref{appsecgeom}, discusses the geometry of the Wasserstein
space of a Riemannian manifold $M$.  
It shows that if $M$ has nonnegative sectional curvature
then $P_2(M)$ has nonnegative Alexandrov curvature.  The tangent
cones at absolutely continuous measures are computed, thereby 
making rigorous the formal Riemannian metric on $P_2(M)$ introduced by Otto.

Appendices \ref{applsc} and \ref{secmoll} 
are the technical core of Theorem \ref{main}.
Appendix \ref{applsc} shows that $U_{\nu}(\mu)$ is lower semicontinuous in 
both $\mu$ and $\nu$, 
and is nonincreasing under pushforward of
$\mu$ and $\nu$. Appendix \ref{secmoll} shows that a measure
$\mu \in P_2(X)$ with $\supp(\mu) \subset \supp(\nu)$ can
be weak-$*$ approximated by measures $\{\mu_k\}_{k=1}^\infty$
with continuous densities (with respect to $\nu$) so that 
$U_{\nu}(\mu) \: = \: \lim_{k \rightarrow \infty} U_{\nu}(\mu_k)$.

Appendix \ref{appformaldc} contains formal computations of the Hessian
of $U_\nu$. Appendix \ref{noncompact} explains how to
extend the results of the paper from the setting of compact measured
length spaces to the setting of 
complete pointed locally compact measured length spaces.
Appendix \ref{bibnotes} has some bibliographic notes on optimal
transport and displacement convexity.

The results of this paper were presented at the workshop ``Collapsing 
and metric
geometry'' in M\"unster, August 1-7, 2004.
After the  writing of the paper was essentially completed
we learned of related work by
Karl-Theodor Sturm \cite{Sturm1,Sturm2}. Also, Ludger R\"uschendorf
kindly pointed out to us that Theorem~\ref{thmlsc} was already proven in
\cite[Chapter~1]{Liese-Vajda (1987)} by different means.
We decided to retain our proof of 
Theorem~\ref{thmlsc} rather than just quoting~\cite{Liese-Vajda (1987)}, 
partly because the method of proof may be of
independent interest, partly for completeness and convenience to 
the reader, and partly because our method of proof is used
in the extension of the theorem considered in Appendix~\ref{noncompact}.

We thank MSRI and the UC-Berkeley mathematics department for their hospitality
while part of this research was performed. 
We also thank the anonymous referees for their suggestions.

\tableofcontents

\section{Notation and basic definitions} \label{notation}

In this section we first recall some facts about convex functions.
We then define gradient norms, length spaces and measured Gromov--Hausdorff
convergence.  Finally, we define the 2-Wasserstein metric $W_2$ on $P(X)$.

\subsection{Convex analysis}

Let us recall a few results from convex analysis.
See \cite[Chapter 2.1]{Villani (2003)} and references therein for
further information. 

Given a convex lower semicontinuous function 
$U \: : \: \R\to \R \cup\{\infty\}$ (which we assume is
not identically $\infty$),
its Legendre transform $U^* \: : \: 
\R\to \R \cup\{\infty\}$ is defined by
\begin{equation} 
U^*(p) = \sup_{r\in\R} \: \bigl[ pr - U(r)\bigr].
\end{equation}
Then $U^*$ is also convex and lower semicontinuous.
We will sometimes identify a convex lower semicontinuous function $U$ defined
on a closed interval $I \subset \R$ with the convex function defined 
on the 
whole of $\R$ by extending $U$ by $\infty$ outside of $I$.

Let $U \: : \: [0, \infty) \rightarrow \R$ be a convex lower
semicontinuous function. Then
$U$ admits a left derivative $U^\prime_- \: : \: (0, \infty) \rightarrow
\R$ and a right derivative $U^\prime_+ \: : \: [0, \infty) \rightarrow
\{-\infty\} \cup \R$, with $U^\prime_+(0,\infty) \subset \R$.
Furthermore, $U^\prime_- \: \le \: U^\prime_+$. They agree almost
everywhere and are both nondecreasing.
We will write 
\begin{equation}
U'(\infty) = \lim_{r\to\infty} U'_+(r) = 
\lim_{r\to\infty} \frac{U(r)}{r} \in \R \cup \{\infty\}.
\end{equation}
If we extend $U$ by $\infty$ on $(-\infty, 0)$ then its Legendre transform
$U^* \: : \: 
\R\to \R \cup\{\infty\}$ becomes
$U^*(p) \:= \: \sup_{r \ge 0} \: \bigl[ pr - U(r)\bigr]$.
It is nondecreasing in $p$,
infinite on $(U^\prime(\infty), \infty)$ and
equals $-U(0)$ on $(- \infty, U^\prime_+(0)]$.
Furthermore, it is continuous on $(-\infty, U^\prime(\infty))$.
For all $r \in [0, \infty)$,
we have $U^*(U^\prime_+(r)) \: = \: r \: U^\prime_+(r) \: - \: U(r)$.

\subsection{Geometry of metric spaces}

\subsubsection{Gradient norms}

Let $(X, d)$ be a compact metric space (with $d$ valued in $[0, \infty)$).
The open ball of radius
$r$ around $x \in X$ will be denoted by $B_r(x)$ and the sphere of radius $r$
around $x$ will be denoted by $S_r(x)$.

Let $L^\infty(X)$ denote the set of bounded measurable functions
on $X$. (We will consider such a function to be defined everywhere.)
Let $\Lip(X)$ denote the set of Lipschitz functions on $X$. 
Given $f \in \Lip(X)$, we define the {\em gradient norm} of $f$ by
\begin{equation} \label{nablaf}
|\nabla f| (x) \: = \:
\limsup_{y\to x} \frac{|f(y)-f(x)|}{d(x,y)}
\end{equation}
if $x$ is not an isolated point, and $|\nabla f| (x) \: = \: 0$ if
$x$ is isolated.
Then $|\nabla f| \in L^\infty(X)$.

On some occasions we will use a finer notion of gradient norm:
 \begin{equation} \label{nablamoins}
|\nabla^- f|(x) \: = \: \limsup_{y\to x} \frac{[f(y)-f(x)]_-}{d(x,y)} \: = \:
\limsup_{y\to x} \frac{[f(x)-f(y)]_+}{d(x,y)}
\end{equation} 
if $x$ is not isolated, and $|\nabla^- f|(x) \: = \: 0$ if $x$ is isolated.
Here $a_+ \: = \: \max(a, 0)$ and $a_- \: = \: \max(-a, 0)$.
Clearly $|\nabla^- f|(x) \: \le \: |\nabla f|(x)$.
Note that $|\nabla^- f|(x)$ is automatically zero if $f$ has a 
local minimum at $x$. In a sense,
$|\nabla^- f|(x)$ measures the downward pointing component of
$f$ near $x$.

\subsubsection{Length spaces} \label{lengthspaces}

If $\gamma$ is a curve in $X$, i.e. a continuous map
$\gamma \: : \: [0, 1] \rightarrow X$, then its length is
\begin{equation}
L(\gamma) \: = 
\: \sup_{J\in \N}\; \sup_{0 = t_0 \le t_1 \le \ldots \le t_J = 1}
\sum_{j=1}^{J} d\bigl(\gamma(t_{j-1}), \gamma(t_{j})\bigr).
\end{equation}
{From} the triangle inequality, $L(\gamma) \: \ge \: d(\gamma(0), \gamma(1))$.

We will assume that $X$ is a {\em length space}, meaning that the
distance between two points $x_0, x_1 \in X$ is the infimum of the
lengths of curves from $x_0$ to $x_1$. Such a space is path connected.

As $X$ is compact, it is a strictly intrinsic length space, 
meaning that we can replace infimum by 
minimum~\cite[Theorem 2.5.23]{Burago-Burago-Ivanov (2001)}. 
That is, for any $x_0, x_1 \in X$, there
is a minimal geodesic (possibly nonunique) from $x_0$ to $x_1$.
We may sometimes write ``geodesic'' instead of ``minimal geodesic''. 

By \cite[Proposition 2.5.9]{Burago-Burago-Ivanov (2001)},
any minimal geodesic $\gamma$ joining $x_0$ to $x_1$ can be
parametrized uniquely by $t\in [0, 1]$ so that
\begin{equation}
d(\gamma(t), \gamma(t^\prime)) \: = \: |t - t^\prime | \: d(x_0,x_1).
\end{equation}
We will often assume that the geodesic has been so parametrized.

By definition, a subset $A \subset X$ is {\em convex} if for 
any $x_0, x_1 \in A$ there is a minimizing geodesic from $x_0$ to $x_1$ 
that lies entirely in $A$. It is {\em totally convex} if for any
$x_0, x_1 \in A$, any minimizing geodesic in $X$ from $x_0$ to
$x_1$ lies in $A$.

Given $\lambda\in\R$, a function $F:X\to \R\cup\{\infty\}$ is said
to be {\em $\lambda$-convex} if for any geodesic
$\gamma \: : \: [0,1] \rightarrow X$ 
and any $t\in [0,1]$, we have
\begin{equation} F(\gamma(t)) \leq t F(\gamma(1)) + (1-t) F(\gamma(0))
- \frac12 \: \lambda \: t(1-t)\, L(\gamma)^2. \end{equation}
In the case when $X$ is a smooth Riemannian manifold with
Riemannian metric $g$, and $F \in C^2(X)$,
this is the same as saying that $\Hess F\geq \lambda g$.

\subsubsection{(Measured) Gromov--Hausdorff convergence} \label{subdefGH}

\begin{definition} \label{defGH}
Given two compact metric spaces $(X_1,d_1)$ and $(X_2,d_2)$, 
an $\epsilon$-Gromov--Hausdorff 
approximation from $X_1$ to $X_2$ is a (not necessarily continuous)
map $f \: : \: X_1 \rightarrow X_2$ so that 

(i)  For all $x_1, x_1^\prime \in X_1$, 
$\bigl|d_{2}(f(x_1), f(x_1^\prime)) \: - \: d_{1}(x_1, x_1^\prime)\bigr|
\: \le \: \epsilon$. 

(ii) For all $x_2 \in X_2$, there is an $x_1 \in X_1$ so that
$d_{2}(f(x_1), x_2) \: \le \: \epsilon$.
\end{definition}

An $\epsilon$-Gromov--Hausdorff 
approximation $f \: : \: X_1 \rightarrow X_2$ has an
approximate inverse $f^\prime \: : \: X_2 \rightarrow X_1$, which
can be constructed as follows:
Given $x_2 \in X_2$, choose $x_1 \in X_1$ so that
$d_{2}(f(x_1), x_2) \: \le \: \epsilon$ and put 
$f^\prime(x_2) \: = \: x_1$. Then $f^\prime$ is
a $3 \epsilon$-Gromov--Hausdorff 
approximation from $X_2$ to $X_1$. Moreover,
for all $x_1 \in X_1$, $d_{1}(x_1, (f^\prime \circ f)(x_1)) 
\: \le \: 2\epsilon$, and for all $x_2 \in X_2$, 
$d_{2}(x_2, (f \circ f^\prime)(x_2)) 
\: \le \: \epsilon$.

\begin{definition}
A sequence of compact metric spaces
$\{X_i\}_{i=1}^{\infty}$ converges to $X$ in the 
{\em Gromov--Hausdorff topology} if there is
a sequence of $\epsilon_i$-approximations $f_i \: : \: X_i \rightarrow
X$ with $\lim_{i \rightarrow \infty} \epsilon_i \: = \: 0$.
\end{definition}

This notion of convergence comes from a metrizable topology on the
space of all compact metric spaces modulo isometries.
If $\{X_i\}_{i=1}^\infty$ are length spaces that converge to $X$ in the
Gromov--Hausdorff topology then $X$ is also a length space 
\cite[Theorem 7.5.1]{Burago-Burago-Ivanov (2001)}.

For the purposes of this paper, 
we can and will assume that the maps $f$ and $f^\prime$ in
Gromov--Hausdorff approximations are Borel.
Let $P(X)$ denote the space of Borel probability measures on $X$.
We give $P(X)$ the weak-$*$ topology,
i.e. $\lim_{i \rightarrow \infty} \mu_i \: = \: \mu$ if and only if for all
$F \in C(X)$, $\lim_{i \rightarrow \infty} \int_X F \: d\mu_i \: = \: 
\int_X F \: d\mu$.

\begin{definition}
Given $\nu \in P(X)$, 
consider the metric-measure space $(X,d,\nu)$. A sequence
$\{(X_i, d_i, \nu_i)\}_{i=1}^\infty$ converges to $(X,d, \nu)$ in the
{\em measured Gromov--Hausdorff topology} 
if there are $\epsilon_i$-approximations $f_i:X_i\to X$, with
$\lim_{i \rightarrow \infty} \epsilon_i \: = \: 0$, so that
$\lim_{i \rightarrow \infty} (f_i)_* \nu_i \: = \: \nu$ in $P(X)$.
\end{definition}

Other topologies on the class of metric-measure spaces are discussed in
\cite[Chapter $3 \frac12$]{Gromov (1999)}.

For later use we note the following 
generalization of the Arzel\`a--Ascoli theorem.

\begin{lemma} (cf. \cite[p. 66]{Gromov (1981)}, 
\cite[Appendix A]{Grove-Petersen (1991)}) 
\label{Arzela-Ascoli}
Let $\{X_i\}_{i=1}^\infty$ be a sequence of compact metric spaces
converging to $X$ in the Gromov--Hausdorff topology, with
$\epsilon_i$-approximations $f_i \: : \: X_i \rightarrow X$.
Let
$\{Y_i\}_{i=1}^\infty$ be a sequence of compact metric spaces
converging to $Y$ in the Gromov--Hausdorff topology, with
$\epsilon_i$-approximations $g_i \: : \: Y_i \rightarrow Y$. 
For each $i$, let $f_i^\prime \: : \: X \rightarrow X_i$ be an approximate
inverse to $f_i$, as in the paragraph following Definition
\ref{defGH}.
Let $\{\alpha_i\}_{i=1}^\infty$ be a sequence of maps
$\alpha_i \: : \: X_i \rightarrow Y_i$ that
are asymptotically equicontinuous in the
sense that for every $\epsilon > 0$, there are $\delta = \delta(\epsilon) > 0$
and $N = N(\epsilon) \in \Z^+$ so that
for all $i \ge N$, 
\begin{equation}
d_{X_i}(x_i, x_i^\prime) < \delta \qquad \Longrightarrow\qquad
d_{Y_i}(\alpha_i(x_i), \alpha_i(x_i^\prime)) < \epsilon.
\end{equation}
Then after passing to a subsequence, the
maps $g_i \circ \alpha_i \circ f_i^\prime \: : \: X \rightarrow Y$ 
converge uniformly
to a continuous map $\alpha \: : \: X \rightarrow Y$.
\end{lemma}

In the conclusion of Lemma \ref{Arzela-Ascoli} the maps 
$g_i \circ \alpha_i \circ f_i^\prime$ may not be continuous, but
the notion of uniform convergence makes sense nevertheless.

\subsection{Optimal transport: basic definitions}

Given $\mu_0, \mu_1 \in P(X)$,
we say that a probability measure $\pi \in P(X\times X)$ is
a {\em transference plan} between $\mu_0$ and $\mu_1$ if
\begin{equation} 
(p_0)_* \pi \: = \: \mu_0,\qquad (p_1)_* \pi \: = \: \mu_1,
\end{equation}
where $p_0, p_1 \: : \: X \times X \rightarrow X$ are projections
onto the first and second factors, respectively. In words,
$\pi$ represents a way to transport the mass from $\mu_0$ to $\mu_1$,
and $\pi(x_0,x_1)$ is the amount of mass which is taken from a point $x_0$
and brought to a point $x_1$.

We will use optimal transport with quadratic cost function
(square of the distance). Namely, given $\mu_0, \mu_1 \in P(X)$, we 
consider the variational problem
\begin{equation} \label{var}
W_2(\mu_0, \mu_1)^2 \: = \:
\inf_{\pi}  \: \int_{X \times X} d(x_0, x_1)^2 \: d\pi(x_0, x_1),
\end{equation}
where $\pi$ ranges over the set of all transference plans between
$\mu_0$ and $\mu_1$. Any minimizer $\pi$ for this variational problem
is called an {\em optimal transference plan}.

In \eqref{var}, one can replace the infimum by the minimum
\cite[Proposition 2.1]{Villani (2003)}, i.e. there always exists
(at least) one optimal transference plan. Since $X$ has finite
diameter, the infimum is obviously finite. The quantity $W_2$ will
be called the {\em Wasserstein distance of order~2} between $\mu_0$
and $\mu_1$; it defines a metric on $P(X)$. 
The topology that it induces on $P(X)$ is the
weak-$*$ topology~\cite[Theorems 7.3 and 7.12]{Villani (2003)}.
When equipped with the metric $W_2$, $P(X)$ is a compact metric space,
which we will often denote by $P_2(X)$.

We remark that there is an isometric embedding $X \rightarrow P_2(X)$ given by
$x \rightarrow \delta_x$. This shows that $\diam(P_2(X)) \: \geq \:
\diam(X)$. Since the reverse inequality follows from the definition
of $W_2$, actually
$\diam(P_2(X)) \: = \: \diam(X)$.

A {\em Monge transport} is a transference plan coming from a map
$F \: : \: X \rightarrow X$ with $F_* \mu_0 \: = \: \mu_1$, given
by $\pi \: = \: (\Id, F)_* \mu_0$. In general an optimal transference
plan does not have to be a Monge transport, although this may
be true under some assumptions (as we will recall below).

A function $\phi \: : \: X \rightarrow [- \infty, \infty)$ is 
{\em $\frac{d^2}{2}$-concave} if it is not identically $- \infty$ and it can
be written in the form
\begin{equation}
\phi(x) \: = \: \inf_{x^\prime \in X} 
\left( \frac{d(x, x^\prime)^2}{2} \: - \: \widetilde{\phi}(x^\prime) \right) 
\end{equation}
for some function $\widetilde{\phi} \: : \: X \rightarrow [-\infty, \infty)$.
Such functions play an important role in the description of optimal
transport on Riemannian manifolds.

\section{Geometry of the Wasserstein space} \label{appsecgeom}

In this section, we investigate some features of 
the Wasserstein space $P_2(X)$ associated to a compact length space $(X,d)$.
(Recall that the subscript 2 in $P_2(X)$ means that $P(X)$ is
equipped with the 2-Wasserstein metric.) We show that $P_2(X)$ is a
length space. We define displacement interpolations
and show that every Wasserstein geodesic comes from a displacement
interpolation. We then recall some facts about optimal transport on
Riemannian manifolds.  

\subsection{Displacement interpolations}

We denote by $\Lip([0,1],X)$ the space of Lipschitz 
continuous maps $c\: : \: [0,1]
\rightarrow X$ with the uniform topology. For any $k > 0$,
\begin{equation} \label{lipk}
\Lip_k([0,1], X) \: = \:
\Bigl\{c \in \Lip([0,1], X) \: : \: d(c(t), c(t^\prime)) \: \le \: 
k |t-t^\prime| \mbox{ for all } t, t^\prime \in [0,1]\Bigr\}
\end{equation}
is a compact subset of  $\Lip([0,1],X)$.

Let $\Gamma$ denote the set of minimizing geodesics on $X$.
It is a closed subspace of $\Lip_{\diam(X)}([0,1],X)$, defined
by the equation $L(c) \: = \: d(c(0), c(1))$.

For any $t\in [0,1]$, the {\em evaluation map} 
$e_t \: : \: \Gamma \rightarrow X$ defined by 
\begin{equation} e_t(\gamma) \: = \: \gamma(t) \end{equation} 
is continuous. Let $E \: : \: \Gamma \rightarrow X \times X$ be the
``endpoints'' map given by $E(\gamma) \: = \: (e_0(\gamma),e_1(\gamma))$. 
A {\em dynamical transference plan} consists of a transference plan $\pi$ 
and a Borel measure $\Pi$ on $\Gamma$ such that $E_* \Pi \: = \: \pi$;
it is said to be optimal if $\pi$ itself is. 
In words, the transference plan $\pi$ tells us how much mass goes from
a point $x_0$ to another point $x_1$, but does not tell us about the actual
path that the mass has to follow. Intuitively, mass should flow
along geodesics, but there may be several possible choices of geodesics
between two given points and the transport may be divided among
these geodesics; this is the information provided
by $\Pi$. 

If $\Pi$ is an optimal dynamical transference plan then for $t \in [0,1]$, 
we put
\begin{equation} 
\mu_t \: = \: (e_t)_* \Pi.
\end{equation} 
The one-parameter family of measures $\{\mu_t\}_{t \in [0,1]}$ 
is called a {\em displacement interpolation}.
In words, $\mu_t$ is what has become of the mass of $\mu_0$ after
it has travelled from time~0 to time~$t$ according to the dynamical
transference plan $\Pi$.

\begin{lemma} \label{Wgeodesic}
The map $c \: : \: [0,1] \rightarrow P_2(X)$ given by
$c(t) \: = \: \mu_t$ has length $L(c) \: = \: W_2(\mu_0, \mu_1)$.
\end{lemma}
\begin{proof}
Given $0 \: \le \: t \: \le \: t^\prime \: \le \: 1$, 
$(e_t, e_{t^\prime})_* \Pi$ is a particular transference plan
from $\mu_t$ to $\mu_{t^\prime}$, and so
\begin{align} \label{sum}
W_2( \mu_t, \mu_{t^\prime})^2 \: & \le \:
\int_{X \times X} d(x_0, x_1)^2 \: d
\left( (e_t, e_{t^\prime})_* \Pi \right)(x_0, x_1)
\: = \:
\int_{\Gamma} d(\gamma(t), \gamma(t^\prime))^2 \: d\Pi(\gamma) \\
& = \:
\int_{\Gamma} (t^\prime - t)^2 \: 
L \left( \gamma \right)^2 \: d\Pi(\gamma) \: = \:
(t^\prime - t)^2 \: \int_{\Gamma} 
d(\gamma(0), \gamma(1))^2 \: d\Pi(\gamma) \notag \\
& = \: 
(t^\prime - t)^2 \: \int_{X \times X} 
d(x_0, x_1)^2 \: (dE_*\Pi)(x_0, x_1) \: = \:
(t^\prime - t)^2 \: W_2( \mu_0, \mu_1)^2.
\notag
\end{align}
Equation (\ref{sum}) implies that $L(c) \: \le \: W_2(\mu_0, \mu_1)$,
and so $L(c) \: = \: W_2(\mu_0, \mu_1)$.
\end{proof}

\subsection{The Wasserstein space as a length space}

\begin{proposition} \label{Wlengthspace}
Let $(X,d)$ be a compact length space.
Then any two points $\mu_0, \mu_1 \in P_2(X)$ can be joined by a
displacement interpolation. 
\end{proposition}

\begin{proof}
The endpoints map $E$ is Borel and surjective. Given
$(x_0, x_1) \in X \times X$, 
$E^{-1}(x_0, x_1)$ is compact.  It follows that there is a Borel map
$S \: : \: X \times X \rightarrow \Gamma$
so that $E \circ S \: = \: \Id_{X \times X}$ 
\cite[Corollary A.6]{Zimmer (1984)}. In words, $S$ is a measurable
way to join points by minimizing geodesics.
Given $\mu_0, \mu_1 \in P_2(X)$, let $\pi$ be an optimal transference plan
between $\mu_0$ and $\mu_1$, and put $\Pi \: = \: S_* (\pi)$.
The corresponding displacement interpolation joins $\mu_0$ and $\mu_1$.
\end{proof}

\begin{corollary}
If $X$ is a compact length space then $P_2(X)$
is a compact length space.
\end{corollary}
\begin{proof}
We already know that $P_2(X)$ is compact.
Given $\mu_0, \mu_1 \in P_2(X)$, Proposition \ref{Wlengthspace}
gives a displacement interpolation $c$ from $\mu_0$ to $\mu_1$.
By Lemma \ref{Wgeodesic}, $L(c) \: = \: W_2(\mu_0, \mu_1)$,
so $P_2(X)$ is also a length space. 
\end{proof}

\begin{remark}
The same argument shows that $(P(X), W_p)$ is a compact length space
for all $p \in [1, \infty)$, where $W_p$ is the
Wasserstein distance of order $p$ \cite[Section 7.1.1]{Villani (2003)}.
\end{remark}

\begin{example} Suppose that $X \: = \: A \cup B \cup C$, where
$A$, $B$ and $C$ are subsets of the plane given by
$A \: = \: \{(x_1, 0) \: : \: -2 \le x_1 \le -1\}$, 
$B \: = \: \{(x_1, x_2) \: : \: x_1^2 + x_2^2 = 1\}$ and
$C \: = \: \{(x_1, 0) \: : \: 1 \le x_1 \le 2\}$. Let $\mu_0$ be
the one-dimensional Hausdorff measure of $A$ and let $\mu_1$ be the 
one-dimensional Hausdorff
measure of $C$.  Then there is an uncountable number of
Wasserstein geodesics from $\mu_0$ to $\mu_1$, given by the
whims of a switchman at the point $(-1, 0)$. 
\end{example}

\subsection{Wasserstein geodesics as displacement interpolations}

The next result states that every Wasserstein geodesic arises from
a displacement interpolation.

\begin{proposition} \label{everygeod} Let $(X,d)$ be a compact length space
and let $\{\mu_t\}_{t \in [0,1]}$ be a geodesic path in $P_2(X)$.
Then there exists some optimal dynamical transference plan $\Pi$
such that $\{\mu_t\}_{t \in [0,1]}$ 
is the displacement interpolation associated to $\Pi$.
\end{proposition}
\begin{proof}
Let $\{\mu_t\}_{t \in [0,1]}$ be a Wasserstein geodesic. 
Up to reparametrization,
we can assume that for all $t,t'\in [0,1]$,
\begin{equation}
W_2(\mu_t,\mu_{t'}) = |t-t'|\, W_2(\mu_0,\mu_{1}).
\end{equation}
Let $\pi^{(0)}_{x_0, x_{1/2}}$ be an optimal transference plan from
$\mu_0$ to $\mu_{1/2}$, and let
$\pi^{(1/2)}_{x_{1/2}, x_1}$ be an optimal transference plan from
$\mu_{1/2}$ to $\mu_1$. Consider the measure obtained by
``gluing together'' $\pi^{(0)}_{x_0, x_{1/2}}$ and
$\pi^{(1/2)}_{x_{1/2}, x_1}$:
\begin{equation}
M^{(1)} \: = \: \frac{d\pi^{(0)}_{x_0, x_{1/2}} \: 
d\pi^{(1/2)}_{x_{1/2}, x_1}}{d\mu_{1/2}(x_{1/2})}
\end{equation}
on $X \times X \times X$.

The precise meaning of this expression is just as in
the ``gluing lemma'' stated in~\cite[Lemma~7.6]{Villani (2003)}:
Decompose $\pi^{(0)}$ with respect to the projection $p_1 \:
: \: X \times X \rightarrow X$ on the second factor
as $\pi^{(0)} \: = \: \sigma^{(0)}_{x_{1/2}} \: 
\mu_{1/2}(x_{1/2})$, where for $\mu_{1/2}$-almost all $x_{1/2}$,
$\sigma^{(0)}_{x_{1/2}} \in P(p_1^{-1}(x_{1/2}))$ is a
probability measure on $p_1^{-1}(x_{1/2})$.
Decompose $\pi^{(1/2)}$ with respect to the projection $p_0 \:
: \: X \times X \rightarrow X$ on the first factor
as $\pi^{(1/2)} \: = \: \sigma^{(1/2)}_{x_{1/2}} \: 
\mu_{1/2}(x_{1/2})$, where for $\mu_{1/2}$-almost all $x_{1/2}$,
$\sigma^{(1/2)}_{x_{1/2}} \in P(p_0^{-1}(x_{1/2}))$.
Then, for $F \in C(X \times X \times X)$,
\begin{equation}
\int_{X \times X \times X} F \: dM^{(1)} \: \equiv 
\: \int_X \int_{p_1^{-1}(x_{1/2}) \times p_0^{-1}(x_{1/2})}
F(x_0, x_{1/2}, x_1) \: 
d\sigma^{(0)}_{x_{1/2}}(x_0) \:  d\sigma^{(1/2)}_{x_{1/2}}(x_1)
\: d\mu_{1/2}(x_{1/2}).
\end{equation}

The formula
\begin{equation}
d\pi_{x_0, x_1} \: = \: \int_X M^{(1)}_{x_0,x_{1/2},x_1}
\end{equation}
defines a transference plan from $\mu_0$ to $\mu_1$ with cost
\begin{align} \label{inegeverygeod}
\int_{X \times X} d(x_0, x_1)^2 \: d\pi_{x_0, x_1} \: & \le \:
\int_{X \times X \times X} (d(x_0, x_{1/2})\: + \: d(x_{1/2}, x_1))^2 \: 
\frac{d\pi^{(0)}_{x_0, x_{1/2}} \: 
d\pi^{(1/2)}_{x_{1/2}, x_1}}{d\mu_{1/2}(x_{1/2})} \\
& \le \:
\int_{X \times X \times X} 2 (d(x_0, x_{1/2})^2 \: + \: d(x_{1/2}, x_1)^2) \: 
\frac{d\pi^{(0)}_{x_0, x_{1/2}} \: 
d\pi^{(1/2)}_{x_{1/2}, x_1}}{d\mu_{1/2}(x_{1/2})} \notag \\
&  = \:
2 \left( \int_{X \times X} d(x_0, x_{1/2})^2 \: d\pi^{(0)}_{x_0, x_{1/2}}
\: + \:
\int_{X \times X} d(x_{1/2}, x_1)^2 \: d\pi^{(1/2)}_{x_{1/2}, x_1} \right) 
\notag \\
& = \:
2 \left( W_2(\mu_0, \mu_{\frac12})^2 \: + \: W_2(\mu_{\frac12}, \mu_1)^2
\right) \: = \: W_2(\mu_0, \mu_1)^2. \notag
\end{align}
Thus $\pi$ is an optimal transference plan and we must have equality
everywhere in~\eqref{inegeverygeod}. Let
\begin{equation}
B^{(1)} \: = \: \left\{ (x_0, x_{1/2}, x_1) 
\in X \times X \times X \: : \: d(x_0, x_{1/2}) \: = \:
d(x_{1/2}, x_1) \: = \: \frac12 \: d(x_0, x_1) \right\};
\end{equation}
then $M^{(1)}$ is supported on $B^{(1)}$.
For $t \in \{0, \frac12, 1\}$,
define $e_t \: : \: B^{(1)} \rightarrow X$ by
$e_t(x_0, x_{1/2}, x_1) \: = \: x_t$.
Then $(e_t)_* \: M^{(1)} \: = \: \mu_t$.

We can repeat the same procedure using a decomposition of the interval 
$[0,1]$ into $2^i$ subintervals.
For any $i \: \ge \: 1$, define
\begin{align}
B^{(i)} \: = \: & \left\{ (x_0, x_{2^{-i}}, x_{2 \cdot 2^{-i}}, \ldots, 
x_{1 \: - \: 2^{-i}}, x_1) 
\in X^{2^i + 1} \: : \right. \\
& \: \: \: \: \: \: \left. d(x_0, x_{2^{-i}}) \: = \:
d(x_{2^{-i}}, x_{2 \cdot 2^{-i}}) \: = \: \ldots \: = \:
d(x_{1 \: - \: 2^{-i}}, x_1) \: = \: 2^{-i} \: d(x_0, x_1) \right\} . \notag
\end{align}
For $0 \: \le \: j \le \: 2^i - 1$, choose an optimal transference plan
$\pi^{(j \cdot 2^{-i})}_{x_{j \cdot 2^{-i}}, x_{(j+1) \cdot 2^{-i}}}$ from
$\mu_{j \cdot 2^{-i}}$ to $\mu_{(j+1) \cdot 2^{-i}}$.
Then as before, we obtain a probability measure $M^{(i)}$ on $B^{(i)}$ by
\begin{equation}
M^{(i)}_{x_0, x_{2^{-i}}, \ldots, x_1} \: = \:
\frac{d\pi^{(0)}(x_0, x_{2^{-i}}) \: d\pi^{(2^{-i})}(x_{2^{-i}}, 
x_{2 \cdot 2^{-i}}) \:
\ldots \: d\pi^{(1 - 2^{-i})}(x_{1 - 2^{-i}}, x_1)}{d\mu_{2^{-i}}(x_{2^{-i}})
\: \ldots \: d\mu_{1-2^{-i}}(x_{1-2^{-i}})}.
\end{equation}
The formula
\begin{equation}
d\pi_{x_0, x_1} \: = \: \int_{X^{2^i - 1}} 
M^{(i)}_{x_0, x_{2^{-i}}, \ldots, x_1}
\end{equation}
defines a transference plan from $\mu_0$ to $\mu_1$.
For $t \: = \: j \cdot 2^{-i}$, $0 \: \le \: j \: \le \: 2^i$, define
$e_t \: : \: B^{(i)} \rightarrow X$ by
$e_t(x_0, \ldots, x_1) \: = \: x_t$; then
$(e_t)_* \: M^{(i)} \: = \: \mu_t$.

Let $S$ be as in the proof of Proposition~\ref{Wlengthspace}.
Given  $(x_0, \ldots, x_1) \in B^{(i)}$, define a map
$p_{x_0, \ldots, x_1} \: : \: [0,1] \: \rightarrow X$ as the
concatenation of the paths $S(x_0, x_{2^{-i}})$, $S(x_{2^{-i}},
x_{2 \cdot 2^{-i}})$, $\ldots$, and $S(x_{1-2^{-i}}, x_1)$.
As $p_{x_0, \ldots, x_1}$ is a normalized continuous curve from
$x_0$ to $x_1$ of length $d(x_0, x_1)$, it is a geodesic.
For each $i$, the linear functional on $C(\Gamma)$ given by
\begin{equation}
F \rightarrow \int_{X^{2^i + 1}} F(p_{x_0, \ldots, x_1}) \:
dM^{(i)}_{x_0,\ldots, x_1}
\end{equation}
defines a probability measure $R^{(i)}$ 
on the compact space
$\Gamma$. Let $R^{(\infty)}$ be the  limit of a weak-$*$ convergent
subsequence of
$\left\{  R^{(i)} \right\}_{i=1}^\infty$. It is also a probability measure 
on $\Gamma$.

For any $t \in \frac{\N}{2^{\N}} \cap [0,1]$ and $f \in C(X)$, we have
$\int_{K} (e_t)^* f \: dR^{(i)} \: = \: \int_X f \: d\mu_t$ 
for large $i$. Then 
$\int_{K} (e_t)^* f \: dR^{(\infty)} \: = \: \int_X f \: d\mu_t$ for
all $f \in C(X)$, or
equivalently,  $(e_t)_* \: R^{(\infty)} \: = \: \mu_t$.
But as in the proof of Lemma \ref{Wgeodesic},
$(e_t)_* R^{(\infty)}$ is weak-$*$ continuous in $t$.
It follows that $(e_t)_* \: R^{(\infty)} \: = \: \mu_t$ for all
$t \in [0,1]$.
\end{proof}

\subsection{Optimal transport on Riemannian manifolds} \label{recall}

In the rest of this section we discuss the case when $X$ is a
smooth compact connected Riemannian manifold $M$ with 
Riemannian metric $g$. (The results are also valid if $g$ is only
$C^3$-smooth).

Given $\mu_0, \mu_1 \in P_2(M)$ which are absolutely continuous with
respect to $\dvol_M$, it is known that there is a unique
Wasserstein geodesic $c$ joining $\mu_0$ to $\mu_1$
\cite[Theorem 9]{McCann (2001)}. Furthermore, for each $t \in [0, 1]$,
$c(t)$ is absolutely continuous with respect to $\dvol_M$
\cite[Proposition 5.4]{Cordero-Erausquin-McCann-Schmuckenschlager (2001)}.
Thus it makes sense to talk about the length space $P^{ac}_2(M)$ of
Borel probability measures on $M$ that are absolutely continuous with
respect to the Riemannian density, equipped with
the metric $W_2$.  It is a dense totally convex subset of $P_2(M)$. 
Note that if $M$ is other than a point then $P^{ac}_2(M)$ is an incomplete 
metric space and is neither open nor closed in $P_2(M)$.

An optimal transference plan in $P^{ac}_2(M)$ turns out to be a Monge 
transport; that is, $c(t) \: = \: (F_t)_* \mu_0$ for a family of Monge 
transports $\{F_t\}_{t \in [0,1]}$ of $M$. For each $m \in M$, 
$\{F_t(m)\}_{t \in [0,1]}$ is a minimizing geodesic.
Furthermore, there is a $\frac{d^2}2$-concave function $\phi$ on $M$
so that for almost all $m \in M$,
$F_t(m) \: = \: \exp_m(- \: t \: \nabla \phi(m))$
\cite[Theorem 3.2 and 
Corollary 5.2]{Cordero-Erausquin-McCann-Schmuckenschlager (2001)}.
This function $\phi$, just as any $\frac{d^2}2$-concave function on
a compact Riemannian manifold, is Lipschitz~\cite[Lemma~2]{McCann (2001)} 
and has a Hessian almost everywhere 
\cite[Proposition 3.14]{Cordero-Erausquin-McCann-Schmuckenschlager (2001)}.
If we only want the Wasserstein geodesic to be defined for an interval
$[0, r^{-1}]$ then we can use the same formula with $\phi$ being
$\frac{rd^2}{2}$-concave.

\section{Functionals on the Wasserstein space} \label{secfunct}

This section is devoted to the study of certain functions on the
Wasserstein space $P_2(X)$.
We first define the functional $U_\nu$ on $P_2(X)$.
We then define $\lambda$-displacement convexity of the functional, along with
its variations : (weak) $\lambda$-(a.c.) displacement convexity.
We give relations among these various notions of displacement
convexity. 
We define the $H$-functionals $H_{N,\nu}$.  
Finally, under certain displacement convexity 
assumptions, we prove HWI functional inequalities. 

The notion of $\lambda$-displacement convexity is more conventional than
that of weak $\lambda$-displacement convexity. However, the ``weak''
notion turns out to be more useful when considering measured Gromov--Hausdorff
limits. We will see that the ``weak'' hypothesis is sufficient for
proving functional inequalities.

\subsection{Weak displacement convexity}

All of our results will involve a distinguished reference measure,
which is not {\em a priori} canonically given.
So by ``measured length space'' we will
mean a triple $(X, d, \nu)$, where $(X,d)$ is a compact length space
and $\nu$ is a Borel probability measure on $X$. These assumptions
automatically imply that $\nu$ is a regular measure.

We write 
\begin{equation}
P_2(X,\nu) \: = \: \bigl\{ \mu \in P_2(X) \: : \: 
\supp(\mu) \subset \supp(\nu)\bigr\}.
\end{equation}
We denote by $P^{ac}_2(X,\nu)$ the elements of $P_2(X,\nu)$
that are absolutely continuous with respect to $\nu$.

\begin{definition} \label{deffunctional}
Let $U$ be a continuous 
convex function on $[0, \infty)$ with $U(0) \: = \: 0$.
Given $\mu,\nu \in P_2(X)$,
we define the functional $U_\nu \: : \: P_2(X) \rightarrow
\R \cup \{ \infty\}$ by
\begeq \label{defUnu}
U_\nu(\mu) = 
\int_X U(\rho(x))\,d\nu(x) \ + \ U'(\infty)\,\mu_s(X), 
\endeq
where 
\begin{equation} 
\mu = \rho \nu + \mu_s 
\end{equation}
is the Lebesgue decomposition of $\mu$ with respect to $\nu$ 
into an absolutely continuous part $\rho\nu$ and a singular part $\mu_s$.
\end{definition}

\begin{remark} If $U'(\infty)=\infty$, then finiteness of $U_\nu(\mu)$
implies that $\mu$ is absolutely continuous with respect to $\nu$.
This is not true if $U'(\infty)<\infty$.
\end{remark}

\begin{lemma} \label{bound}
$U_\nu(\mu) \: \ge \: U_\nu(\nu)= U(1)$.
\end{lemma}

\begin{remark}
The lemma says that as a function of $\mu$, $U_\nu$ is minimized at
$\nu$.
If $\mu$ is absolutely continuous with respect to $\nu$ then the lemma
is just Jensen's inequality in the form
\begin{equation}
\int_X U(\rho(x))\,d\nu(x) \geq U\left( \int_X \rho(x)\,d\nu(x)\right).
\end{equation}
The general case could be proven using this particular case together with
an approximation argument such as Theorem~\ref{thmapproxU}.
However, we give a direct proof below.
\end{remark}

\begin{proof}[Proof of Lemma~\ref{bound}]
As $U$ is convex, for any $\alpha \in (0,1)$ we have
\begin{equation}
U(\alpha r \: + \: 1-\alpha) \: \le \: \alpha \: U(r) \: + \: (1-\alpha) \: U(1),
\end{equation}
or
\begin{equation}
U(r) \: - \: U(1) \: \ge \: \frac{1}{\alpha} \: \left[ U(\alpha r \: + \: 1 \: - \: \alpha)
\: - \: U(1) \right].
\end{equation}
Then
\begin{equation} \label{breakup}
\int_X U(\rho) \: d\nu \: - \: U(1) \: \ge \: \int_X 
\frac{U(\alpha \rho \: + \: 1 \: - \: \alpha) \: - \: U(1)}{\alpha 
\rho \: - \: \alpha} \: (\rho \: - \: 1) \: d\nu,
\end{equation}
where we take the integrand of the right-hand-side
to vanish at points $x \in X$ where 
$\rho(x) \: = \: 1$. We break up the
right-hand-side of (\ref{breakup}) according to whether 
$\rho(x) \: \le \: 1$ or $\rho(x) \: > \: 1$. 
From monotone convergence, for
$\rho \: \le \: 1$ we have
\begin{equation}
\lim_{\alpha \rightarrow 0^+} \int_X 
\frac{U(\alpha \rho \: + \: 1 \: - \: \alpha) \: - \: U(1)}{\alpha 
\rho \: - \: \alpha} \: (\rho \: - \: 1) \: 
1_{\rho \le 1} \: d\nu  \: = \:
U^\prime_-(1) \: \int_X (\rho \: - \: 1) \: 1_{\rho \le 1} \: d\nu,
\end{equation}
while for $\rho \: > \: 1$ we have
\begin{equation}
\lim_{\alpha \rightarrow 0^+} \int_X 
\frac{U(\alpha \rho \: + \: 1 \: - \: \alpha) \: - \: U(1)}{\alpha 
\rho \: - \: \alpha} \: (\rho \: - \: 1) \: 
1_{\rho > 1} \: d\nu  \: = \:
U^\prime_+(1) \: \int_X (\rho \: - \: 1) \: 1_{\rho > 1} \: d\nu.
\end{equation}
Then
\begin{align}
\int_X U(\rho) \: d\nu \: - \: U(1) \: & \ge \: U^\prime_-(1) \: 
\int_X (\rho \: - \: 1) \: d\nu \: + \:
(U^\prime_+(1) \: - \: U^\prime_-(1)) \: \int_X (\rho \: - \: 1) \: 1_{\rho > 1} \: d\nu \\
& \ge \: U^\prime_-(1) \: 
\int_X (\rho \: - \: 1) \: d\nu \: \ge \: U^\prime(\infty) \: 
\int_X (\rho \: - \: 1) \: d\nu \: = \: - \: U^\prime(\infty) \: \mu_s(X). \notag
\end{align}
As $U_\nu(\nu) \: = \: U(1)$, the lemma follows.
\end{proof}

\begin{definition} \label{dcdef} 
Given a compact measured length space $(X,d,\nu)$ and 
a number $\lambda\in\R$, we say that $U_\nu$ is

\begin{itemize}
\item {\em $\lambda$-displacement convex} if for
all Wasserstein geodesics $\{\mu_t\}_{t \in [0,1]}$ with
$\mu_0, \mu_1 \in P_2(X,\nu)$, we have
\begin{equation} \label{defldc}
U_\nu(\mu_t) \: \le \: t \: U_\nu(\mu_1) \: + \: (1-t) \: 
U_\nu(\mu_0) \: - \:
\frac12 \: \lambda \: t(1-t) W_2(\mu_0, \mu_1)^2
\end{equation}
for all $t \in [0,1]$;
\smallskip

\item {\em weakly $\lambda$-displacement convex} if 
for all $\mu_0, \mu_1 \in P_2(X,\nu)$,
there is {\em some} Wasserstein geodesic
from $\mu_0$ to $\mu_1$ along which~\eqref{defldc} is satisfied;
\smallskip

\item {\em (weakly) $\lambda$-a.c. displacement convex}
if the condition is satisfied when we just assume that
$\mu_0,\mu_1 \in P^{ac}_2(X, \nu)$.
\end{itemize}
\end{definition}

\begin{remark}
In Definition \ref{dcdef} we assume that
$\supp(\mu_0) \subset \supp(\nu)$ and $\supp(\mu_1) \subset \supp(\nu)$,
but we do not assume that $\supp(\mu_t) \subset \supp(\nu)$ for
$t \in (0,1)$.
\end{remark}

\begin{remark}
If $U_\nu$ is $\lambda$-displacement convex and $\supp(\nu) \: = \: X$ then
the function
$t \rightarrow U_\nu(\mu_t)$ is $\lambda$-convex on $[0,1]$, i.e. 
for all $0 \: \le \: s \: \le \: s^\prime \: \le 1$ and
$t \in [0,1]$,
\begin{equation} \label{convexity}
U_\nu(\mu_{t s^\prime + (1-t) s}) \: \le \: t \: U_\nu(\mu_{s^\prime}) 
\: + \: (1-t) \: 
U_\nu(\mu_s) \: - \:
\frac12 \: \lambda \: t(1-t) (s^\prime - s)^2 \: W_2(\mu_0, \mu_1)^2.
\end{equation}
This is not {\em a priori} the case if we only assume that
$U_\nu$ is weakly $\lambda$-displacement convex.
\end{remark}

We may sometimes write ``displacement convex'' instead of 
$0$-displacement convex. 
In short, {\em weakly} means that we require a condition to hold only for
{\em some} geodesic between two measures, as opposed to all geodesics,
and {\em a.c.} means that we only require the condition to hold when the 
two measures are absolutely continuous.

There are obvious implications
\begin{equation} \label{implications}
\begin{matrix}
\text{$\lambda$-displacement convex} & \Longrightarrow &
\text{weakly $\lambda$-displacement convex} \\
\Downarrow & & \Downarrow \\
\text{$\lambda$-a.c. displacement convex} & \Longrightarrow &
\text{weakly $\lambda$-a.c. displacement convex.}
\end{matrix}
\end{equation}

The next proposition reverses the right 
vertical implication in~\eqref{implications}.

\begin{proposition} \label{lemsuffdc}
Let $U$ be a continuous  convex function on $[0, \infty)$
with $U(0) \: = \: 0$. 
Let  $(X,d,\nu)$
be a compact measured length space. 
Then $U_\nu$ is weakly $\lambda$-displacement convex
if and only if it is weakly $\lambda$-a.c. displacement convex.
\end{proposition}
\begin{proof}
We must show that if 
$U_\nu$ is weakly $\lambda$-a.c. displacement convex then it is
weakly $\lambda$-displacement convex. That is, for $\mu_0, \mu_1 \in
P_2(X,\nu)$, we must show that there is some Wasserstein geodesic 
$\{\mu_t\}_{t \in [0,1]}$ from $\mu_0$ to $\mu_1$ along which
\begin{equation} \label{toprovelambdaac}
U_\nu(\mu_t) \: \le \: t \: U_\nu(\mu_0) \: + \: (1-t) \:
U_\nu(\mu_1) \: - \: \frac{1}{2} \: \lambda \: t \: (1-t) \:
W_2(\mu_0, \mu_1)^2. 
\end{equation} 
We may assume that 
$U_\nu(\mu_0) < \infty$ and
$U_\nu(\mu_1) < \infty$, as otherwise~\eqref{toprovelambdaac} is
trivially true for any Wasserstein geodesic from $\mu_0$ to $\mu_1$.
{From} Theorem~\ref{thmapproxU}
in Appendix~\ref{secmoll}, there are
sequences $\{\mu_{k, 0}\}_{k=1}^\infty$ and $\{\mu_{k, 1}\}_{k=1}^\infty$
in $P^{ac}_2(X, \nu)$ (in fact with continuous densities) so that 
$\lim_{k \rightarrow \infty} \mu_{k, 0} \: = \: \mu_0$,
$\lim_{k \rightarrow \infty} \mu_{k, 1} \: = \: \mu_1$, 
$\lim_{k \rightarrow \infty} U_\nu(\mu_{k,0}) \: = \: U_\nu(\mu_{0})$ and
$\lim_{k \rightarrow \infty} U_\nu(\mu_{k,1}) \: = \: U_\nu(\mu_{1})$.
Let $c_k \: : \: [0,1] \rightarrow P_2(X)$ be a minimal
geodesic from $\mu_{k, 0}$ to $\mu_{k, 1}$ such that for
all $t \in [0,1]$,
\begin{equation}
U_\nu(c_k(t)) \: \le \: t \: 
U_\nu(\mu_{k,1}) \: + \:
(1- t) \: U_\nu(\mu_{k,0}) \: - \: \frac12 \: \lambda
t (1- t) W_2(\mu_{k,0}, \mu_{k,1})^2.
\end{equation}
After taking a subsequence,
we may assume that the geodesics $\{c_k\}_{k=1}^\infty$ 
converge uniformly (i.e. in $C([0,1], P_2(X))$) to
a geodesic $c \: : \: [0,1] \rightarrow P_2(X)$ from
$\mu_0$ to $\mu_1$
\cite[Theorem 2.5.14 and Proposition 2.5.17]{Burago-Burago-Ivanov (2001)}. 
The lower semicontinuity of $U_\nu$, Theorem~\ref{thmlsc}(i) in
Appendix~\ref{applsc}, implies that 
$U_\nu(c(t)) \: \le \: \liminf_{k \rightarrow \infty} 
U_\nu(c_k(t))$. The proposition follows.
\end{proof}

In fact, the proof of Proposition \ref{lemsuffdc} gives the following slightly
stronger result.

\begin{lemma} \label{suffcontdens}
Let $U$ be a continuous convex function on $[0, \infty)$
with $U(0) \: = \: 0$.
Let $(X,d,\nu)$ be
a compact measured length space. Suppose that for all
$\mu_0, \mu_1 \in P^{ac}_2(X, \nu)$ with continuous
densities, there is some Wasserstein geodesic from
$\mu_0$ to $\mu_1$ along which~\eqref{defldc} is satisfied.
Then $U_\nu$ is weakly $\lambda$-displacement convex.
\end{lemma}

The next lemma gives sufficient conditions for the horizontal
implications in~\eqref{implications} to be reversed. We recall
the definition of total convexity from Section \ref{lengthspaces}.

\begin{lemma} \label{reversehoriz}
(i) Suppose that $X$ has the property that for
each minimizing geodesic $c \: : \: [0, 1] \rightarrow P_2(X)$,
there is some $\delta_c > 0$ so
that the minimizing geodesic between $c(t)$ and $c(t^\prime)$ is unique 
whenever $|t - t^\prime| \: \le \: \delta_c$. Suppose that
$\supp(\nu) \: = \: X$. If
$U_\nu$ is weakly $\lambda$-displacement convex
then it is $\lambda$-displacement convex.
\smallskip

(ii) Suppose that $P^{ac}_2(X, \nu)$ is totally convex in $P_2(X)$.
Suppose that $X$ has the property that for
each minimizing geodesic $c \: : \: [0, 1] \rightarrow P^{ac}_2(X, \nu)$,
there is some $\delta_c > 0$ so
that the minimizing geodesic between $c(t)$ and $c(t^\prime)$ is unique 
whenever $|t - t^\prime| \: \le \: \delta_c$. Suppose that
$\supp(\nu) \: = \: X$.
If $U_\nu$ is weakly $\lambda$-a.c.
displacement convex
then it is $\lambda$-a.c. displacement convex.
\end{lemma}

\begin{proof}
For part (i), 
suppose that $U_\nu$ is weakly $\lambda$-displacement convex.
Let $c \: : \: [0, 1] \rightarrow P_2(X)$ be a minimizing geodesic.  
We want to 
show the $\lambda$-convexity of $U_\nu$ along $c$.
By assumption, for all $0 \: \le \: s \: \le \: s^\prime \: \le 1$ 
there is some geodesic from
$c(s)$ to $c(s^\prime)$ so that (\ref{convexity}) is satisfied
for all $t \in [0,1]$.
If $|s - s^\prime| \: \le \: \delta_c$ then this geodesic must be
$c \big|_{[s, s^\prime]}$. It follows that the function
$s \rightarrow U_\nu(c(s))$ is $\lambda$-convex on each
interval $[s, s^\prime]$ with $|s - s^\prime| \: \le \: \delta_c$,
and hence on $[0,1]$. This proves part (i).

The same argument works for (ii) provided that
we restrict to absolutely continuous measures.
\end{proof}

\subsection{Important examples}

The following functionals will play an important role.

\begin{definition}
Put
\begin{equation} \label{Udef}
U_N(r) \: = \: 
\begin{cases}
Nr (1 - r^{-1/N}) & \text{ if $1 < N < \infty$}, \\
r \log{r} & \text{ if $N = \infty$}.
\end{cases}
\end{equation}
\end{definition}

\begin{definition}
Let $H_{N,\nu} \: : \: P_2(X) \rightarrow [0, \infty]$ 
be the functional associated to $U_N$, via
Definition~\ref{defUnu}. More explicitly:

\quad - For $N \in (1, \infty)$,
 \begin{equation} 
H_{N,\nu} = N - N \int_X \rho^{1-\frac{1}{N}} \,d\nu, 
\end{equation} 
where $\rho \nu$ is the absolutely continuous part in the Lebesgue
decomposition of $\mu$ with respect to $\nu$.

\quad - For $N=\infty$, the functional $H_{\infty,\nu}$ is defined as 
follows: if $\mu$ is absolutely continuous with respect to $\nu$, with
$\mu \: = \: \rho \nu$, then
\begin{equation} 
H_{\infty,\nu}(\mu) = \int_X \rho\log\rho \,d\nu,
\end{equation}
while if $\mu$ is not absolutely continuous with respect to $\nu$ then 
$H_{\infty,\nu}(\mu)=\infty$.
\end{definition}
\smallskip

To verify that $H_{N,\nu}$ is indeed the functional associated to
$U_N$, we note that $U_N^\prime(\infty) \: = \: N$ and write
\begin{align}
N \int_X \rho \left( 1 \: - \: \rho^{- \: \frac{1}{N}} \right) \: d\nu \: + N \:
\mu_s(X) \: = \: &
N \int_X \rho \left( 1 \: - \: \rho^{- \: \frac{1}{N}} \right) \: d\nu \: + \\
& N \:
\left( 1 \: - \: \int_X \rho \: d\nu \right)  \notag \\
= \: & N \: - \: N \int_X \rho^{1 - \frac{1}{N}} \: d\nu. \notag
\end{align}

Of course, the difference of treatment of the singular part of $\nu$
according to whether $N$ is finite or not reflects the fact
that $U_N$ grows at most linearly when $N<\infty$, but superlinearly when
$N=\infty$. Theorem~\ref{thmlsc} in Appendix \ref{applsc} ensures
that $H_{N,\nu}$ is lower semicontinuous on $P_2(X)$.

\begin{remark}  
Formally extending (\ref{Udef}) to the case $N = 1$ would give
$U_1(r) \: = \: r-1$, which does not satisfy the condition $U(0) = 0$.
This could be ameliorated by instead considering the function
$U(r) \: = \: r$. However, the corresponding entropy functional
$U_\nu$ is identically one, which is not of much use.  We will
deal with the case $N=1$ separately.
\end{remark}

\begin{remark} The quantity $H_{\infty, \nu}(\mu)$ is variously called the
Boltzmann $H$-functional, the negative entropy or
the relative Kullback information of $\mu$ with respect to $\nu$.
As a function of $\mu$, $H_{N, \nu}(\mu)$ attains a minimum when $\mu = \nu$,
which can be considered to be the measure with the least information
content with respect to $\nu$. In some sense, $H_{N, \nu}(\mu)$
is a way of measuring the nonuniformity of $\mu$ with respect to $\nu$.
\end{remark}

\subsection{HWI inequalities} \label{HWI inequalities}

\begin{definition} \label{genfisher}
Let $(X,d,\nu)$ be a compact measured length space.
Let $U$ be a continuous
convex function on $[0, \infty)$, with $U(0)=0$, which is $C^2$-regular on 
$(0,\infty)$. Given $\mu \in P^{ac}_2(X, \nu)$ with
$\rho \: = \: \frac{d\mu}{d\nu}$ a positive Lipschitz function on $X$,
define the ``generalized Fisher information'' $I_U$ by
\begin{equation} 
I_U(\mu) = \int_X U''(\rho)^2\, |\nabla^- \rho|^2\,d\mu =
\int_X \rho \, U''(\rho)^2 \, |\nabla^- \rho|^2\,d\nu.
\end{equation} 
\end{definition}
\noindent
(See Remark~\ref{rqfisher} about the terminology.)

The following estimates generalize the ones that underlie
the HWI inequalities in~\cite{Otto-Villani (2000)}.

\begin{proposition} \label{displconvineq}
Let $(X,d,\nu)$ be a compact measured length space.
Let $U$ be a continuous 
convex function on $[0, \infty)$ with $U(0)=0$.
Given $\mu \in P_2(X,\nu)$, let 
$\{\mu_t\}_{t \in [0,1]}$ be a Wasserstein geodesic from
$\mu_0=\mu$ to $\mu_1=\nu$.
Given $\lambda \in \R$, suppose that 
(\ref{defldc}) is satisfied. Then
\begin{equation} \label{gental}
\frac{\lambda}{2} \: W_2(\mu,\nu)^2 \: \leq \:
U_\nu(\mu) \: - \: U_\nu(\nu).
\end{equation}

Now suppose in addition that $U$ is $C^2$-regular on $(0, \infty)$ and
that $\mu \in P^{ac}_2(X, \nu)$ is such that
$\rho \: = \: \frac{d\mu}{d\nu}$ is a positive Lipschitz function
on $X$. 
Suppose that $U_\nu(\mu) < \infty$
and $\mu_t \in P^{ac}_2(X, \nu)$ for all $t \in [0,1]$.
Then
\begin{equation} \label{dci}
U_\nu(\mu) \: - \: U_\nu(\nu) \: \leq  W_2(\mu,\nu) \sqrt{I_U(\mu)} - 
\frac{\lambda}{2} W_2(\mu,\nu)^2. 
\end{equation}
\end{proposition}
\begin{proof} 
Consider the function $\phi(t) \: = \: U_\nu(\mu_t)$.
Then $\phi(0)=U_\nu(\mu)$ and
$\phi(1) \: = \: U_\nu(\nu)$.
By assumption,
\begin{equation} \label{phiassump}
\phi(t) \: \le \:  t \: \phi(1) \: + \: (1-t) \: \phi(0) \: - \:
\frac{1}{2} \: \lambda \: t \: (1-t) \: W_2(\mu, \nu)^2.
\end{equation}
If $\phi(0) \: - \: \phi(1) \: < \: \frac{1}{2} \: \lambda \: W_2(\mu, \nu)^2$ then as
$\phi(t) \: - \: \phi(1) \: \leq \:
(1-t) \: \left( \phi(0) \: - \: \phi(1) \: - \: \frac12 \: \lambda \: t \: W_2(\mu, \nu)^2
\right)$, we conclude that $\phi(t) \: - \: \phi(1)$ is negative
for $t$ close to $1$, which contradicts Lemma
\ref{bound}.
Thus $\phi(0) \: - \: \phi(1) \: \ge \: \frac{1}{2} \: \lambda \: W_2(\mu, \nu)^2$,
which proves (\ref{gental}).

To prove (\ref{dci}), put $\rho_t= \frac{d\mu_t}{d\nu}$. Then
$\phi(t) \: = \: \int_X U(\rho_t) \: d\nu$.
{From} (\ref{phiassump}), for $t > 0$ we have
\begin{equation}
\phi(0) \: - \: \phi(1) \: \le \:
- \: \frac{\phi(t)-\phi(0)}{t}
 \: - \: \frac{1}{2} \: \lambda \: (1-t) \: W_2(\mu, \nu)^2.
\end{equation}
To prove the inequality \eqref{dci},
it suffices to prove that
\begin{equation} \liminf_{t \rightarrow 0} \left(
- \: \frac{\phi(t)-\phi(0)}{t} \right) \leq W_2(\mu, \nu) 
\sqrt{I_U(\mu)}.
\end{equation}

The convexity of $U$ implies that
\begin{equation} U(\rho_t) - U(\rho_0)  
\geq U'(\rho_0)(\rho_t-\rho_0).
\end{equation}
Integrating with respect to $\nu$ and dividing by $-t<0$, we infer
\begin{equation} 
- \: \frac{1}{t} \left [ \phi(t)-\phi(0) \right ] \: \leq
\: - \:
\frac1{t} \int_X U'(\rho_0(x)) [d\mu_t(x) - d\mu_0(x)].
\end{equation}

By Proposition~\ref{everygeod}, 
$\mu_t = (e_t)_*\Pi$, where $\Pi$ is a certain probability
measure on the space $\Gamma$ of minimal geodesics in $X$. In particular,
\begin{equation} \label{Urhomu}
- \: \frac1{t} \int_X U'(\rho_0(x)) [d\mu_t(x) - d\mu_0(x)] \: = \: - 
\: \frac1{t} \int_\Gamma [ U'(\rho_0(\gamma(t))) - 
U'(\rho_0(\gamma(0))) ]\, d\Pi(\gamma). 
\end{equation}

 Since 
$U'$ is nondecreasing and $t d(\gamma(0), \gamma(1)) = d(\gamma(0), 
\gamma(t))$, we have
\begin{align} \label{tobound}
& - \: \frac{1}{t} \: \int_\Gamma \bigl[ U'(\rho_0(\gamma(t))) - 
U'(\rho_0(\gamma(0))) \bigr ]\, d\Pi(\gamma) \: \le \\
& - \: \frac{1}{t} \: \int_\Gamma 1_{\rho_0(\gamma(t))\leq
\rho_0(\gamma(0))} \: \bigl[ U'(\rho_0(\gamma(t))) - 
U'(\rho_0(\gamma(0))) \bigr]\, d\Pi(\gamma) \: = \notag \\
& \int_\Gamma
\frac{U'(\rho_0(\gamma(t)))-
U'(\rho_0(\gamma(0)))}
{\rho_0(\gamma(t)) - \rho_0(\gamma(0))} \:
\frac{[\rho_0(\gamma(t)) - \rho_0(\gamma(0))]_-}{d(\gamma(0),
\gamma(t))} \: d(\gamma(0),
\gamma(1)) \: d\Pi(\gamma), \notag
\end{align}
where strictly speaking we define the integrand of the last
term to be zero when $\rho_0(\gamma(t)) = \rho_0(\gamma(0))$.
Applying the Cauchy--Schwarz inequality, we can bound the last
term above by
\begin{equation} 
\sqrt{
\int_\Gamma
\frac{[U'(\rho_0(\gamma(t)))-
U'(\rho_0(\gamma(0)))]^2}
{[\rho_0(\gamma(t)) - \rho_0(\gamma(0))]^2} \:
\frac{[\rho_0(\gamma(t)) - \rho_0(\gamma(0))]_-^2}{d(\gamma(0),
\gamma(t))^2} \: d\Pi(\gamma)
}
\; \sqrt{\int_\Gamma d(\gamma(0),\gamma(1))^2 \:
d\Pi(\gamma)}.
\end{equation}
The second square root is just $W_2(\mu_0,\mu_1)$. To conclude the argument,
it suffices to show that
\begin{equation} \label{lll}
\liminf_{t\to 0} \int_\Gamma
\frac{[U'(\rho_0(\gamma(t)))-
U'(\rho_0(\gamma(0)))]^2}
{[\rho_0(\gamma(t)) - \rho_0(\gamma(0))]^2} \:
\frac{[\rho_0(\gamma(t)) - \rho_0(\gamma(0))]_-^2}{d(\gamma(0),
\gamma(t))^2} \: d\Pi(\gamma)
 \leq {I_U}(\mu).
\end{equation}

The continuity of $\rho_0$ implies that 
$\lim_{t \rightarrow 0} \rho_0(\gamma(t)) \: = \: \rho_0(\gamma(0))$. So 
\begin{equation}
\lim_{t \rightarrow 0} \:
 \frac{[U'(\rho_0(\gamma(t)))-U'(\rho_0(\gamma(0)))]^2}
{[\rho_0(\gamma(t)) - \rho_0(\gamma(0))]^2}
\: = \: U''(\rho_0(\gamma(0)))^2.
\end{equation}
On the other hand, the definition of the gradient implies
\begin{equation}
\limsup_{t\to 0} 
\frac{[\rho_0(\gamma(t)) - \rho_0(\gamma(0))]_-^2}{d(\gamma(0),\gamma(t))^2}
\leq |\nabla^- \rho_0|^2 (\gamma(0)).
\end{equation}
As $\rho_0$ is a positive Lipschitz function on $X$,
and $U^\prime$ is $C^1$-regular on $(0,\infty)$, 
$U^\prime \circ \rho_0$ is also Lipschitz
on $X$. Then $\frac{[ U'(\rho_0(\gamma(t))) - 
U'(\rho_0(\gamma(0))) ]^2}
{d(\gamma(0), \gamma(t))^2}$ is uniformly bounded on $\Gamma$,
with respect to $t$, and
dominated convergence implies that
\begin{align} 
& \liminf_{t\to 0} \int_\Gamma
\frac{[U'(\rho_0(\gamma(t)))-
U'(\rho_0(\gamma(0)))]^2}
{[\rho_0(\gamma(t)) - \rho_0(\gamma(0))]^2} \:
\frac{[\rho_0(\gamma(t)) - \rho_0(\gamma(0))]_-^2}{d(\gamma(0),
\gamma(t))^2} \: d\Pi(\gamma) \: \le \\
& \int_\Gamma U''(\rho_0(\gamma(0)))^2 |\nabla^- \rho_0|^2 
(\gamma(0))\, d\Pi(\gamma)
\: = \:  \int_X U''(\rho_0(x))^2 |\nabla^- \rho_0|^2 (x)\, d\mu(x). 
\notag
\end{align}
This concludes the proof of the inequality on the right-hand-side
of~\eqref{dci}. 
\end{proof}

\begin{remark}
Modulo the notational burden caused by the nonsmooth setting,
the proof of Proposition \ref{displconvineq}
is somewhat simpler than the ``standard'' Euclidean proof
because we used a convexity inequality to avoid computing $\phi'(0)$ 
explicitly (compare with~\cite[p.~161]{Villani (2003)}).
\end{remark}

\begin{partics} \label{partcases} 
Taking $U \: = \: U_N$, with $\mu \: = \: \rho \: \nu$ and
$\rho \in \Lip(X)$ a positive function, define
 \begin{equation} \label{fisherdef}
I_{N,\nu} (\mu) = 
\begin{cases} \dps
\left( \frac{N-1}{N} \right)^2 \:
\int_X \frac{|\nabla^- \rho|^2}{\rho^{\frac{2}{N}+1}} \,d\nu &
\text{if $1 < N < \infty$}, \\ \\
\dps \int_X \frac{|\nabla^- \rho|^2}{\rho}\,d\nu &
\text{if $N = \infty$}. 
\end{cases}
\end{equation} 
\smallskip
Proposition (\ref{displconvineq}) implies the following inequalities :

\qquad - If $\lambda > 0$ then
\begin{equation} 
\frac{\lambda}{2} \: W_2(\mu, \nu)^2 \: \leq \:
H_{N,\nu}(\mu) \: \leq \: W_2(\mu, \nu) \: \sqrt{I_{N,\nu}(\mu)} \: - \:  
\frac{\lambda}{2} \: W_2(\mu, \nu)^2 
\: \leq \: \frac{1}{2\lambda} I_{N,\nu}(\mu).
\end{equation}

\qquad - If $\lambda \le 0$ then
\begin{equation} 
H_{N,\nu}(\mu) \: \leq  \: \diam(X) \: \sqrt{I_{N,\nu}(\mu)} \: - \:  
\frac{\lambda}{2} \: \diam(X)^2. 
\end{equation}
\end{partics}

\begin{remark} \label{rqfisher} $I_{\infty,\nu}(\mu)$ is the
classical Fisher information of $\mu$ relative to the reference
measure $\nu$, which is why we call $I_U$ a ``generalized Fisher information''.
\end{remark}

\section{Weak displacement convexity and
measured Gromov--Hausdorff limits} \label{secGH}

In this section we first show that if a sequence of compact metric
spaces converges in the Gromov--Hausdorff topology then their
associated Wasserstein spaces also converge in the Gromov--Hausdorff
topology.  Assuming the results of Appendices \ref{applsc} and \ref{secmoll},
we show that weak displacement convexity of $U_\nu$ is preserved by
measured Gromov--Hausdorff limits. Finally, we define the notion of
weak $\lambda$-displacement convexity for a family ${\cal F}$ of
functions $U$.

\subsection{Gromov--Hausdorff convergence of the Wasserstein space}

\begin{proposition} \label{propGHW}
If $f \: : \: (X_1,d_1) \rightarrow (X_2,d_2)$ is an 
$\epsilon$-Gromov--Hausdorff approximation then 
$f_* \: : \: P_2(X_1) \rightarrow P_2(X_2)$ is an
$\widetilde{\epsilon}$-Gromov--Hausdorff approximation, where
\begin{equation}
\widetilde{\epsilon} \: = \: 4 \epsilon \: + \: 
\sqrt{\epsilon \: (2 \diam(X_2) \: + \: \epsilon)}.
\end{equation}
\end{proposition}

\begin{corollary} \label{GHcor}
If a sequence of compact metric spaces $\{(X_i,d_i)\}_{i=1}^\infty$ 
converges in the
Gromov--Hausdorff topology to a compact metric space $(X,d)$ then
$\{P_2(X_i)\}_{i=1}^\infty$ converges in the Gromov--Hausdorff 
topology to $P_2(X)$. 
\end{corollary}

\begin{proof}[Proof of Proposition~\ref{propGHW}]
Given $\mu_1, \mu_1^\prime \in P_2(X_1)$, let $\pi_1$ be an optimal
transference plan for $\mu_1$ and $\mu_1^\prime$. Put
$\pi_2 \: = \: (f \times f)_* \pi_1$. Then $\pi_2$ is a transference
plan for $f_* \mu_1$ and $f_* \mu_1^\prime$.
We have
\begin{equation}
W_2(f_* \mu_1, f_* \mu_1^\prime)^2 \: \le \:
\int_{X_2 \times X_2} d_{2}(x_2, y_2)^2 \: d\pi_2(x_2, y_2) \: = \:
\int_{X_1 \times X_1} d_{2}(f(x_1), f(y_1))^2 \: d\pi_1(x_1, y_1).
\end{equation}
As
\begin{align}
|d_{2}(f(x_1), f(y_1))^2 \: - \: d_{1}(x_1, y_1)^2| \: & = \:
|d_{2}(f(x_1), f(y_1)) \: - \: d_{1}(x_1, y_1)| \: \cdot \\
& \: \: \: \: \: \: \:
\bigl (d_{2}(f(x_1), f(y_1)) \: + \: d_{1}(x_1, y_1)\bigr ) \notag
\end{align}
we have
\begin{equation}
|d_{2}(f(x_1), f(y_1))^2 \: - \: d_{1}(x_1, y_1)^2| \: \le \:
\epsilon \: (2 \diam(X_1) \: + \: \epsilon)
\end{equation}
and
\begin{equation}
|d_{2}(f(x_1), f(y_1))^2 \: - \: d_{1}(x_1, y_1)^2| \: \le \:
\epsilon \: (2 \diam(X_2) \: + \: \epsilon).
\end{equation}
It follows that
\begin{equation} \label{GHW0}
W_2(f_* \mu_1, f_* \mu_1^\prime)^2 \: \le \:
W_2(\mu_1, \mu_1^\prime)^2 \: + \:
\epsilon \: (2 \diam(X_1) \: + \: \epsilon)
\end{equation}
and
\begin{equation} \label{GHW1}
W_2(f_* \mu_1, f_* \mu_1^\prime)^2 \: \le \:
W_2(\mu_1, \mu_1^\prime)^2 \: + \:
\epsilon \: (2 \diam(X_2) \: + \: \epsilon).
\end{equation}
It follows from this last inequality that
\begin{equation} \label{bbb}
W_2(f_* \mu_1, f_* \mu_1^\prime) \: \le \:
W_2(\mu_1, \mu_1^\prime) \: + \:
\sqrt{\epsilon \: (2 \diam(X_2) \: + \: \epsilon)}.
\end{equation}

We now exchange the roles of $X_1$ and $X_2$. We correspondingly apply
(\ref{GHW0}) instead of (\ref{GHW1}), to 
the map $f^\prime$ and the measures $f_* \mu_1$ and $f_* \mu_1^\prime$,
to obtain
\begin{equation} \label{eqW2ff'}
W_2(f^\prime_* (f_* \mu_1), f^\prime_* (f_* \mu_1^\prime)) \: \le \:
W_2(f_* \mu_1, f_* \mu_1^\prime) \: + \:
\sqrt{\epsilon \: (2 \diam(X_2) \: + \: \epsilon)}.
\end{equation}
Since $f^\prime \circ f$ is an admissible Monge transport between
$\mu_1$ and $(f'\circ f)_*\mu_1$, or between $\mu'_1$ and 
$(f'\circ f)_*\mu'_1$, which moves points by a distance at most $2\var$,
we have
\begin{equation}
W_2((f^\prime \circ f)_* \mu_1, \mu_1) \: \le \: 2 \epsilon,\qquad
W_2((f^\prime \circ f)_* \mu_1^\prime, \mu_1^\prime) \: \le \: 2 \epsilon.
\end{equation}
Thus by~\eqref{eqW2ff'} and the triangle inequality,
\begin{equation} \label{ccc}
W_2(\mu_1, \mu_1^\prime) \: \le \:
W_2(f_* \mu_1, f_* \mu_1^\prime) \: + \: 4 \epsilon \: + \: 
\sqrt{\epsilon \: (2 \diam(X_2) \: + \: \epsilon)}.
\end{equation}
Equations (\ref{bbb}) and (\ref{ccc}) show 
that condition (i) of Definition~\ref{defGH} is satisfied.

Finally, given $\mu_2 \in P_2(X_2)$, consider the Monge transport
$f \circ f^\prime$ from
$\mu_2$ to $(f \circ f^\prime)_* \mu_2$. Then
$W_2(\mu_2, f_*( f^\prime_* \mu_2)) \: \le \: \epsilon$.
Thus condition (ii) of Definition~\ref{defGH} is satisfied as well.
\end{proof}

\begin{remark} The map $f_*$ is generally discontinuous.  In fact,
it is continuous if and only if $f$ is continuous.
\end{remark}

\subsection{Stability of weak displacement convexity}

\begin{theorem} \label{thmstabGH}
Let $\{(X_i, d_i, \nu_i)\}_{i=1}^\infty$ be a sequence of
compact measured length spaces so that
$\lim_{i \rightarrow \infty} (X_i, d_i, \nu_i) \: = \: (X, d, \nu_\infty)$
in the measured Gromov--Hausdorff topology.
Let $U$ be a continuous 
convex function on $[0, \infty)$ with $U(0)=0$. Given 
 $\lambda \in \R$,  
suppose that for all $i$, 
$U_{\nu_i}$ is weakly $\lambda$-displacement convex 
for $(X_i, d_i, \nu_i)$. 
Then $U_{\nu_\infty}$ is weakly $\lambda$-displacement convex
for $(X,d,\nu)$.
\end{theorem}
\begin{proof}
By Lemma~\ref{suffcontdens},
it suffices to show that for any
$\mu_0, \mu_1 \in P_2(X)$ with continuous densities with respect to 
$\nu_\infty$, 
there is a Wasserstein geodesic joining them along 
which inequality~\eqref{defldc} holds for $U_{\nu_\infty}$. We may assume that
$U_{\nu_\infty}(\mu_0) < \infty$ and $U_{\nu_\infty}(\mu_1) < \infty$,
as otherwise any Wasserstein geodesic works.

Write $\mu_0 \: = \: \rho_0 \: \nu_\infty$ and
$\mu_1 \: = \: \rho_1 \: \nu_\infty$. Let $f_i \: : \: X_i \rightarrow X$
be an $\epsilon_i$-approximation, with $\lim_{i \rightarrow \infty} 
\epsilon_i \: = \: 0$ and
$\lim_{i \rightarrow \infty} (f_i)_* \nu_i \: = \: \nu_\infty$. If $i$ is
sufficiently large then $\int_X \rho_0 \: d(f_i)_* \nu_i \: > \: 0$ and
$\int_X \rho_1 \: d(f_i)_* \nu_i \: > \: 0$. For such $i$,
put $\mu_{i,0} \: = \:
\frac{(f_i^* \rho_0) \: \nu_i}{\int_X \rho_0 \: d(f_i)_* \nu_i}$ and $\mu_{i,1} \: = \:
\frac{(f_i^* \rho_1) \: \nu_i}{\int_X \rho_1 \: d(f_i)_* \nu_i}$. Then 
$(f_i)_*  \mu_{i,0} \: = \: 
\frac{\rho_0 \: (f_i)_* \nu_i}{\int_X \rho_0 \: d(f_i)_* \nu_i}$ and
$(f_i)_*  \mu_{i,1} \: = \: 
\frac{\rho_1 \: (f_i)_* \nu_i}{\int_X \rho_1 \: d(f_i)_* \nu_i}$.
Now choose geodesics $c_i \: : \: [0, 1] \rightarrow P_2(X_i)$ with
$c_i(0) \: = \: \mu_{i,0}$ and $c_i(1) \: = \: \mu_{i,1}$ so that
for all $t \in [0,1]$, we have
\begin{equation} \label{topasslimitGHW}
U_{\nu_i}(c_i(t)) \: \le \: t \: U_{\nu_i}(\mu_{i,1}) \: + \: 
(1-t) \: U_{\nu_i}(\mu_{i,0}) \: - \:
\frac12 \: \lambda \: t(1-t) W_2(\mu_{i,0}, \mu_{i,1})^2.
\end{equation}
{From} Lemma \ref{Arzela-Ascoli} and Corollary \ref{GHcor},
after passing to a subsequence,
the maps $(f_i)_* \circ c_i \: : \: [0,1] \rightarrow P_2(X)$
converge uniformly to a 
continuous map $c \: : \: [0, 1] \rightarrow P_2(X)$.
As $W_2(c_i(t), c_i(t^\prime)) \: = \: |t - t^\prime| \: 
W_2(\mu_{i,0}, \mu_{i,1})$, it follows that
$W_2(c(t), c(t^\prime)) \: = \: |t - t^\prime| \: 
W_2(\mu_{0}, \mu_{1})$. Thus $c$ is a Wasserstein geodesic.
The problem is to pass to the limit in~\eqref{topasslimitGHW} as $i\to\infty$.

Given $F \in C(X)$, the fact that $\rho_0 \in C(X)$ implies that 
\begin{equation}
\lim_{i \rightarrow \infty}
\int_X F \: d(f_i)_* \mu_{i,0} \: = \:
\lim_{i \rightarrow \infty}
\int_X F \rho_0 \: \frac{d(f_i)_*\nu_i}{\int_X \rho_0 \: d(f_i)_*\nu_i} \: = \: 
\int_X F \rho_0\: d\nu_\infty.
\end{equation}
Thus $\lim_{i \rightarrow \infty}
(f_i)_*  \mu_{i,0} \: = \: \mu_0$. Similarly,
$\lim_{i \rightarrow \infty}
(f_i)_*  \mu_{i,1} \: = \: \mu_1$. It follows from Corollary \ref{GHcor} that 
\begin{equation} \label{limW2}
\lim_{i \rightarrow \infty} W_2 (\mu_{i,0},
\mu_{i,1}) \: = \: W_2(\mu_0, \mu_1).
\end{equation} 

Next,
\begin{equation}
U_{\nu_i}(\mu_{i,0}) \: = \: \int_{X_i} U \left( 
\frac{f_i^* \rho_0}{\int_X \rho_0 \: d(f_i)_* \nu_i} \right) \: d\nu_i
\: = \: \int_X U \left(
\frac{\rho_0}{\int_X \rho_0 \: d(f_i)_* \nu_i} \right) \: d(f_i)_* \nu_i.
\end{equation}
As
\begin{equation}
\lim_{i \rightarrow \infty} U \left(
\frac{\rho_0}{\int_X \rho_0 \: d(f_i)_* \nu_i} \right) \: = \:
U ( \rho_0 )
\end{equation}
uniformly on $X$, it follows that
\begin{equation}
\lim_{i \rightarrow \infty} \int_X U \left(
\frac{\rho_0}{\int_X \rho_0 \: d(f_i)_* \nu_i} \right) \: d(f_i)_* \nu_i
\: = \: \lim_{i \rightarrow \infty} \int_X U
(\rho_0) \: d(f_i)_* \nu_i \: = \: 
\int_X U(\rho_0) \: d\nu_\infty.
\end{equation}
Thus $\lim_{i \rightarrow \infty}
U_{\nu_i}(\mu_{i,0}) \: = \: U_{\nu_\infty}(\mu_0)$. 
Similarly, 
$\lim_{i \rightarrow \infty}
U_{\nu_i}(\mu_{i,1}) \: = \: U_{\nu_\infty}(\mu_1)$. 

It follows from Theorem~\ref{thmlsc}(ii) in Appendix~\ref{applsc} that
\begeq 
U_{(f_i)_*\nu_i}((f_i)_* c_i(t)) \: \le \: U_{\nu_i}(c_i(t)).
\endeq
Then, for any $t \in [0, 1]$, we can combine this with
the lower semicontinuity of 
$(\mu,\nu)\rightarrow U_\nu(\mu)$ (Theorem~\ref{thmlsc}(i)
in Appendix~\ref{applsc}) to obtain
\begin{equation}
U_{ \nu_\infty}(c(t))
\: \le \: \liminf_{i \rightarrow \infty} U_{(f_i)_* \nu_i} ((f_i)_* c_i(t)) \: \le \:
\liminf_{i \rightarrow \infty}U_{ \nu_i}(c_{i} (t)).
\end{equation}
Combining this with (\ref{limW2}) and the preceding results,
we can take $i \rightarrow \infty$ in~\eqref{topasslimitGHW} and find
\begin{equation}
U_{ \nu_\infty}(c(t)) \: \le \: t \: 
U_{ \nu_\infty}(\mu_1) \: + \:
(1- t) \: U_{ \nu_\infty}(\mu_0) \: - \: \frac12 \: \lambda
t (1- t) W_2(\mu_0, \mu_1)^2.
\end{equation}
This concludes the proof.
\end{proof}

\begin{definition} \label{families}
Let ${\cal F}$ be a family of continuous 
convex functions $U$ on $[0, \infty)$ with $U(0)=0$.
Given a function $\lambda \: : \: {\cal F} \rightarrow \R \cup \{- \infty\}$,
we say that a compact measured length space $(X, d, \nu)$ is 
weakly $\lambda$-displacement convex for the family ${\cal F}$ if
for any $\mu_0, \mu_1 \in P_2(X,\nu)$, 
one can find a Wasserstein geodesic $\{\mu_t\}_{t \in [0,1]}$
from $\mu_0$ to $\mu_1$ so that for each $U \in {\cal F}$,
$U_{\nu}$ satisfies 
\begin{equation}
U_{\nu}(\mu_t) \: \le \: t \: U_{\nu}(\mu_1) \: + \: (1-t) \: 
U_{\nu}(\mu_0) \: - \:
\frac12 \: \lambda(U) \: t(1-t) W_2(\mu_0, \mu_1)^2
\end{equation}
for all $t \in [0,1]$.
\end{definition}

There is also an obvious definition of 
``weakly $\lambda$-a.c. displacement convex for the family ${\cal F}$'',
in which one just requires the condition to hold when 
$\mu_0, \mu_1 \in P^{ac}_2(X, \nu)$.
Note that in Definition \ref{families}, 
the same Wasserstein geodesic $\{\mu_t\}_{t \in [0,1]}$ is
supposed to work for all of the functions $U \in {\cal F}$. 
Hence if $(X, d, \nu)$ is weakly $\lambda$-displacement convex for the family 
${\cal F}$ then it is 
weakly $\lambda(U)$-displacement convex for each $U \in {\cal F}$, but
the converse is not {\em a priori} true.

The proof of Theorem \ref{thmstabGH} establishes the
following result.

\begin{theorem} \label{variantGH}
Let $\{(X_i, d_i, \nu_i)\}_{i=1}^\infty$ be a sequence of
compact measured length spaces with
$\lim_{i \rightarrow \infty} (X_i, d_i, \nu_i) \: = \: (X, d, \nu_\infty)$
in the measured Gromov--Hausdorff toology.
Let ${\cal F}$ be a family of continuous 
convex functions $U$ on $[0, \infty)$ with $U(0)=0$.
Given a function $\lambda \: : \: {\cal F} \rightarrow \R \cup \{-\infty\}$,
suppose that each $(X_i, d_i, \nu_i)$ is weakly $\lambda$-displacement
convex for the family ${\cal F}$.
Then $(X, d, \nu_\infty)$ is weakly $\lambda$-displacement
convex for the family ${\cal F}$.
\end{theorem}

For later use, we note that the proof of Proposition \ref{lemsuffdc}
establishes the following result.

\begin{proposition} \label{lemsuffdc2}
Let ${\cal F}$ be a family of continuous 
convex functions $U$ on $[0, \infty)$ with $U(0)=0$.
Given a function $\lambda \: : \: {\cal F} \rightarrow \R \cup \{-\infty\}$,
$(X, d, \nu)$ is weakly $\lambda$-displacement convex for the family
${\cal F}$
if and only if it is weakly $\lambda$-a.c. displacement convex for
the family ${\cal F}$.
\end{proposition}

\section{$N$-Ricci curvature for measured length spaces} \label{secricci}

This section deals with $N$-Ricci curvature and its
basic properties.  We first define certain classes $\DC_N$ of convex functions
$U$. We use these to 
define the notions of a measured length space $(X,d,\nu)$ having
nonnegative $N$-Ricci curvature, or $\infty$-Ricci curvature bounded
below by $K \in \R$. We show that these properties pass to totally convex
subsets of $X$. We prove that the Ricci curvature definitions are
preserved by measured Gromov--Hausdorff limits.  We show that
nonnegative $N$-Ricci curvature for $N < \infty$
implies a Bishop--Gromov-type inequality.
We show that in certain cases, lower Ricci curvature bounds are
preserved upon quotienting by compact group actions.
Finally, we show that under the assumption of nonnegative $N$-Ricci curvature
with $N < \infty$,
any two measures that are absolutely continuous with respect to $\nu$
can be joined by a Wasserstein geodesic all of whose points are
absolutely continuous measures with respect to $\nu$.

\subsection{Displacement convex classes}

We first define a suitable class of convex functions,
introduced by McCann~\cite{McCann (1997)}. Consider a
continuous
convex function $U \: : \: [0, \infty) \rightarrow \R$ with
$U(0)=0$.
We define the nonnegative function
\begin{equation} \label{pdef}
\p(r) = r U_+'(r) - U(r),
\end{equation}
with $\p(0)=0$.
If one thinks of $U$ as defining an internal energy for a continuous
medium then $\p$ can be thought of as a pressure. By analogy,
if $U$ is $C^2$-regular 
on $(0,\infty)$
then we define the ``iterated pressure''
\begin{equation} 
\pp(r) = r \p'(r) - \p(r).
\end{equation}

\begin{definition} For $N\in [1,\infty)$,
we define $\DC_N$ to be the set of all continuous convex 
functions $U$ on $[0, \infty)$, with $U(0)=0$, such that the function
\begin{equation} 
\psi(\lambda) = \lambda^N \: U(\lambda^{-N})
\end{equation}
is convex on $(0, \infty)$.

We further define $\DC_\infty$ 
to be the set of all continuous convex 
functions $U$ on $[0, \infty)$, with $U(0)=0$, such that the function
\begin{equation}
\psi(\lambda) = e^\lambda \: U(e^{-\lambda}) 
\end{equation} 
is convex on $(- \infty, \infty)$.
\end{definition}

We note that the convexity of $U$ implies that $\psi$ is 
nonincreasing in $\lambda$, as $\frac{U(x)}{x}$ is
nondecreasing in $x$. Below are some useful facts about the classes $\DC_N$.

\begin{lemma} \label{growing}
If $N \le N^\prime$ then $\DC_{N^\prime} \subset \DC_{N}$.
\end{lemma}
\begin{proof}
If $N^\prime < \infty$, let $\psi_{N}$ and $\psi_{N^\prime}$ denote the
corresponding functions. Then $\psi_{N}(\lambda) \: = \:
\psi_{N^\prime} \left( \lambda^{N/N^\prime} \right)$. 
The conclusion follows from the
fact that the function $x \rightarrow x^{N/N^\prime}$ is concave on
$[0, \infty)$, along with the fact that the composition of a
nonincreasing convex function and a concave function is convex.
The case $N^\prime = \infty$ is similar.
\end{proof}

\begin{lemma} \label{p2lem}
For $N \in [1, \infty]$, \\
(a) If $U$ is a continuous convex 
function on $[0, \infty)$ with $U(0)=0$ then
$U \in \DC_N$ if and only if the function
$r\longmapsto p(r)/r^{1-\frac{1}{N}}$ 
is nondecreasing on $(0, \infty)$. \\
(b) If $U$ is a continuous convex 
function on $[0, \infty)$ that is $C^2$-regular on
$(0, \infty)$, with $U(0)=0$, then
$U \in \DC_N$ if and only if 
$\pp \: \geq \: - \: \frac{p}{N}$.
\end{lemma}
\begin{proof}
Suppose first that $U$ is a continuous convex
function on $[0, \infty)$ and $N \in [1, \infty)$.
Putting $r(\lambda) \: = \: \lambda^{-N}$, one can check that 
\begin{equation} 
\psi_-'(\lambda) \: = \: - \: N \: p(r)/r^{1- \frac{1}{N}}.
\end{equation}
Then $\psi$ is convex if and only if $\psi_-^\prime$ is nondecreasing, which 
is the case if and only if the function 
$r\longmapsto p(r)/r^{1-\frac{1}{N}}$ is
nondecreasing (since the map $\lambda \rightarrow \lambda^{-N}$ is
nonincreasing).
Next, suppose that $U$ is $C^2$-regular on $(0,\infty)$.
One can check that
\begin{equation} 
\psi''(\lambda) \: = \: N^2 \: r^{\frac{2}{N}-1} \left( p_2(r) + \frac{p(r)}{N}
\right).
\end{equation}
Then $\psi$ is convex if and only if $\psi^{\prime \prime} \ge 0$, which 
is the case if and only if $p_2 \: \ge \: - \: \frac{p}{N}$.

The proof in the case $N = \infty$ is similar.
\end{proof}

\begin{lemma} \label{increasing}
Given $U\in \DC_\infty$, either $U$ is linear or there exist $a, b>0$
such that 
$U(r)\geq a\, r\log r \, - \, br$.
\end{lemma}
\begin{proof}
The function $U$ can be reconstructed from $\psi$ by the formula
\begin{equation}
U(x) \: =  \: x \: \psi(\log(1/x)).
\end{equation}
As $\psi$ is convex and nonincreasing, either $\psi$ is constant or there
are constants $a, b > 0$ such that $\psi(\lambda) \: \ge \: - 
\: a \lambda - b$ for all $\lambda \in \R$. 
In the first case, $U$ is linear. In the second case, 
we have $U(x)\geq -a x\log (1/x) -bx$,
as required.
\end{proof}

\subsection{Ricci curvature via weak displacement convexity}

We recall from Definition \ref{families} the notion of
a compact measured length space $(X,d,\nu)$ 
being weakly $\lambda$-displacement convex for a family of 
convex functions ${\cal F}$.

\begin{definition} \label{nonnegdef}
Given $N\in [1,\infty]$, we say that a compact measured length space
$(X,d,\nu)$ has nonnegative $N$-Ricci curvature
if it is weakly displacement convex for the family
$\DC_N$.
\end{definition}

By Lemma \ref{growing}, if $N \le N^\prime$ and $X$ has nonnegative
$N$-Ricci curvature then it has nonnegative $N^\prime$-Ricci curvature.
In the case $N=\infty$, we can define a more precise notion.

\begin{definition} \label{lambdadef}
Given $K\in\R$, define
$\lambda \: : \: \DC_\infty \rightarrow \R \cup \{-\infty\}$ by
\begeq\label{lambdainfty} 
\lambda(U) = \inf_{r > 0} K \: \frac{p(r)}r = 
\begin{cases}
K \lim_{r \rightarrow 0^+} \frac{p(r)}r & \text{if $K > 0$}, \\
0 & \text{if $K = 0$}, \\
K \lim_{r \rightarrow \infty} \frac{p(r)}r & \text{if $K < 0$},
\end{cases}
\endeq
where $p$ is given by (\ref{pdef}).
We say that a compact measured length space $(X,d,\nu)$
has $\infty$-Ricci curvature bounded
below by $K$
if it is weakly $\lambda$-displacement convex for the family $\DC_\infty$.
\end{definition}

If $K \: \le \: K^\prime$ and
$(X,d,\nu)$
has $\infty$-Ricci curvature bounded
below by $K^\prime$ then it
has $\infty$-Ricci curvature bounded
below by $K$.

The next proposition shows that our definitions localize on
totally convex subsets.

\begin{proposition} \label{HNrestr}
Suppose that a closed set $A \subset X$ is totally convex. 
Given $\nu \in P_2(X)$ with $\nu(A) > 0$, put
$\nu^\prime \: = \: \frac{1}{\nu(A)} \: \nu \big|_A \in P_2(A)$.\\
(a) If
$(X, d, \nu)$ has nonnegative $N$-Ricci curvature then
$(A, d, \nu^\prime)$ has nonnegative $N$-Ricci curvature. \\
(b) If
$(X,d,\nu)$ has $\infty$-Ricci curvature bounded
below by $K$ then $(A, d, \nu^\prime)$ has 
$\infty$-Ricci curvature bounded
below by $K$.
\end{proposition}

\begin{proof}
By Proposition \ref{everygeod}, 
$P_2(A)$ is a totally convex subset of $P_2(X)$.
Given $\mu \in P_2(A) \subset P_2(X)$, let $\mu \: = \: \rho \: \nu \: + \:
\mu_s$ be its Lebesgue decomposition with respect
to $\nu$. Then
$\mu \: = \: \rho^\prime \: \nu^\prime \: + \: \mu_s$ 
is the Lebesgue decomposition
of $\mu$ with respect to $\nu^\prime$, where
$\rho^\prime \: = \: \nu(A) \: \rho \big|_A$.
Given a continuous 
convex function $U \: : \: [0, \infty) \rightarrow \R$ with
$U(0) \: = 0$, define 
\begin{equation}
\widetilde{U}(r) \: = \: \frac{U(\nu(A) r)}{\nu(A)}.
\end{equation} 
Then $\widetilde{U}^\prime(\infty) \: = \: 
{U}^\prime(\infty)$, and $U \in \DC_N$ if and only if $\widetilde{U} \in \DC_N$ .
Now
\begin{align}
U_{\nu^\prime}(\mu) \: & = \: \int_A  U(\rho^\prime) \: d\nu^\prime \: + \:
U^\prime(\infty) \: \mu_s(A) \\
& = \: \frac{1}{\nu(A)} \:
\int_A  U(\nu(A) \rho) \: d\nu \: + \:
U^\prime(\infty) \: \mu_s(A) \notag \\
& = \int_X  \widetilde{U}(\rho) \: d\nu
\: + \: \widetilde{U}^\prime(\infty)  \: \mu_s(X) \: = \: 
\widetilde{U}_{\nu}(\mu).
\notag
\end{align}
As $P_2(A, \nu^\prime) \subset P_2(X, \nu)$, 
part (a) follows.

Letting $\widetilde{p}$ denote the pressure of $\widetilde{U}$, one finds
that
\begin{equation}
\frac{\widetilde{p}(r)}{r} \: = \: 
\frac{{p}(\nu(A)r)}{\nu(A)r}.
\end{equation}
Then with reference to Definition \ref{lambdadef}, 
$\lambda(\widetilde{U}) \: = \:
\lambda(U)$. Part (b) follows.
\end{proof}

\subsection{Preservation of $N$-Ricci curvature bounds}

The next theorem can be considered to be the main result of this paper.

\begin{theorem} \label{stabRicci}
Let $\{(X_i,d_i,\nu_i)\}_{i=1}^\infty$ be a sequence of
compact measured 
length spaces with $\lim_{i \rightarrow \infty} (X_i, d_i, \nu_i) \: = \:
(X,d,\nu)$ in the measured Gromov--Hausdorff topology.

\qquad - If each $(X_i, d_i, \nu_i)$ has nonnegative $N$-Ricci curvature
then $(X, d, \nu)$ has nonnegative $N$-Ricci curvature.

\qquad - If each $(X_i, d_i, \nu_i)$ has $\infty$-Ricci curvature bounded below
by $K$, for some $K\in\R$,
then $(X, d, \nu)$ has $\infty$-Ricci curvature bounded below
by $K$.
\end{theorem}

\begin{proof} 
If $N < \infty$ then the theorem follows from Theorem~\ref{variantGH} with the
family ${\cal F}=\DC_N$ and $\lambda=0$.
If $N = \infty$ then it follows from Theorem~\ref{variantGH} with the
family ${\cal F}=\DC_\infty$ and $\lambda$ given by Definition
\ref{lambdadef}.
\end{proof}

In what we have presented so far, 
the concept of $(X, d, \nu)$ having nonnegative $N$-Ricci curvature,
or having $\infty$-Ricci curvature bounded below by $K$, 
may seem somewhat abstract.
In Section~\ref{Riemcase} we will show that in the setting of
Riemannian manifolds, it can be expressed in terms of classical
tensors related to the Ricci tensor.

\subsection{Bishop--Gromov inequality}

We first show that a weak displacement convexity
assumption implies that the measure $\nu$ either is a delta function
or is nonatomic.

\begin{proposition} \label{nonatom}
Let $(X,d,\nu)$ be a compact measured length space.
For all $N \in (1, \infty]$,
if $H_{N, \nu}$ is weakly $\lambda$-displacement convex then $\nu$ 
either is a delta function or is nonatomic.
\end{proposition}

\begin{proof}
We will assume that $\nu(\{x\}) \in(0,1)$ for some $x \in X$ and derive
a contradiction.

Suppose first that $N \in (1,\infty)$.
Put $\mu_0 \: = \: \delta_x$ and $\mu_1 \: = \:
\frac{\nu \: - \: \nu(\{x\}) \delta_x}{1 - \nu(\{x\})}$. 
By the hypothesis and Proposition \ref{everygeod}, there is a
displacement interpolation $\{ \mu_t \}_{t \in [0,1]}$ from $\mu_0$ to
$\mu_1$ along which (\ref{defldc}) is satisfied with 
$U_\nu \: = \: H_{N, \nu}$ . Now
$H_{N, \nu}(\mu_0) \: = \: N \: - \: N \: (\nu(\{x\}))^{1/N}$ and
$H_{N, \nu}(\mu_1) \: = \: N \: - \: N \: (1-\nu(\{x\}))^{1/N}$. Hence
\begin{equation} \label{contrad1}
H_{N, \nu}(\mu_t) \: \le \: N \: - \: (1-t) \: N \: (\nu(\{x\}))^{1/N}
\: - \: t \: N \: (1-\nu(\{x\}))^{1/N} \: - \: \frac12 \: \lambda \:
t(1-t) \: W_2(\mu_0, \mu_1)^2.
\end{equation}
Put $D=\diam(X)$. As we have a displacement interpolation, it
follows that if $t > 0$ then $\supp(\mu_t) \subset 
\overline{B_{t D}(x)}$ and
$\mu_t(\{x\}) = 0$. 
Letting $\mu_t \: = \: \rho_t \: \nu \: + \: (\mu_t)_s$ be the
Lebesgue decomposition of $\mu_t$ with respect to $\nu$, 
H\"older's inequality implies that
\begin{align}
\int_X \rho_t^{1-\frac{1}{N}} \: d\nu \: & = \:
\int_{\overline{B_{t D}(x)} - \{x\}} \rho_t^{1-\frac{1}{N}} \: d\nu \\
& \le \:
\left( \int_{\overline{B_{t D}(x)} - \{x\}}
\rho_t \: d\nu \right)^{1-\frac{1}{N}} \: 
\nu \left( \overline{B_{t D}(x)} - \{x\} 
\right)^{\frac{1}{N}} \notag \\
& \le \:
\nu \left( \overline{B_{t D}(x)} - \{x\} \right)^{\frac{1}{N}}. \notag
\end{align}
Then
\begin{equation}
H_{N, \nu}(\mu_t) \: \ge \:
N \: - \: N \: \left( \nu \left( \overline{B_{t D}(x)} \right) 
- \nu(\{x\}) \right)^{1/N}.
\end{equation}
As $\lim_{t \rightarrow 0^+} \nu \left( \overline{B_{t D}(x)} \right) 
= \nu(\{x\})$, we obtain a contradiction with (\ref{contrad1}) 
when $t$ is small.
\smallskip

If $N = \infty$ then $H_{\infty, \nu}(\mu_0) \: = \: \log
\frac{1}{\nu(\{x\})}$ and
$H_{\infty, \nu}(\mu_1) \: = \: \log \frac{1}{1-\nu(\{x\})}$. 
Hence
\begin{equation} \label{contrad2}
H_{\infty, \nu}(\mu_t) \: \le \: (1-t) \: \log \frac{1}{\nu(\{x\})}
\: + \: t \: \log \frac{1}{1-\nu(\{x\})} \: - \: \frac12 \: \lambda \:
t(1-t) \: W_2(\mu_0, \mu_1)^2.
\end{equation}
In particular, $\mu_t$ is absolutely
continuous with respect to $\nu$. Write $\mu_t \: = \: \rho_t \: \nu$.
Jensen's inequality implies that for $t > 0$,
\begin{align}
& \int_{\overline{B_{t D}(x)} - \{x\}} \rho_t \: \log(\rho_t) \: 
\frac{d\nu}{\nu \left( \overline{B_{t D}(x)} - \{x\} \right)} \: \ge \\
& \left( \int_{\overline{B_{t D}(x)} - \{x\}} \rho_t \: 
\frac{d\nu}{\nu \left( \overline{B_{t D}(x)} - \{x\} \right)} 
\right) \cdot 
\log \left(
\int_{\overline{B_{t D}(x)} - \{x\}} \rho_t \: 
\frac{d\nu}{\nu \left( \overline{B_{t D}(x)} - \{x\} \right)} 
\right) \: = \notag \\ 
& \frac{1}{\nu \left( \overline{B_{t D}(x)} 
- \{x\} \right)} \: \log \left(
\frac{1}{\nu \left( \overline{B_{t D}(x)} 
- \{x\} \right)} \right). \notag
\end{align}
Then
\begin{align}
H_{\infty, \nu}(\mu_t) \: & = \: 
\int_{X} \rho_t \: \log(\rho_t) \: 
d\nu \: = \:
\int_{\overline{B_{t D}(x)} - \{x\}} \rho_t \: \log(\rho_t) \: 
d\nu \\
& \ge \: \log \left(
\frac{1}{\nu \left( \overline{B_{t D}(x)} - \{x\} \right)} \right).
\notag
\end{align}
As $\lim_{t \rightarrow 0^+} \nu \left( \overline{B_{t D}(x)} \right) 
= \nu(\{x\})$, we obtain a contradiction with (\ref{contrad2}) 
when $t$ is small.
\end{proof}

We now prove a Bishop--Gromov-type inequality.

\begin{proposition} \label{BGprop}
Let $(X,d,\nu)$ be a compact measured length space.
Assume that $H_{N, \nu}$ is weakly displacement convex 
on $P_2(X)$, for some $N \in (1,  \infty)$. Then for all 
$x\in \supp(\nu)$ and all $0<r_1\leq r_2$,
\begin{equation} 
\nu(B_{r_2}(x)) \leq \left ( \frac{r_2}{r_1} \right )^N \nu(B_{r_1}(x)).
\end{equation}
\end{proposition}
\begin{proof}
{From} Proposition \ref{nonatom}
we may assume that $\nu$ is nonatomic, as the theorem is
trivially true when $\nu = \delta_x$.
Put $\mu_0 \: = \: \delta_x$ and
$\mu_1 \: = \: 
\frac{1_{B_{r_2}(x)}}{\nu(B_{r_2}(x))} \: \nu$. By the hypothesis and
Proposition \ref{everygeod}, there is a
displacement interpolation $\{ \mu_t \}_{t \in [0,1]}$ from $\mu_0$ to
$\mu_1$ along which (\ref{defldc}) is satisfied with
$U_\nu \: = \: H_{N, \nu}$ and $\lambda \: = \: 0$.  Now
$H_{N, \nu}(\mu_0) \: = \: N$ and
$H_{N, \nu}(\mu_1) \: = \: N \: - \: N \: (\nu(B_{r_2}(x)))^{1/N}$. Hence
\begin{equation}
H_{N, \nu}(\mu_t) \: \le \: 
N \: - \: t \: N \: (\nu(B_{r_2}(x)))^{1/N}.
\end{equation}
Let $\mu_t \: = \: \rho_t \: \nu \: + \: (\mu_t)_s$ be the
Lebesgue decomposition of $\mu_t$ with respect to $\nu$.
As we have a displacement interpolation, $\rho_t$ vanishes outside of
$B_{tr_2}(x)$. 
Then from H\"older's inequality,
\begin{equation}
H_{N, \nu}(\mu_t) \: \ge \:
N \: - \: N \: (\nu(B_{tr_2}(x)))^{1/N}.
\end{equation}
The theorem follows by taking
$t \: = \: \frac{r_1}{r_2}$.
\end{proof}

\begin{theorem} \label{thmBG}
If a compact measured length space
$(X, d, \nu)$ has nonnegative $N$-Ricci curvature for some
$N \in [1,\infty)$ then for all 
$x\in \supp(\nu)$ and all $0<r_1\leq r_2$,
\begin{equation} \label{nuBB}
\nu(B_{r_2}(x)) \leq \left ( \frac{r_2}{r_1} \right )^N \nu(B_{r_1}(x)).
\end{equation}
\end{theorem}
\begin{proof}
If $N \in (1,\infty)$ then the theorem follows from Proposition
\ref{BGprop}. If $N = 1$ then $(X, d, \nu)$ has nonnegative $N^\prime$-Ricci curvature 
for all $N^\prime \in (1,\infty)$. The theorem now follows by replacing $N$ in
(\ref{nuBB}) by $N^\prime$ and taking $N^\prime \rightarrow 1$.
\end{proof}

\begin{corollary} \label{MGHcompactness}
Given $N \in [1, \infty)$ and $D \: \ge \: 0$, the space of compact measured 
length spaces $(X, d, \nu)$ with nonnegative $N$-Ricci curvature,
$\diam(X, d) \: \le \: D$ and $\supp(\nu) \: = \: X$ is sequentially
compact in the measured Gromov--Hausdorff topology.
\end{corollary}
\begin{proof}
Let $\{(X_i, d_i, \nu_i)\}_{i=1}^\infty$ be a sequence of such spaces.
Using the Bishop--Gromov inequality of Theorem \ref{thmBG}, along with the
fact that $\supp(\nu_i) \: = \: X_i$, it follows as in
\cite[Theorem 5.3]{Gromov (1999)} that after passing to a subsequence we may
assume that $\{(X_i, d_i) \}_{i=1}^\infty$ converges in the
Gromov--Hausdorff topology to a compact length space $(X, d)$, 
necessarily with $\diam(X, d) \: \le \: D$. Let
$f_i \: : \: X_i \rightarrow X$ be Borel $\epsilon_i$-approximations,
with $\lim_{i \rightarrow \infty} \epsilon_i \: = \: 0$.
From the compactness of $P_2(X)$, after passing to a subsequence we may
assume that $\lim_{i \rightarrow \infty} (f_i)_* \nu_i \: = \: \nu$
for some $\nu \in P_2(X)$. From Theorem \ref{stabRicci}, 
$(X, d, \nu)$ has nonnegative $N$-Ricci curvature.

It remains to show that $\supp(\nu) \: = \: X$.
Given $x \in X$, the measured Gromov--Hausdorff convergence of
$\{(X_i, d_i, \nu_i)\}_{i=1}^\infty$ to $(X, d, \nu)$ implies that
there is a sequence of points $x_i \in X_i$ with
$\lim_{i \rightarrow \infty} f_i(x_i) \: = \: x$ so that for all $r > 0$ and
$\epsilon \in (0, r)$, we have
$\limsup_{i \rightarrow \infty} \nu_i(\overline{B_{r - \epsilon}(x_i)})
\: \le \: \nu({B_{r}(x)})$.
By Theorem \ref{thmBG}, $(r-\epsilon)^{-N} \: \nu_i(B_{r-\epsilon}(x_i)) 
\: \ge \:
\diam(X_i, d_i)^{-N}$. Then $\nu({B_{r}(x)}) \: \ge \:
\left( \frac{r}{\diam(X, d)} \right)^{N}$, which proves the claim.
\end{proof}

\begin{remark}
Corollary \ref{MGHcompactness} shows that it is consistent in some
sense to restrict to the case $\supp(\nu) = X$, at least when
$N$ is finite; see also Theorem \ref{thmredspt}. The analog
of Corollary \ref{MGHcompactness} does not hold in the case
$N = \infty$, as can be seen by taking $X \: =\: [-1,1]$,
$d(x,y) \: = \: |x-y|$,
$\nu \: = \: \frac{e^{- t x^2} \: dx}{\int_{-1}^1 e^{- t x^2} \: dx}$
and letting $t$ go to infinity. 
\end{remark}

\subsection{Compact group actions}

In this section we show that in certain cases,
lower Ricci curvature bounds are preserved upon
quotienting by a compact group action.

\begin{theorem} \label{Gaction}
Let $(X, d, \nu)$ be a compact measured length space.
Suppose that any two $\mu_0, \mu_1 \in P_2^{ac}(X, \nu)$ are joined
by a unique Wasserstein geodesic, that lies in
$P_2^{ac}(X, \nu)$. 
Suppose that a compact topological group $G$ acts continuously and
isometrically on $X$,
preserving $\nu$. Let $p \: : \: X \rightarrow X/G$ be the
quotient map and let $d_{X/G}$ be the quotient metric. 
We have the following implications: \\
\quad a. For $N \in [1, \infty)$, if $(X, d, \nu)$ has nonnegative $N$-Ricci
curvature then $(X/G, d_{X/G}, p_* \nu)$ 
has nonnegative $N$-Ricci curvature. \\
\quad b. If $(X, d, \nu)$ has $\infty$-Ricci curvature bounded below by
$K$ then $(X/G, d_{X/G}, 
p_* \nu)$ has $\infty$-Ricci curvature bounded below by
$K$.
\end{theorem}

The proof of this theorem will be an easy consequence of the following
lemma, which does not involve the length space structure.

\begin{lemma} \label{lemaction}
The map $p_*:P_2(X)\rightarrow P_2(X/G)$ restricts to an
isometric isomorphism between the set $P_2(X)^G$ of $G$-invariant elements
in $P_2(X)$, and $P_2(X/G)$.
\end{lemma}

\begin{proof}
Let $dh$ be the normalized Haar measure on $G$.
The map $p_* \: : \: P_2(X) \rightarrow P_2(X/G)$ restricts to an
isomorphism $p_* \: : \: P_2(X)^G \rightarrow P_2(X/G)$; the problem
is to show that it is an isometry.

Let $\widetilde{\pi}$ be a transference plan between $\widetilde{\mu}_0, 
\widetilde{\mu}_1 \in P_2(X)^G$. Then
$\widetilde{\pi}^\prime \: = \:
\int_G g \cdot \widetilde{\pi} \: dh(g)$ is also a transference plan between
$\widetilde{\mu}_0$ and $\widetilde{\mu}_1$, with
\begin{align}
\int_{X \times X} d_X(\widetilde{x}, \widetilde{y})^2 \: 
d \widetilde{\pi}^\prime (\widetilde{x}, \widetilde{y}) \: & = \:
\int_G \int_{X \times X} d_X(\widetilde{x} g, \widetilde{y} g)^2 \: d 
\pi(\widetilde{x}, \widetilde{y}) \: dh(g) \\
&= \:
\int_{X \times X} 
d_X(\widetilde{x}, \widetilde{y})^2 \: 
d\pi(\widetilde{x}, \widetilde{y}). \notag
\end{align}
Thus there is a $G$-invariant optimal transference 
plan $\widetilde{\pi}$ between
$\widetilde{\mu}_0$ and 
$\widetilde{\mu}_1$. As $\pi \: = \:
(p \times p)_* \widetilde{\pi}$ is a transference plan
between $p_* \widetilde{\mu}_0$ and 
$p_* \widetilde{\mu}_1$, with
\begin{align}
\int_{(X/G) \times (X/G)}  d_{X/G}(x, y)^2 \: d\pi(x,y) \: & = \:
\int_{X \times X}  d_{X/G}(p(\widetilde{x}), p(\widetilde{y}))^2 \: 
d\widetilde{\pi}(\widetilde{x}, \widetilde{y}) \\
& \le \: \int_{X \times X}  d_X(\widetilde{x}, \widetilde{y})^2 \: 
d\widetilde{\pi}(\widetilde{x}, \widetilde{y}), \notag
\end{align}
it follows that the map
$p_* \: : \: P_2(X)^G \rightarrow P_2(X/G)$ is metrically nonincreasing.

Conversely, let $s \: : \: (X/G) \times (X/G) \rightarrow X \times X$ 
be a Borel map such that $(p \times p) 
\circ s \: = \: \Id$ and $d_X \circ s \: = \: d_{X/G}$.
That is, given $x, y \in X/G$, the map $s$ picks points
$\widetilde{x} \in p^{-1}(x)$ and $\widetilde{y} \in p^{-1}(y)$ 
in the corresponding orbits so that the distance between $\widetilde{x}$
and $\widetilde{y}$ is minimized among all pairs of points in 
$p^{-1}(x) \times p^{-1}(y)$.
(The existence of $s$ follows from applying
\cite[Corollary A.6]{Zimmer (1984)} to the restriction of
$p \times p$ to 
$\{ (\widetilde{x}, \widetilde{y}) \in X \times X 
\: : \: d_X(\widetilde{x}, \widetilde{y}) \: = \: 
d_{X/G} (p(\widetilde{x}), p(\widetilde{y})) \}$.
The restriction map is a
surjective Borel map with compact preimages.)
Given an optimal transference plan $\pi$ between $\mu_0, \mu_1 \in
P_2(X/G)$, define a measure $\widetilde{\pi}$ on $X \times X$ by
saying that for all $\widetilde{F} \in C(X \times X)$,
\begin{equation}
\int_{X \times X} \widetilde{F} \: d\widetilde{\pi} \: = \: 
\int_G \int_{(X/G) \times (X/G)} \widetilde{F}(s(x, y)\cdot (g, g)) \:
d\pi(x, y) \: dh(g).
\end{equation}
Then for $F \in C((X/G) \times (X/G))$,
\begin{align}
\int_{(X/G) \times (X/G)} F \: d(p \times p)_* \widetilde{\pi} \: & = \: 
\int_{X \times X} (p \times p)^* F \: d\widetilde{\pi} \\
& = \: 
\int_G \int_{(X/G) \times (X/G)} \left( 
(p \times p)^*{F} \right)(s(x, y)\cdot (g, g)) \:
d\pi(x, y) \: dh(g) \notag \\
& \: = \: 
\int_G \int_{(X/G) \times (X/G)} 
{F} \left( (p \times p)(s(x, y)\cdot (g, g)) \right) \:
d\pi(x, y) \: dh(g) \notag \\ 
& = \: 
\int_{(X/G) \times (X/G)} 
{F}(x,y) \: d\pi(x, y). \notag
\end{align}
Thus $(p \times p)_* \widetilde{\pi} \: = \: \pi$. As
$\widetilde{\pi}$ is $G$-invariant, it follows that it is
a transference plan between
$(p_*)^{-1}(\mu_0), (p_*)^{-1}(\mu_1) \in P_2(X)^G$.
Now
\begin{align}
\int_{X \times X} d_X(\widetilde{x}, \widetilde{y})^2 \: 
d\widetilde{\pi}(\widetilde{x}, \widetilde{y}) \: & = \: 
\int_G \int_{(X/G) \times (X/G)} d_X(s(x, y)\cdot (g, g))^2 \:
d\pi(x, y) \: dh(g) \\
& = \: \int_{(X/G) \times (X/G)} d_{X/G}(x, y)^2 \:
d\pi(x, y).\notag
\end{align}
Thus $p_*$ and $(p_*)^{-1}$ are metrically nonincreasing,
which shows that $p_*$ defines an isometric isomorphism between
$P_2(X)^G$ and $P_2(X/G)$.
\end{proof}

\begin{proof}[Proof of Theorem~\ref{Gaction}]
The proofs of parts a. and b. of the theorem are similar,
so we will be content with proving just part a.

First, $(X/G, d_{X/G})$ is a length space. (Given $x, y \in X/G$, let
$\widetilde{x} \in p^{-1}(x)$ and $\widetilde{y} \in p^{-1}(y)$ satisfy
$d_X(\widetilde{x}, \widetilde{y}) \: = \: d_{X/G}(x, y)$. If
$c$ is a geodesic from $\widetilde{x}$ to $\widetilde{y}$ then
$p \circ c$ is a geodesic from ${x}$ to ${y}$.)

Given $\mu_0, \mu_1 \in P_2^{ac}(X/G, p_* \nu)$, write
$\mu_0 \: = \: \rho_0 \: p_* \nu$ and
$\mu_1 \: = \: \rho_1 \: p_* \nu$. Put
$\widetilde{\mu}_0 \: = \: (p^* \rho_0) \: \nu$ and
$\widetilde{\mu}_1 \: = \: (p^* \rho_1) \: \nu$. From Lemma~\ref{lemaction},
$W_2(\widetilde{\mu}_0, \widetilde{\mu}_1) \: = \:W_2(\mu_0, \mu_1)$.
By hypothesis,
there is a Wasserstein geodesic $\{ \widetilde{\mu}_t \}_{t \in [0,1]}$
from $\widetilde{\mu}_0$ to $\widetilde{\mu}_1$ so that for all
$U \in \DC_N$, equation (\ref{defldc}) is satisfied
along $\{ \widetilde{\mu}_t \}_{t \in [0,1]}$,
with $\lambda \: = \: 0$. The geodesic 
$\{ \widetilde{\mu}_t \}_{t \in [0,1]}$ is $G$-invariant, as otherwise
by applying an appropriate 
element of $G$ we would obtain two distinct Wasserstein
geodesics between $\widetilde{\mu}_0$ and $\widetilde{\mu}_1$.
Put $\mu_t \: = \: p_* \widetilde{\mu}_t$. It follows from the above discussion
that $\{\mu_t\}_{t \in [0, 1]}$ is a curve with length
$W_2(\mu_0, \mu_1)$, and so is a Wasserstein geodesic. 
As $\widetilde{\mu}_t \in P_2^{ac}(X, \nu)$, we have
${\mu}_t \in P_2^{ac}(X/G, p_* \nu)$.
Write
$\mu_t \: = \: \rho_t \: p_*\nu$. Then
$\widetilde{\mu}_t \: = \: (p^* \rho_t) \: \nu$.
As
\begin{equation}
U_{p_* \nu}(\mu_t) \: = \:
\int_{X/G} U(\rho_t) \: dp_* \nu \: = \:
\int_{X} p^*U(\rho_t) \: d\nu \: = \:
\int_{X} U(p^* \rho_t) \: d\nu \: = \: 
U_\nu(\widetilde{\mu}_t),
\end{equation}
it follows
that equation (\ref{defldc}) is satisfied
along $\{{\mu}_t \}_{t \in [0,1]}$,
with $\lambda \: = \: 0$. Along with
Proposition \ref{lemsuffdc}, this concludes the proof of part a.
\end{proof}

\subsection{Uniform integrability and absolute continuity}

In what has been done so far, it would be logically consistent to
make our definition of nonnegative $N$-Ricci curvature to mean
weak displacement convexity of just $H_{N, \nu}$, and not 
necessarily all of $\DC_N$. The reasons to require weak
displacement convexity for $\DC_N$ are first that we can,
in the sense of being consistent with the classical definitions
in the case of a Riemannian manifold, and second that
we thereby obtain a useful absolute continuity property for the 
measures appearing in a Wasserstein geodesic joining two
absolutely continuous measures.  This last property will
feed into Proposition \ref{displconvineq}, when proving 
Theorem \ref{geninequality}.

\begin{lemma} \label{lemexistU} Let $\{\mu_i\}_{i=1}^m$ be a finite
subset of $P^{ac}_2(X, \nu)$,
with densities $\rho_i \: = \: \frac{d\mu_i}{d\nu}$.
If $N<\infty$ then there is a function $U\in \DC_N$ such that
\begin{equation} 
\lim_{r\to\infty} \frac{U(r)}r = \infty
\end{equation}
and
\begin{equation} \sup_{1\leq i\leq m} \int_X U(\rho_i(x))\,d\nu(x) < \infty.
\end{equation}
\end{lemma}
\begin{proof}
As a special case of the Dunford-Pettis theorem
\cite[Theorem 2.12]{Buttazzo-Giaquinta-Hildebrandt (1998)},
there is an increasing function $\Phi \: : \: (0, \infty) 
\rightarrow \R$ such that
\begin{equation} 
\lim_{r\to\infty} \frac{\Phi(r)}{r} = \infty
\end{equation}
and 
\begin{equation} 
\sup_{1\leq i\leq m} \int_X \Phi(\rho_i(x))\,d\nu(x) <\infty.
\end{equation}
We may assume that $\Phi$ is identically zero on $[0,1]$.

Consider the function $\psi \: : \: (0, \infty) \rightarrow \R$ given by
\begin{equation} 
\psi(\lambda) \: = \: \lambda^N \Phi(\lambda^{-N}).
\end{equation}
Then $\psi\equiv 0$ on $[1,\infty)$, and 
$\lim_{\lambda \rightarrow 0^+} \psi(\lambda) \: = \: \infty$.
Let $\tilde{\psi}$ be
the lower convex hull of $\psi$ on $(0, \infty)$, 
i.e. the supremum of the linear
functions bounded above by $\psi$. Then
$\tilde{\psi} \equiv 0$ on $[1,\infty)$ and $\tilde{\psi}$ is
nonincreasing.
We claim that $\lim_{\lambda\rightarrow 0^+} \tilde{\psi}(\lambda)
\: = \: \infty$. If not, suppose that 
$\lim_{\lambda\rightarrow 0^+} \tilde{\psi}(\lambda)
\: = \: M \: < \: \infty$. Let $a \: = \: \sup_{\lambda \ge 0} \:
\frac{M+1-\psi(\lambda)}{\lambda} < \infty$ (because this quantity
is $\leq 0$ when $\lambda$ is small enough). Then $\psi(\lambda) \: \ge \:
M+1-a\lambda$, so 
$\lim_{\lambda\rightarrow 0^+} \tilde{\psi}(\lambda) \: \ge \: M+1$, 
which is a contradiction.

Now set
\begin{equation} 
U(r) \: = \:  r \: \tilde{\psi}(r^{-1/N}).
\end{equation}
Since $\tilde{\psi}\leq \psi$ and $\Phi(r) = r \psi(r^{-1/N})$,
we see that $U\leq \Phi$. Hence
\begin{equation} 
\sup_{1\leq i\leq m} \int_X U(\rho_i(x))\,d\nu(x) <\infty.
\end{equation}
Since $\lim_{\lambda\rightarrow 0^+} \tilde{\psi}(\lambda)
\: = \: \infty$,
we also know that 
\begin{equation} 
\lim_{r \rightarrow \infty} \frac{U(r)}{r} \: = \: \infty.
\end{equation}
Clearly $U$ is continuous with $U(0) \: = \: 0$.
As $\tilde{\psi}$ is convex and nonincreasing, it follows that
$U$ is convex. Hence $U\in \DC_N$.
\end{proof}

\begin{theorem} \label{thmgeodcvx}
If $(X,d,\nu)$ has nonnegative $N$-Ricci curvature 
for some $N \in [1, \infty)$
then $P^{ac}_2(X, \nu)$ is a convex subset of $P_2(X)$.
\end{theorem}
\begin{proof} Given $\mu_0, \mu_1 \in P^{ac}_2(X, \nu)$, put
$\rho_0 \: = \: \frac{d\mu_0}{d\nu}$ and
$\rho_1 \: = \: \frac{d\mu_1}{d\nu}$.
By Lemma~\ref{lemexistU}, there is a
$U\in \DC_N$ with $U^\prime(\infty) \: = \: \infty$
such that $U_\nu(\mu_0) < \infty$ and $U_\nu(\mu_1) < \infty$.
As $(X,d,\nu)$ has nonnegative $N$-Ricci curvature, 
there is a Wasserstein geodesic 
$\{\mu_t\}_{t \in [0,1]}$ from $\mu_0$
to $\mu_1$ so that (\ref{defldc}) is satisfied
with $\lambda \: = \: 0$.
In particular, $U_\nu(\mu_t) \: < \: \infty$ for all $t \in [0,1]$.
As $U^\prime(\infty) \: = \: \infty$, it follows that
$\mu_t \in P^{ac}_2(X, \nu)$ for each $t$.
\end{proof}

We now clarify the relationship between $(X, d, \nu)$ having 
nonnegative $N$-Ricci curvature and the analogous statement for
$\supp(\nu)$. We recall the notion of
a subset $A \subset X$ being convex or totally convex, 
from Section \ref{lengthspaces}, and
we note that $d \big|_A$ defines a length space structure on
a closed subset $A$ if and only if $A$ is convex in $X$.

\begin{theorem} \label{thmredspt}
a. Given $N\in [1,\infty)$, suppose that a compact measured
length space $(X,d,\nu)$  has nonnegative
$N$-Ricci curvature.
Then $\supp(\nu)$ is a convex subset of $X$ 
(although not necessarily totally convex) 
and $(\supp(\nu), d |_{\supp(\nu)}, \nu)$ has nonnegative $N$-Ricci
curvature.  Conversely, if $\supp(\nu)$ is a convex subset of $X$ and
$(\supp(\nu), d |_{\supp(\nu)}, \nu)$ has nonnegative $N$-Ricci curvature
then $(X,d,\nu)$ has nonnegative $N$-Ricci curvature.
\smallskip

b. Given $K \in \R$, the analogous statement holds 
when one replaces ``nonnegative $N$-Ricci curvature'' by ``$\infty$-Ricci curvature 
bounded below by $K$''.
\end{theorem}
\begin{proof}a. 
Let $(X,d,\nu)$ be a compact measured
length space with nonnegative $N$-Ricci curvature.
Let $\mu_0$ and $\mu_1$ be
elements of $P_2(X,\nu)$. By Theorem~\ref{thmapproxU}
in Appendix~\ref{secmoll},
there are sequences $\{\mu_{k,0}\}_{k=1}^\infty$ and
$\{\mu_{k,1}\}_{k=1}^\infty$ in $P^{ac}_2(X,\nu)$ (in fact with
continuous densities) such that
$\lim_{k \rightarrow \infty} \mu_{k,0} \: =\: \mu_0$,
$\lim_{k \rightarrow \infty} \mu_{k,1} \: =\: \mu_1$ and
for all $U\in\DC_N$,
$\lim_{k \rightarrow \infty} U_\nu(\mu_{k,0}) \: = \: U_\nu(\mu_0)$
and
$\lim_{k \rightarrow \infty} U_\nu(\mu_{k,1}) \: = \: U_\nu(\mu_1)$.
From the definition of nonnegative $N$-Ricci,
for each $k$ there is a Wasserstein geodesic $\{\mu_{k,t}\}_{t \in [0,1]}$
such that
\begin{equation}
 U_\nu(\mu_{k,t})\leq t\, U_\nu(\mu_{k,1}) \: +\:  
(1-t)\, U_\nu(\mu_{k,0})
\end{equation}
for all $U\in\DC_N$ and $t\in [0,1]$.
By repeating the proof of Theorem~\ref{thmgeodcvx}, each $\mu_{k,t}$
is absolutely continuous with respect to $\nu$. In particular,
it is supported in $\supp(\nu)$. By the same reasoning as in
the proof of Proposition~\ref{lemsuffdc}, after passing to a subsequence
we may assume that as $k \rightarrow \infty$, the geodesics 
$\{\mu_{k,t}\}_{t \in [0,1]}$ converge uniformly to a Wasserstein geodesic
$\{\mu_t\}_{t \in [0,1]}$ that satisfies
\begin{equation} \label{char}
 U_\nu(\mu_{t})\leq t\, U_\nu(\mu_{1}) \: +\:  
(1-t)\, U_\nu(\mu_{0}).
\end{equation}
For each $t \in [0,1]$, the measure $\mu_t$ is the
weak-$*$ limit of the probability measures 
$\{\mu_{k,t}\}_{k=1}^\infty$, which are all supported in the
closed set $\supp(\nu)$. Hence $\mu_t$ is also supported in $\supp(\nu)$.
To summarize, we have shown that 
$\{\mu_t\}_{t\in [0,1]}$ is a Wasserstein geodesic lying in
$P_2(X,\nu)$ that satisfies (\ref{char}) for
all $U\in\DC_N$ and $t\in [0,1]$.

We now check that $\supp(\nu)$ is convex.
Let $x_0$ and $x_1$ be any two points in $\supp(\nu)$.
Applying the reasoning above to $\mu_0=\delta_{x_0}$ and
$\mu_1=\delta_{x_1}$, one obtains the existence of a Wasserstein
geodesic $\{\mu_t\}_{t \in [0,1]}$ joining $\delta_{x_0}$ to $\delta_{x_1}$
such that each $\mu_t$ is supported in $\supp(\nu)$.
By Proposition~\ref{everygeod}, there is an optimal dynamical transference plan
$\Pi \in P(\Gamma)$ such that $\mu_t=(e_t)_*\Pi$ for all
$t \in [0,1]$. For each $t \in [0,1]$,
we know that $\gamma(t)\in \supp(\nu)$ holds 
for $\Pi$-almost all $\gamma$.
It follows that for $\Pi$-almost all $\gamma$, we have
$\gamma(t)\in \supp(\nu)$ for all 
$t\in\Q\cap [0,1]$. As $\gamma \in \Gamma$
is continuous, this is the same as saying that 
for $\Pi$-almost all $\gamma$, the geodesic
$\gamma$ is entirely contained in $\supp(\nu)$. Also,
for $\Pi$-almost all $\gamma$
we have $\gamma(0)=x_0$ and $\gamma(1)=x_1$. Thus $x_0$ and $x_1$
are indeed joined by a geodesic path contained in $\supp(\nu)$.

This proves the direct implication in part a. The converse is
immediate. \\
b. The proof of part b. follows the same lines as that of part a. 
We construct the
approximants $\{\mu_{k,0}\}_{k=1}^\infty$ and
$\{\mu_{k,1}\}_{k=1}^\infty$, with continuous densities, and
the geodesics $\{\mu_{k,t}\}_{t \in [0,1]}$. As
$H_{\infty,\nu}(\mu_{0,k})<\infty$ and 
$H_{\infty,\nu}(\mu_{1,k})<\infty$, we can apply 
inequality (\ref{defldc}) with $U=H_\infty$ and $\lambda \: =\: K$, 
to deduce that
$H_{\infty,\nu}(\mu_{t,k})<\infty$ for all $t\in [0,1]$.
This implies that
$\mu_{t,k}$ is absolutely continuous with respect to $\nu$.
The rest of the argument is similar to that of part a.
\end{proof}

\begin{remark} Corollary~\ref{MGHcompactness} and Theorem~\ref{thmredspt}.a 
together show that in the case $N < \infty$, 
we do not lose much by assuming
that $X=\supp(\nu)$.
\end{remark}

\section{Log Sobolev, Talagrand and Poincar\'e inequalities} \label{secineq}

In this section we study several functional inequalities with geometric
content that are associated to optimal transport and concentration
of measure : log Sobolev inequalities, Talagrand
inequalities and Poincar\'e inequalities. We refer to
\cite{Toulouse (2000)} and \cite[Chapter~9]{Villani (2003)}
for concise surveys about previous work on  these inequalities.

We first write some general functional inequalities.  In the case of
$\infty$-Ricci curvature bounded below by $K$, we make explicit
the ensuing log Sobolev
inequalities, Talagrand inequalities and Poincar\'e inequalities.
We then write out explicit functional
inequalities in the case of nonnegative $N$-Ricci curvature.
Finally, we prove a weak Bonnet--Myers theorem, following 
\cite[Section 6]{Otto-Villani (2000)}.

\subsection{The general inequalities}

We recall the generalized Fisher information of
(\ref{fisherdef}), where $\rho \in \Lip(X)$ is positive
and $\mu \: = \: \rho  \nu$ is the corresponding measure.

\begin{theorem} \label{geninequality}
Suppose that 
$(X, d, \nu)$ has $\infty$-Ricci curvature bounded below by $K > 0$.
Then for all $\mu \in P_2(X,\nu)$,
\begin{equation} 
\frac{K}{2} \: W_2(\mu, \nu)^2 \: \leq \:
H_{\infty,\nu}(\mu).
\end{equation}

If now $\mu \in P^{ac}_2(X, \nu)$ and its density 
$\rho \: = \: \frac{d\mu}{d\nu}$ is a 
positive Lipschitz function on $X$ then
\begin{equation} 
H_{\infty,\nu}(\mu) \: \leq \: W_2(\mu, \nu) \: \sqrt{I_{\infty,\nu}(\mu)} 
\: - \:  
\frac{K}{2} \: W_2(\mu, \nu)^2 
\: \leq \: \frac{1}{2K} \: I_{\infty,\nu}(\mu).
\end{equation}
If on the other hand
$(X, d, \nu)$ has $\infty$-Ricci curvature bounded below by $K \le 0$
then
\begin{equation} 
H_{\infty,\nu}(\mu) \: \leq  \: \diam(X) \: \sqrt{I_{\infty,\nu}(\mu)} 
\: - \:  
\frac{K}{2} \: \diam(X)^2. 
\end{equation}
If $(X, d, \nu)$ has nonnegative $N$-Ricci curvature then
\begin{equation} \label{finiteN}
 H_{N,\nu}(\mu) \: \leq  \: \diam(X) \: \sqrt{I_{N,\nu}(\mu)}.
\end{equation}
\end{theorem}
\begin{proof}
We wish to apply Proposition \ref{displconvineq} to the cases described in
Particular Cases 
\ref{partcases}. Under the assumption that $U_\nu(\mu) < \infty$,
we have to show that there is a Wasserstein geodesic as in the statement of
Proposition \ref{displconvineq} with $\mu_t \in P^{ac}_2(X, \nu)$
for all $t \in [0,1]$. If $N = \infty$ then there is some
Wasserstein geodesic $\{\mu_t\}_{t \in [0,1]}$
from $\mu$ to $\nu$ which in particular satisfies 
equation 
(\ref{defldc}) with $U_\nu \: = \: H_{\infty,\nu}$ and $\lambda \: = \: K$.
Hence $H_{\infty,\nu}(\mu_t) < \infty$ for all $t \in [0,1]$ and 
the claim follows from the fact that 
$U_\infty^\prime(\infty) \: = \: \infty$.
If $N \in [1, \infty)$ then the claim follows from Theorem
\ref{thmgeodcvx}.
\end{proof}

We now express the conclusion of Theorem \ref{geninequality} 
in terms of more standard
inequalities, starting with the case $N = \infty$.

\subsection{The case $N = \infty$}

\begin{definition}
Suppose that $K > 0$.

$\bullet$ We say that $\nu$ satisfies a
log Sobolev inequality with constant $K$, $\LSI(K)$,
if for all $\mu \in P^{ac}_2(X, \nu)$ whose density 
$\rho \: = \: \frac{d\mu}{d\nu}$ is Lipschitz and positive, we have
\begin{equation} 
H_{\infty,\nu} (\mu) \leq \frac1{2K}\, I_{\infty,\nu}(\mu).
\end{equation}

$\bullet$ We say that $\nu$ satisfies a Talagrand inequality with
constant $K$, $T(K)$, if for all $\mu \in P_2(X,\nu)$,
\begin{equation} 
W_2(\mu,\nu) \leq \sqrt{\frac{2 H_{\infty,\nu} (\mu)}{K}}. 
\end{equation}

$\bullet$ We say that $\nu$ satisfies a Poincar\'e inequality
with constant $K$, $P(K)$, if for 
all $h \in \Lip(X)$ with $\int_X h\,d\nu = 0$, we have
\begin{equation} 
\int_X h^2\,d\nu \leq 
\frac1{K} \int_X |\nabla^- h|^2\,d\nu. 
\end{equation}
\end{definition}

\begin{remark}
Here we used the gradient norm defined in~\eqref{nablamoins},
instead of the one defined in~\eqref{nablaf}. Accordingly, our
log Sobolev inequality and Poincar\'e inequalities are slightly 
stronger statements than those discussed by some other authors.
\end{remark}

All of these inequalities are associated with concentration of 
measure~\cite{Toulouse (2000),Bobkov-Gentil-Ledoux 
(2001),Bobkov-Gotze (1999),Ledoux (1999), Ledoux (2001)}. 
For example, $T(K)$
implies a Gaussian-type concentration of measure.
The following chain of implications, none of which is an equivalence,
is well-known in the context of smooth Riemannian manifolds :
\begin{equation} 
[\Ric \geq K ] \Longrightarrow \LSI(K)
\Longrightarrow T(K) \Longrightarrow P(K). 
\end{equation}

In the context of length spaces, we see from Theorem \ref{geninequality} that
having $\infty$-Ricci curvature bounded below by $K > 0$ implies
$\LSI(K)$ and $T(K)$. The next corollary makes the statement of the
log Sobolev inequality more explicit.

\begin{corollary} \label{lscor}
Suppose that $(X,d,\nu)$ has $\infty$-Ricci curvature bounded below by
$K \in \R$.
If $f \in \Lip(X)$ satisfies $\int_X f^2 \: d\nu \: = \: 1$
then
\begin{equation}
\int_X f^2 \: \log(f^2) \: d\nu \: \le \: 2 \: W_2(f^2 \: \nu, \nu) \:
\sqrt{\int_X |\nabla^- f|^2 \: d\nu} \: - \: 
\frac{K}2 \: W_2(f^2 \: \nu, \nu)^2.
\end{equation}

In particular, if $K > 0$ then
\begin{equation} \label{ls}
\int_X f^2 \: \log(f^2) \: d\nu \: \le \: \frac{2}{K} \:
\int_X |\nabla^- f|^2 \: d\nu,
\end{equation}
while if $K \: \le \: 0$ then
\begin{equation}
\int_X f^2 \: \log(f^2) \: d\nu \: \le \: 2 \: \diam(X) \:
\sqrt{\int_X |\nabla^- f|^2 \: d\nu} \: - \: 
\frac{K}2 \: \diam(X)^2.
\end{equation}
\end{corollary}
\begin{proof}
For any $\epsilon > 0$, put $\rho_\epsilon \: = \:
\frac{f^2 +  \epsilon}{1 \: + \: \epsilon}$.
{From} Theorem \ref{geninequality},
\begin{equation} \int_X \rho_\epsilon \: \log(\rho_\epsilon) \: d\nu
\: \leq  \:
W_2(\rho_\epsilon \: \nu,\nu) \sqrt{\int_X 
\frac{|\nabla^- \rho_\epsilon|^2}{\rho_\epsilon} \: d\nu} \: - \: 
\frac{K}2 \: W_2(\rho_\epsilon \nu,\nu)^2. 
\end{equation}
As 
\begin{equation}
\frac{|\nabla^- \rho_\epsilon|^2}{\rho_\epsilon} \: = \:
\frac{1}{1+\epsilon} \: \frac{4f^2}{f^2 \: + \: \epsilon} \: |\nabla^- f|^2,
\end{equation}
the corollary follows by taking $\epsilon \rightarrow 0$.
\end{proof}

We now recall the standard fact that $\LSI(K)$ implies $P(K)$.

\begin{theorem} \label{LSIP}
Let $(X,d,\nu)$ be a compact measured length space
satisfying $\LSI(K)$ for 
some $K>0$. Then it also satisfies $P(K)$.
\end{theorem}
\begin{proof}
Suppose that $h \in \Lip(X)$ satisfies $\int_X h\,d\nu = 0$.
For $\epsilon \in [0, \frac{1}{\parallel h \parallel_\infty})$, put 
$f_\epsilon \: = \: \sqrt{1 \: + \: \epsilon \: h} \: > \: 0$.
As $2 f_\epsilon \nabla^- f_\epsilon \: = \: \epsilon \: \nabla^- h$, 
it follows that
\begin{equation}
\lim_{\epsilon \rightarrow 0^+} \left(\frac{1}{\epsilon^2} \:
\int_X |\nabla^- f_\epsilon|^2 \: d\nu\right) \: = \:
\frac14 \: \int_X |\nabla^- h|^2 \: d\nu.
\end{equation}
As the Taylor expansion of $x \log(x) \: - \: x \: + \: 1$ around
$x =1$ is $\frac12 (x-1)^2 \: + \: \ldots$, it follows that
\begin{equation} 
\lim_{\epsilon \rightarrow 0^+} \frac{1}{\epsilon^2} \:
\int_X f_\epsilon^2 \: \log(f_\epsilon^2) \: d\nu \: = \: \frac12 \:
\int_X h^2\ \: d\nu.
\end{equation}
Then the conclusion follows from (\ref{ls}).
\end{proof}

\begin{remark}
If $Q(h) \: = \: \int_X |\nabla h|^2 \: d\nu$ defines a quadratic form on
$\Lip(X)$, which in addition is closable in $L^2(X, \nu)$, then
there is a (nonpositive) self-adjoint Laplacian $\triangle_\nu$
associated to $Q$. In this case, $P(K)$ implies that $-\triangle_\nu
\: \ge \: K$ on the orthogonal complement of the constant
functions.

We do not claim to show that there are such Laplacians on $(X, d, \nu)$ in
general.  In the case of a limit space arising from a sequence of 
manifolds with Ricci curvature bounded below, Cheeger and Colding used 
additional structure on the limit space in order to show
the Laplacian does exist \cite{Cheeger-Colding (2000)}.
\end{remark}

As mentioned above, in the case of smooth Riemannian manifolds there are
stronger implications: $T(K)$ implies $P(K)$, and $\LSI(K)$ implies $T(K)$.
We will show elsewhere that the former is always true,
while the latter is true under the additional assumption of a lower bound
on the Alexandrov curvature:

\begin{theorem} \label{LSIT}
Let $(X,d,\nu)$ be a compact measured length space.
\sm

(i) If $\nu$ satisfies $T(K)$ for some $K>0$, then it also satisfies $P(K)$.
\sm

(ii) If $X$ is a finite-dimensional Alexandrov space with Alexandrov
curvature bounded below, and $\nu$ satisfies $\LSI(K)$ for some
$K > 0$, then it also satisfies $T(K)$.
\end{theorem}

\begin{remark} The Alexandrov curvature bound in (ii) essentially serves as
a regularity assumption. One can ask whether it can be weakened.
\end{remark}

\begin{remark}
We have only discussed {\em global} Poincar\'e inequalities. There is also
a notion of a metric-measure space admitting a {\em local} Poincar\'e inequality, 
as considered for example in
\cite{Cheeger-Colding (2000)}. If a measured length space $(X, d, \nu)$ has
nonnegative $N$-Ricci curvature, with $N < \infty$, then it admits a
local Poincar\'e inequality, at least if one assumes
almost-everywhere uniqueness of geodesics. We will discuss this in detail
elsewhere.
\end{remark}

\subsection{The case $N < \infty$}

We now write an analog of Corollary \ref{lscor} in the case $N < \infty$.
Suppose that $(X,d,\nu)$ has nonnegative $N$-Ricci curvature.
Then if $\rho$ is a positive Lipschitz function on $X$,
(\ref{finiteN}) says that
\begin{equation} \label{finiteNineq}
N \: - \: N \int_X \rho^{1- \frac{1}{N}} \: d\nu \: \le \: \frac{N-1}{N} \:
\diam(X) \: 
\sqrt{\int_X \frac{|\nabla^- \rho|^2}{\rho^{\frac{2}{N}+1}} \,d\nu}.
\end{equation}
If $N > 2$, put
$f \: = \: \rho^{\frac{N-2}{2N}}$. Then 
$\int_X f^{\frac{2N}{N-2}} \: d\nu \: = \: 1$ and
one finds that
(\ref{finiteNineq}) is equivalent to
\begin{equation} \label{finiteNineq2}
1 \: - \: \int_X f^{\frac{2(N-1)}{N-2}} \: d\nu \: \le \: 
\frac{2(N-1)}{N(N-2)} \:
\diam(X) \: 
\sqrt{\int_X |\nabla^- f|^2 \,d\nu}.
\end{equation}
As in the proof of Corollary \ref{lscor},
equation (\ref{finiteNineq2}) holds for all
$f \in \Lip(X)$ satisfying $\int_X f^{\frac{2N}{N-2}} \: d\nu \: = \: 1$.
From H\"older's inequality,
\begin{equation}
\int_X f^{\frac{2(N-1)}{N-2}} \: d\nu \: \le \:
\left( \int_X f \: d\nu \right)^{\frac{2}{N+2}} \:
\left( \int_X f^{\frac{2N}{N-2}} \: d\nu \right)^{\frac{N}{N+2}} \: = \:
\left( \int_X f \: d\nu \right)^{\frac{2}{N+2}}. 
\end{equation}
Then (\ref{finiteNineq2}) implies
\begin{equation} \label{finiteNineq4}
1 \: \le \: \frac{2(N-1)}{N(N-2)} \:
\diam(X) \: \sqrt{\int_X |\nabla^- f|^2 \,d\nu} \: + \:
\left( \int_X f \: d\nu \right)^{\frac{2}{N+2}}.
\end{equation}
Writing~\eqref{finiteNineq4} in a homogeneous form, 
one sees that its content is as follows: for a function $F$ 
on $X$, bounds on $\parallel \nabla^- F \parallel_{2}$ and
$\parallel F \parallel_{1}$ imply a bound on $\parallel F 
\parallel_{\frac{2N}{N-2}}$.
This is of course an instance of Sobolev embedding.

If $N = 2$, putting $f \: = \: \log(\frac{1}{\rho})$, one finds that
$\int_X e^{-f} \: d\nu \: = \: 1$ implies
\begin{equation}
1 \: - \: \int_X e^{-\frac{f}{2}} \: d\nu \: \leq \: \frac14 \:
\diam(X) \: 
\sqrt{\int_X |\nabla^- f|^2 \,d\nu}.
\end{equation}

\subsection{Weak Bonnet--Myers theorem}

The classical Bonnet--Myers theorem says that
if $M$ is a smooth connected complete
$N$-dimensional Riemannian manifold with $\Ric_M \: \geq K g_M \: > \: 0$, 
then $\diam(M) \leq \pi \sqrt{\frac{N-1}{K}}$.

We cannot give an immediate generalization of this theorem to
a measured length space $(X,d,\nu)$, as we have not defined what it means to
have $N$-Ricci curvature bounded below by $K$ for $N < \infty$ and
$K > 0$. However, it does make sense
to state a weak version of the 
Bonnet--Myers theorem under the assumptions that $(X,d,\nu)$ has
nonnegative $N$-Ricci curvature and has $\infty$-Ricci curvature
bounded below by $K > 0$. 

\begin{theorem}\label{thmBM}
There is a constant $C > 0$ with the following property.
Let $(X,d,\nu)$ be a compact measured length space with 
nonnegative $N$-Ricci curvature, and $\infty$-Ricci curvature
bounded below by $K > 0$. Suppose that $\supp(\nu) \:= \: X$.
Then
\begin{equation} 
\diam(X) \leq C \sqrt{\frac{N}{K}}.
\end{equation}
\end{theorem}

\begin{proof} {From} Theorem~\ref{thmBG}, $\nu$ satisfies the growth estimate
\begin{equation} 
\frac{\nu(B_r(x))}{\nu(B_{\alpha r}(x))} \leq \alpha^{-N}, \qquad
0<\alpha\leq 1.
\end{equation}
{From} Theorem~\ref{geninequality}, $\nu$ satisfies $T(K)$. The result follows
by repeating verbatim the proof of~\cite[Theorem~4]{Otto-Villani (2000)}
with $R=0$, $n=N$ and $\rho=K$.
\end{proof}

\begin{remark} The remark at the end of~\cite[Section~6]{Otto-Villani (2000)}
shows that $C= 7.7$ is admissible.
\end{remark}

\section{The case of Riemannian manifolds} \label{Riemcase}

In this section we look at the case of a smooth Riemannian manifold $(M,g)$
equipped with a smooth measure $\nu$. We define the tensor $\Ric_N$ and
show the equivalence of lower bounds on $\Ric_N$ to various displacement
convexity conditions.  In particular, we show that the measured
length space $(M,g,\nu)$ has nonnegative $N$-Ricci curvature if and only
if $\Ric_N \: \ge \: 0$, and that it has $\infty$-Ricci curvature
bounded below by $K$ if and only if $\Ric_\infty \: \ge \: K \: g$.

We use this, along with Theorem \ref{stabRicci}, to 
characterize measured Ricci limit spaces that happen to be smooth.
We give some consequences concerning their metric structure.
We then
show that for Riemannian manifolds, lower $N$-Ricci curvature bounds are
preserved under taking compact quotients.
Finally, we use displacement convexity
to give a ``synthetic'' proof
of a part of the Ricci O'Neill theorem from \cite{Lott (2003)}.

\subsection{Formulation of $N$-Ricci curvature in classical terms}

Let $(M, g)$ be a smooth compact connected $n$-dimensional
Riemannian manifold. Let $\Ric$ denote its Ricci tensor.  

Given $\Psi \in C^\infty(M)$ with $\int_M e^{-\Psi} \: \dvol_M \: = \: 1$, 
put $d\nu \: = \: e^{- \Psi} \: \dvol_M$.

\begin{definition}
For $N \in [1, \infty]$, the $N$-Ricci tensor of $(M, g, \nu)$
is
\begin{equation}
\Ric_N \: = 
\begin{cases}
\Ric \: + \: \Hess(\Psi)  & \text{ if $N = \infty$}, \\
\Ric \: + \: \Hess(\Psi) \: - \: \frac{1}{N-n} \: d \Psi \otimes
d \Psi & \text{ if $n \: < \: N \: < \: \infty$}, \\
\Ric \: + \: \Hess(\Psi) \: - \: \infty \: (d \Psi \otimes
d \Psi) & \text{ if $N = n$}, \\
- \infty & \text{ if $N < n$,}
\end{cases}
\end{equation}
where by convention $\infty \cdot 0 \: = \: 0$.
\end{definition}

The expression for $\Ric_\infty$ is the Bakry-\'Emery tensor
\cite{Bakry-Emery (1985)}. The expression for $\Ric_N$ with
$n < N < \infty$ was considered in \cite{Lott (2003),Qian (1997)}.
The statement $\Ric_N\geq K g$ is equivalent to the statement
that the operator $L \: = \: \triangle \: - \:
(\nabla \Psi) \cdot \nabla$ satisfies Bakry's
curvature-dimension condition ${\rm CD}(K,N)$
\cite[Proposition 6.2]{Bakry (1994)}.

Given $K \in \R$,
we recall the definition of $\lambda \: : \: \DC_\infty \rightarrow \R
\cup \{- \infty\}$
from Definition \ref{lambdadef}.

\begin{theorem} \label{Riccimfolds}
a. For $N \in (1, \infty)$, the following are equivalent : \\
(1) $\Ric_N \: \ge \: 0$. \\
(2) The measured length space
$(M, g, \nu)$ has nonnegative $N$-Ricci curvature. \\
(3) For all $U \in \DC_N$, $U_\nu$ is weakly displacement convex on 
$P_2(M)$.\\
(4) For all $U \in \DC_N$, $U_\nu$ is weakly a.c. displacement convex on 
$P^{ac}_2(M)$.\\
(5) $H_{N,\nu}$ is weakly a.c. displacement convex on $P^{ac}_2(M)$.\\ \\
b. For any $K \in \R$, the following are equivalent : \\
(1) $\Ric_\infty \: \ge \: K \: g$. \\
(2) The measured length space
$(M, g, \nu)$ has $\infty$-Ricci curvature bounded below by $K$. \\
(3) For all $U \in \DC_\infty$, $U_\nu$ is weakly 
$\lambda(U)$-displacement convex on 
$P_2(M)$.\\
(4) For all $U \in \DC_\infty$, 
$U_\nu$ is weakly $\lambda(U)$-a.c. displacement convex on 
$P^{ac}_2(M)$.\\
(5) $H_{\infty,\nu}$ is weakly $K$-a.c. 
displacement convex on $P^{ac}_2(M)$.
\end{theorem}

For both parts (a) and (b), the nontrivial implications are
$(1) \Rightarrow (2)$ and $(5) \Rightarrow (1)$. The proof that
$(1) \Rightarrow (2)$ will be along the lines of
\cite[Theorem 6.2]{Cordero-Erausquin-McCann-Schmuckenschlager (2001)}, 
with some differences. One ingredient is the following lemma.

\begin{lemma} \label{lemD''} 
Let $\phi:M\to\R$ be a $\frac{d^2}{2}$-concave function. We recall that $\phi$ is
necessarily Lipschitz and hence $(\nabla \phi)(y)$ exists for almost all
$y \in $M. For such $y$, define
\begin{equation}
F_t(y)\equiv \exp_y(-t\nabla\phi(y)).
\end{equation}
Assume furthermore that $y\in M$ is such that 
\sm

(i) $\phi$ admits a Hessian at $y$ (in the sense of Alexandrov),
\sm

(ii) $F_t$ is differentiable at $y$ for all $t\in [0,1)$ and
\sm

(iii) $dF_t(y)$ is nonsingular for all $t\in [0,1)$.
\sm

Then $D(t)\equiv \det^{\frac1{n}}(dF_t(y))$ satisfies the 
differential inequality
\begin{equation}
\frac{D''(t)}{D(t)} \leq -\frac1{n} \Ric (F'_t(y), F'_t(y))
\qquad t\in (0,1).
\end{equation}
\end{lemma}

\begin{proof}
Let $\{e_i\}_{i=1}^n$ be an orthonormal basis of $T_yM$. For each
$i$, let $J_i(t)$ be defined by
\begin{equation}
 J_i(t) = (dF_t)_y \: (e_i).
 \end{equation}
Then $\{J_i(t)\}_{i=1}^n$ is a Jacobi field with
$J_i(0)\: = \: e_i$. Next, we note that $d\phi$ is differentiable at $y$,
and that $d(d\phi)_y$ coincides with $\Hess_y(\phi)$,
up to identification. This is not so obvious (indeed, the existence
of a Hessian only means the existence of a second-order Taylor expansion)
but can be shown as a consequence of the semiconcavity of $\phi$,
as in~\cite[Proposition~4.1 (b)]
{Cordero-Erausquin-McCann-Schmuckenschlager (2001)}.
(The case of a convex function in $\R^n$ is treated
in~\cite[Theorems~3.2 and~7.10]{Alberti-Ambrosio (1999)}.)
It follows that
\begin{equation}
J'_i(0)\: = \: - \Hess(\phi)(y)\: e_i.
\end{equation}

Let now $W(t)$ be the $(n \times n)$-matrix with
\begin{equation}
W_{ij}(t) \: = \: \langle J_i(t), J_j(t) \rangle;
\end{equation}
then ${\det}^{\frac{1}{n}}(dF_t)(y) \: = \: 
{\det}^{\frac{1}{2n}} W(t)$. 

Since $W(t)$ is nonsingular for $t\in [0,1)$, $\{J_i(t)\}_{i=1}^n$ 
is a basis of $T_{F_t(y)}M$. Define a matrix $R(t)$ by
$J^\prime_i(t) \: = \: \sum_j R(t)_i^{\: \: j} \: J_j(t)$.
It follows from the equation
\begin{equation}
\frac{d}{dt} \left( \langle J_i^\prime(t), J_j(t) \rangle
\: - \: \langle J_i(t), J^\prime_j(t) \rangle \right) \: = \: 0
\end{equation}
and the self-adjointness of $\Hess(\phi)(y)$
that $RW \: - \: W R^T \: = \: 0$ for all $t \in [0,1)$, or 
equivalently, $R \: = \:
W R^T W^{-1}$. (More intrinsically, the linear operator on
$T_{F_t(y)}M$ defined by $R$ satisfies
$R \: = \: R^*$, where $R^*$ is the
dual defined using the inner product on $T_{F_t(y)}M$.)

Next,
\begin{equation} \label{use1}
W^\prime \: = \: RW \: + \: W R^T. 
\end{equation}
Applying the Jacobi equation to
\begin{equation}
W_{ij}^{\prime \prime} \: = \: 
\langle J_i^{\prime \prime}(t), J_j(t) \rangle\: + \:
\langle J_i(t), J_j^{\prime \prime}(t) \rangle\: + \: 2 \:
\langle J_i^\prime(t), J_j^\prime(t) \rangle
\end{equation}
gives
\begin{equation} \label{use2}
W^{\prime \prime} \: = \: - \: 2 \Riem(\cdot, F_t^\prime(y), \cdot,
F_t^\prime(y)) \: + \: 2 R W R^T.
\end{equation}

Now
\begin{equation}
\frac{d}{dt} \: {\det}^{\frac{1}{2n}} W(t) \:  = \: 
\frac{1}{2n} \: {\det}^{\frac{1}{2n}} W(t) \: 
\Tr \left( W^\prime W^{-1} \right)
\end{equation}
and
\begin{align}
\frac{d^2}{dt^2} \: {\det}^{\frac{1}{2n}} W(t) \:  = \: 
& \frac{1}{4n^2} \: {\det}^{\frac{1}{2n}} W(t) \: 
\left( \Tr \left( W^\prime W^{-1} \right) \right)^2 \: - \:
\frac{1}{2n} \: {\det}^{\frac{1}{2n}} W(t) \: 
\Tr \left( (W^\prime W^{-1})^2 \right) \: + \\ 
&\frac{1}{2n} \: {\det}^{\frac{1}{2n}} W(t) \: 
\Tr \left( W^{\prime \prime} W^{-1}\right). \notag
\end{align}
Then by (\ref{use1}) and (\ref{use2}),
\begin{equation} \label{C2-1}
D^{-1} \: \frac{d^2 D}{dt^2} \: = \: \frac{1}{n^2} \: 
(\Tr(R))^2 \: - \: \frac{2}{n} \: \Tr(R^2) \: - \:
\frac{1}{n} \: \Ric(F_t^\prime(y), F_t^\prime(y)) \: + \:
\frac{1}{n} \: \Tr(R^2).
\end{equation}
As $R$ is self-adjoint,
\begin{equation}
\frac{1}{n} \: 
(\Tr(R))^2 \: - \:  \Tr(R^2) \: \le \: 0,
\end{equation}
from which the conclusion follows.
\end{proof}

\begin{proof}[Proof of Theorem~\ref{Riccimfolds}, part (a)]
To show $(1) \Rightarrow (2)$, suppose that 
$\Ric_N \: \ge \: 0$. By the definition of $\Ric_N$, we must have
$n \: < \: N$, or $n = N$ and $\Psi$ is constant. 
Suppose first that $n < N$. We can write
\begin{equation}
\Ric_N \: = \: \Ric \: - \: (N-n) \: e^{\frac{\Psi}{N-n}} \:
\Hess \left( e^{- \: \frac{\Psi}{N-n}} \right).
\end{equation}

Given $\mu_0, \mu_1 \in P^{ac}_2(M)$,
let $\{\mu_t\}_{t \in [0,1]}$ be
the unique Wasserstein geodesic from $\mu_0$ to $\mu_1$.
{From} Proposition \ref{lemsuffdc2}, in order to prove (2)
it suffices to show that
for all such $\mu_0$ and $\mu_1$, and all
$U \in \DC_N$, the inequality (\ref{defldc}) is satisfied with
$\lambda \: = \: 0$.

We recall facts from Subsection \ref{recall} about optimal transport on
Riemannian manifolds. In particular, $\mu_t$ is absolutely
continuous with respect to $\dvol_M$ for all $t$, and takes the form
$(F_t)_*\mu_0$, where $F_t(y) = \exp_y (-t\nabla\phi(y))$ for some
$\frac{d^2}{2}$-concave 
function $\phi$. Put $\eta_t \: = \: \frac{d\mu_t}{\dvol_M}$.
Using the nonsmooth change-of-variables formula proven 
in~\cite[Corollary~4.7]{Cordero-Erausquin-McCann-Schmuckenschlager (2001)}
(see also~\cite[Theorem 4.4]{McCann (1997)}), we can write
\begin{align} \label{Ueqn}
U_{\nu}(\mu_t) \: & = \: \int_M U(e^{\Psi(m)} \: 
\eta_t(m)) \: 
e^{- \: \Psi(m)} \: \dvol_M(m) \\ 
& = \:
\int_M U \left(
e^{\Psi(F_t(y))} \: \frac{\eta_0(y)}{\det(dF_t)(y)} \right) 
\:
e^{- \: \Psi(F_t(y))} \: \det(dF_t)(y) \: \dvol_M(y). \notag
\end{align}
Putting 
\begin{equation}
C(y,t) 
\: = \: e^{- \: \frac{\Psi(F_t(y))}{N}} \: {\det}^{\frac{1}{N}}(dF_t)(y),
\end{equation}
we can write
\begin{equation} \label{integrand}
U_{\nu}(\mu_t) \: 
 = \:
\int_M  C(y,t)^N \: U \left(
\eta_0(y) \: C(y,t)^{-N} \right)
\: \dvol_M(y). 
\end{equation}
Suppose that we can
show that $C(y,t)$ is concave in $t$ for almost all $y \in M$. 
Then for $y \in \supp(\mu_0)$, as the map 
\begin{equation}
\lambda \: \rightarrow \: \eta_0^{-1}(y) \: \lambda^N \: U \left(
\eta_0(y) \: \lambda^{-N} \right)
\end{equation}
is nonincreasing and convex, and the composition of a nonincreasing
convex function with a concave function is convex, 
it follows that the integrand of (\ref{integrand}) is convex in $t$.
Hence $U_{\nu}(\mu_t)$ will be convex in $t$.

To show that $C(y,t)$ is concave in $t$, fix $y$.
Put
\begin{equation}
C_1(t) \: = \: e^{- \: \frac{\Psi(F_t(y))}{N-n}}
\end{equation}
and
\begin{equation}
C_2(t) \: = \: {\det}^{\frac{1}{n}}(dF_t)(y),
\end{equation}
so
$C(y,t) \: = \: 
C_1(t)^{\frac{N-n}{N}} \: 
C_2(t)^{\frac{n}{N}}$.
We have
\begin{align} \label{C-1}
N \: C^{-1} \: \frac{d^2 C}{dt^2} \: & = \: 
(N - n) \:
C_1^{-1} \: \frac{d^2 C_1}{dt^2} \: + \:
n \:
C_2^{-1} \: \frac{d^2 C_2}{dt^2} \: - \: 
\frac{n(N-n)}{N} \: \left( C_1^{-1} \: \frac{dC_1}{dt} \: - \:
C_2^{-1} \: \frac{dC_2}{dt} \right)^2 \\
& \le \: (\Ric \: - \: \Ric_N) \left( F_t^\prime(y), F_t^\prime(y) \right)
\: + \:
n \:
C_2^{-1} \: \frac{d^2 C_2}{dt^2}. \notag
\end{align}
We may assume that the function $\phi$ has a 
Hessian at $y$
\cite[Theorem 4.2(a)]{Cordero-Erausquin-McCann-Schmuckenschlager (2001)},
and that $dF_t$ is well-defined and nonsingular at $y$ for all $t \in [0, 1)$
\cite[Claim 4.3(a-b)]{Cordero-Erausquin-McCann-Schmuckenschlager (2001)}.
Then Lemma~\ref{lemD''} shows that
\begin{equation}
n C_2^{-1}\: \frac{d^2 C_2}{dt^2} \leq -\Ric (F'_t(y),F'_t(y)).
\end{equation}
So $N C^{-1}(t)\, C''(t) \leq - \Ric_N(F'_t(y),F'_t(y))\leq 0$.
This shows that $(M, g, \nu)$ is weakly displacement convex for the
family $\DC_N$. 

The proof in the case $N = n$ follows the same lines, 
replacing $C_1$ by $1$ and $C_2$ by $C$.
\med

We now prove the implication
$(5) \Rightarrow (1)$. Putting $U \: = \: U_N$ in (\ref{integrand}), we obtain
\begin{equation} \label{special}
H_{N,\nu}(\mu_t) \: = \: N \: - \: N
\int_M C(y,t) \:\eta_0(y)^{1-\frac{1}{N}} \: \dvol_M(y).
\end{equation}
Suppose first that $n \: < \: N$ and
$H_{N, \nu}$ is weakly a.c.
displacement convex. Given $m \in M$ and $v \in T_mM$,
we want to show that $\Ric_N(v, v) \: \ge \: 0$. Choose a smooth 
function $\phi$, defined in a neighborhood of $m$, so that $v \: = \:
- \: (\nabla  \phi)(m)$, $\Hess(\phi)(m)$ is proportionate to $g(m)$ and
\begin{equation}
\frac{1}{N-n} \: v \Psi \: = \:
\frac{1}{n} \: (\triangle \phi )(m).
\end{equation}
Consider the geodesic segment $t \rightarrow \exp_m(tv)$. Then
\begin{equation}
C_1^{-1}(0) \: C_1^\prime(0) \: = \: - \: \frac{1}{N-n} \: v \Psi
\end{equation}
and
\begin{align}
C_2^{-1}(0) \: C_2^\prime(0) \: & = \: \frac{1}{2n} \: \Tr(W^\prime(0) \:
W^{-1}(0)) \: = \: \frac{1}{n} \: \Tr(R(0)) \\
& = \: - \:
\frac{1}{n} \: \Tr(\Hess(\phi)(m)) \: = \: \: - \:
\frac{1}{n} \: (\triangle \phi)(m). \notag
\end{align}
Hence by construction, $C_1^{-1}(0) \: C_1^\prime(0) \: = \:
C_2^{-1}(0) \: C_2^\prime(0)$. {From} (\ref{C-1}), it follows that
\begin{equation}
N \: C^{-1}(0) \: C^{\prime \prime}(0) \: = \:
(\Ric \: - \: \Ric_N) \left( v, v \right)
\: + \:
n \:
C_2^{-1}(0) \: C_2^{\prime \prime}(0).
\end{equation}
As $R(0)$ is a multiple of the identity, (\ref{C2-1}) now implies that 
\begin{equation}
N \: C^{-1}(0) \: C^{\prime \prime}(0) \: = \:
 - \: \Ric_N \left( v, v \right).
\end{equation}
For small numbers $\epsilon_1, \epsilon_2 > 0$, consider a smooth probability
measure $\mu_0$ with support in an $\epsilon_1$-ball around $m$.
Put $\mu_1 \: = \: (F_{\epsilon_2})_* \mu_0$ where $F_t$ is defined by
$F_t(y) \: = \: \exp_y (- \: t \: \nabla \phi(y))$. 
If $\epsilon_2$ is small enough then $\epsilon_2\phi$ is 
$\frac{d^2}{2}$-concave.
As $\mu_0$ is absolutely continuous, it follows that $F_{\epsilon_2}$ 
is the unique optimal transport between $\mu_0$ and $(F_{\epsilon_2})_*\mu_0$.
As a consequence, $\mu_t \equiv (F_{t\epsilon_2})_*\mu_0$ is the unique
Wasserstein geodesic from $\mu_0$ to $\mu_1$.
Taking
$\epsilon_1 \rightarrow 0$ and then $\epsilon_2 \rightarrow 0$,
if $H_{N, \nu}$ is to satisfy (\ref{defldc}) for
all such $\mu_0$ then we must have $C^{\prime \prime}(0) \: \le \: 0$.
Hence $\Ric_N(v, v) \: \ge \: 0$. Since $v$ was arbitrary, this shows that
$\Ric_N \: \ge \: 0$.

Now suppose that $N = n$ and 
$H_{N, \nu}$ is weakly a.c. 
displacement convex. Given $m \in M$ and $v \in T_mM$,
we want to show that $v \Psi \: = \: 0$ and
$\Ric(v, v) \: \ge \: 0$. Choose a smooth 
function $\phi$, defined in a neighborhood of $m$, so that $v \: = \:
- \: (\nabla  \phi)(m)$, and $\Hess(\phi)(m)$ is proportionate to $g(m)$.
We must again have $C^{\prime \prime}(0) \: \le \: 0$, where now
$C(t) \: = \: e^{- \: \frac{\Psi(F_t(y))}{n}} \: {\det}^{\frac{1}{n}}(dF_t)(y)$. 
By direct computation,
\begin{equation}
\frac{C^{\prime \prime}(0)}{C(0)} \: = \:
- \: \frac{1}{n} \: (\Ric \: + \: \Hess(\Psi))(v, v)
\: + \: \frac{(v \Psi)^2}{n^2} \: + \: \frac{2 \: (v \Psi) \:
(\triangle \phi)(m)}{n^2}.  
\end{equation}
If $v \Psi \: \neq \: 0$ then we can make $C^{\prime \prime}(0) \: > \: 0$
by an appropriate choice of $\triangle \phi$. Hence $\Psi$ must be constant
and then we must have $\Ric(v, v) \: \ge \: 0$.

Finally, if $N < n$ then (\ref{C-1}) gives
\begin{equation}
N \: \frac{C^{\prime \prime}(0)}{C(0)} \: = \:
- \: (\Ric \: + \: \Hess(\Psi))(v, v)
\: + \: \frac{(v \Psi)^2}{N-n} \: - \: 
\frac{n(N-n)}{N} \left( - \: \frac{v \Psi}{N-n} \: + \: 
\frac{(\triangle \phi)(m)}{n} \right)^2.
\end{equation}
One can always choose $(\triangle \phi)(m)$ to make 
$C^{\prime \prime}(0)$ positive, so $H_{N, \nu}$ cannot be 
weakly a.c. displacement convex.
\end{proof}

\begin{proof}[Proof of Theorem~\ref{Riccimfolds}, part (b)]
We first show $(1) \Rightarrow (2)$.
Suppose that $\Ric_\infty \: \ge \: K \: g$. 
Given $\mu_0, \mu_1 \in P^{ac}_2(M)$, 
we again use (\ref{Ueqn}), with  $U \in \DC_\infty$.
Putting 
\begin{equation}
C(y,t) \: = \: - \: \Psi(F_t(y)) \: \: + \: \log {\det}(dF_t)(y),
\end{equation}
we have
\begin{equation} 
U_{\nu}(\mu_t) \:  = \:
\int_M  e^{C(y,t)} \: U \left(
\eta_0(y) \: e^{-C(y,t)} \right) \: \dvol_M(y). 
\end{equation}
As in the proof of (a), the condition $\Ric_\infty \: \ge \:
K \: g$ implies that 
\begin{equation}
\frac{d^2C}{dt^2} \: \le \: - \: K \: 
|F_t^\prime(y)|^2
\: = \: - \: K \: 
|\nabla \phi|^2(y),
\end{equation}
where the last equality comes from the constant speed of
the geodesic $t \rightarrow F_t(y)$.
By assumption, the map 
\begin{equation}
\lambda \rightarrow \eta_0^{-1}(y) \: e^\lambda \:
U(\eta_0(y) \: e^{-\lambda})
\end{equation}
is nonincreasing and convex in $\lambda$, with derivative
$- \: \frac{p(\eta_0(y) \: e^{-\lambda})}{\eta_0(y) \: e^{-\lambda}}$.
It follows that the composition
\begin{equation}
t \rightarrow \eta_0^{-1}(y) \: e^{C(y,t)} \: U \left(
\eta_0(y) \: e^{-C(y,t)} \right)
\end{equation}
is $\lambda(U) \: |\nabla \phi|^2(y)$-convex in $t$.
Then
\begin{align}
e^{C(y,t)} \: U \left(
\eta_0(y) \: e^{-C(y,t)} \right)
\: \le \: & t \: 
e^{C(y,1)} \: U \left(
\eta_0(y) \: e^{-C(y,1)} \right) \: + \\
&  (1-t) \: 
e^{C(y,0)} \: U \left(
\eta_0(y) \: e^{-C(y,0)} \right)
\: - \notag \\
&  \frac12 \: \lambda(U) \: |\nabla \phi|^2(y) \: \eta_0(y) \: t(1-t). \notag
\end{align}
Integrating with respect to $\dvol_M(y)$ and using the fact that
\begin{equation}
W_2(\mu_0, \mu_1)^2 \: = \: 
\int_M |\nabla \phi|^2(y) \: \eta_0(y) \: \dvol_M(y)
\end{equation}
shows that (\ref{defldc}) is satisfied with 
$\lambda \: = \: \lambda(U)$.
The implication $(1) \Rightarrow (2)$ now follows from
Proposition \ref{lemsuffdc2}.

The proof that $(5) \Rightarrow (1)$ is similar to the proof in part (a).
\end{proof}

The case $N=1$ is slightly different  because $H_{1,\nu}$ is not defined.
However, the rest of Theorem \ref{Riccimfolds}.a carries through.

\begin{theorem} \label{Riccimfolds1}
a. The following are equivalent : \\
(1) $\Ric_1 \: \ge \: 0$. \\
(2) The measured length space
$(M, g, \nu)$ has nonnegative $1$-Ricci curvature. \\
(3) For all $U \in \DC_1$, $U_\nu$ is weakly displacement convex on 
$P_2(M)$.\\
(4) For all $U \in \DC_1$, $U_\nu$ is weakly a.c. displacement convex on 
$P^{ac}_2(M)$.
\end{theorem}
\begin{proof}
The proofs of $(1) \Rightarrow (2) \Rightarrow (3) \Rightarrow (4)$ are as
in the proof of Theorem \ref{Riccimfolds}.a. It remains to show that 
$(4) \Rightarrow (1)$. 
Since $\DC_N\subset \DC_1$ for all $N>1$,
condition (4) implies that $U_\nu$ is weakly a.c. displacement convex
on $P^{ac}_2(M)$ for all $U\in \DC_N$. So by Theorem~\ref{Riccimfolds}.a,
$M$ satisfies $\Ric_N\:\ge \: 0$ for all $N>1$. It follows that
$n\: \le\: 1$. If $n = 0$, i.e. $M$ is a point, then
$\Ric_1 \:\ge \: 0$ holds automatically. If $n =1$,
i.e. $M$ is a circle, then
taking $N \rightarrow 1^+$ shows that $\Ric_1 \: \ge \: 0$, i.e.
$\Psi$ is constant. 
\end{proof}

\begin{remark} \label{acequiv}
In the Riemannian case there is a unique
Wasserstein geodesic joining $\mu_0, \mu_1 \in P^{ac}_2(M)$.
Hence we could add two more equivalences to
Theorem \ref{Riccimfolds}. Namely, a.(4) 
is equivalent to saying that for all 
$U \in \DC_N$, $U_\nu$ is a.c. displacement convex on $P^{ac}_2(M)$,
and b.(4) is equivalent to saying that for all 
$U \in \DC_\infty$, $U_\nu$ is $\lambda(U)$-a.c. displacement convex on 
$P^{ac}_2(M)$.
\end{remark}

\begin{remark}
Theorem \ref{Riccimfolds} also holds under weaker regularity assumptions.
For example, if $\Psi$ is a continuous function on Euclidean $\R^n$ then
$\left( \R^n, \: e^{- \Psi} \: \dvol_{\R^n} \right)$ has 
$\infty$-Ricci curvature bounded below by zero if and only if $\Psi$ is convex.
\end{remark}

\subsection{Geometric corollaries}

We have shown that our abstract notion of a lower Ricci curvature
bound is stable under measured Gromov--Hausdorff convergence
(Theorem~\ref{stabRicci}) and, in the Riemannian setting, coincides with
a classical notion (Theorem~\ref{Riccimfolds}).
This subsection is devoted to various geometric applications.

We first give a characterization of the {\em smooth} elements in the set of
measured Gromov--Hausdorff limits of manifolds with Ricci curvature
bounded below.

\begin{corollary} \label{limitspacecor}
Let $(B,g_B)$ be a smooth compact connected
Riemannian manifold, equipped with the
Riemannian density
$\dvol_B$, and let $\Psi$ be a $C^2$-regular function on $B$
which is normalized 
by an additive constant so that
$e^{-\Psi}\dvol_B$ is a probability measure on $B$.
We have the following implications: \\
(i) If $(B, g_B, e^{-  \Psi} \dvol_B)$ is a measured 
Gromov--Hausdorff limit of Riemannian manifolds with nonnegative
Ricci curvature and dimension at most $N$ then $\Ric_N(B) \: \ge \:0$. \\
(i') If $(B, g_B, e^{-  \Psi} \dvol_B)$ is a measured 
Gromov--Hausdorff limit of Riemannian manifolds with Ricci curvature
bounded below by $K \in \R$ then $\Ric_\infty(B) \: \ge \: K \: g_B$. \\
(ii) As a partial converse, if $(B, g_B, e^{-  \Psi} \dvol_B)$ 
has $\Ric_N(B) \: \ge \: 0$ with $N \: \ge \: \dim(B) \: + \: 2$ then 
$(B,  g_B, e^{-  \Psi} \dvol_B)$ is a measured 
Gromov--Hausdorff limit of Riemannian manifolds with nonnegative
Ricci curvature and dimension at most $N$. \\
(ii') If $(B, g_B, e^{-  \Psi} \dvol_B)$ has $\Ric_\infty(B) \: \ge \: 
K \: g_B$
then $(B, g_B, e^{-  \Psi} \dvol_B)$ is a measured 
Gromov--Hausdorff limit of Riemannian manifolds $M_i$ with 
$\Ric(M_i) \: \ge \: (K \:- \: \frac{1}{i}) \: g_{M_i}$.
\end{corollary}

\begin{proof}
Parts (i) and (i') are a direct consequence of Theorems \ref{stabRicci},
\ref{Riccimfolds} and \ref{Riccimfolds1}. Part (ii) follows from the warped 
product construction 
of~\cite[Theorem 3.1]{Lott (2003)}. The proof of (ii') is similar.
\end{proof}

\begin{remark} In Corollary~\ref{limitspacecor}(ii'), if $K\neq 0$ then
one can use a rescaling argument to transform the condition
$\Ric(M_i)\, \ge\, \left( K- \frac{1}{i} \right)\, 
g_{M_i}$ into the more stringent condition
$\Ric(M_i)\,\ge\, K\,g_{M_i}$.
\end{remark}

The next two corollaries give some consequences of Corollary 
\ref{limitspacecor} for the {\em metric} structure of smooth limit spaces,
i.e. for the aspects of the limit metric-measure spaces that are
independent of the measure.
In general, one cannot change the conclusion of
Corollary \ref{limitspacecor}(i) to obtain a lower bound on $\Ric$ instead
of $\Ric_N$. However, one does obtain such a lower bound in
the noncollapsing case.

\begin{corollary} \label{noncollapse}
a. Suppose that $(X, d)$ 
is a Gromov--Hausdorff limit of $n$-dimensional Riemannian 
manifolds with nonnegative
Ricci curvature. If $(X, d)$ has Hausdorff dimension $n$, and
$\nu_H$ is its normalized $n$-dimensional Hausdorff measure, then
$(X, d, \nu_H)$ has nonnegative $n$-Ricci curvature. \\
b. If in addition $(X, d)$ happens to be a smooth $n$-dimensional
Riemannian manifold $(B, g_B)$ then $\Ric(B) \: \ge \:0$.
\end{corollary}

\begin{proof}
a. If $\{M_i\}_{i=1}^\infty$ is a sequence of $n$-dimensional Riemannian
manifolds with nonnegative Ricci curvature 
and $\{f_i\}_{i=1}^\infty$ is a sequence
of $\epsilon_i$-approximations $f_i \: : \: M_i \rightarrow X$, with
$\lim_{i \rightarrow \infty} \epsilon_i \: = \: 0$, then
$\lim_{i \rightarrow \infty} (f_i)_* \dvol_{M_i} \: = \: \nu_H$ in the
weak-$*$ topology
\cite[Theorem 5.9]{Cheeger-Colding (1997)}. (This also shows
that the $n$-dimensional Hausdorff measure on $X$ can be normalized
to be a probability measure.) Then part a. follows from
Theorems \ref{stabRicci} and \ref{Riccimfolds}.\\
b. If $(X, d) \: = \: (B, g_B)$ then $\nu_H \: = \: \frac{\dvol_B}{\vol(B)}$
and the claim follows from Theorem \ref{Riccimfolds}, along with
the definition of $\Ric_n$.
\end{proof}

\begin{remark}
A special case of Corollary \ref{noncollapse}.a is when $(X, d)$ is an
$n$-dimensional Gromov--Hausdorff limit of a sequence of $n$-dimensional
Riemannian manifolds with nonnegative sectional curvature.  In this
case, $(X, d)$ has nonnegative Alexandrov curvature and
$(X, d, \nu_H)$ has nonnegative $n$-Ricci curvature.  More generally,
we expect that for an
$n$-dimensional compact length space $(X, d)$ with Alexandrov curvature
bounded below, equipped with the normalized $n$-dimensional
Hausdorff measure $\nu_H$, \\
1. If $(X, d)$ has nonnegative Alexandrov
curvature then $(X, d, \nu_H)$ has nonnegative $n$-Ricci curvature, and \\
2. For $n > 1$, 
if $(X, d)$ has Alexandrov curvature bounded below by $\frac{K}{n-1}$
then $(X, d, \nu_H)$ has $\infty$-Ricci curvature bounded below by $K$.

It is possible that
the proof of Theorem \ref{Riccimfolds} can be adapted to show this.
\end{remark}

As mentioned above, in the collapsing case the lower bound in the
conclusion of
Corollary \ref{limitspacecor}(i) (or
Corollary \ref{limitspacecor}(i')) would generally fail if we replaced
$\Ric_N$ (or $\Ric_\infty$) by $\Ric$.
However, one does obtain a lower
bound on the average scalar curvature of $B$.

\begin{corollary} \label{meanscalar}
If $\left( B, g_B, e^{-\Psi} \: \dvol_B \right)$ is a smooth $n$-dimensional
measured Gromov--Hausdorff limit of Riemannian manifolds 
(of arbitrary dimension), each with Ricci curvature bounded below by 
$K \in \R$, then the scalar curvature $S$ of $(B, g_B)$ satisfies 
\begin{equation}
\frac{\int_B S \: \dvol_B}{\vol(B)} 
\: \ge \: nK.
\end{equation}
\end{corollary}

\begin{proof}
From Corollary \ref{limitspacecor}(iii), $\Ric(B) \: + \: \Hess(\Psi) \:
\ge \: K \: g_B$. Tracing gives $S \: + \: \triangle \Psi \: \ge nK$.
Integrating gives $\int_B S \: \dvol_B \: \ge \: nK \: \vol(B)$. 
\end{proof}

\med

Next, we show
that for Riemannian manifolds, lower $N$-Ricci curvature bounds are
preserved upon taking quotients by compact Lie group actions.

\begin{corollary} \label{Gquotient}
Let $M$ be a compact connected Riemannian manifold.
Let $G$ be a compact Lie group that acts
isometrically on $M$, preserving a function $\Psi \in C^\infty(M)$
that satisfies $\int_M e^{-\Psi} \: \dvol_M \: = \: 1$. 
Let $p \: : \: M \rightarrow M/G$
be the quotient map. \\
a. For $N \in [1, \infty)$, if $\left( M, \: e^{-\Psi} \: \dvol_M
\right)$ has 
$\Ric_N \: \ge \:  0$
then $(M/G, d_{M/G}, p_* (e^{-\Psi} \: \dvol_M))$ 
has nonnegative $N$-Ricci curvature. \\
b. If $\left( M, \: e^{-\Psi} \: \dvol_M
\right)$ has $\Ric_\infty \: \ge \: K \: g_M$ then
$(M/G, d_{M/G}, p_* (e^{-\Psi} \: \dvol_M))$ 
has $\infty$-Ricci curvature bounded below by $K$.
\end{corollary}
\begin{proof}
This follows from Theorem \ref{Gaction}.
\end{proof}

Corollary \ref{Gquotient} provides many examples of singular spaces
with lower Ricci curvature bounds.  Of course, the main case is when
$\Psi$ is constant.

\med

We conclude this section by giving a ``synthetic'' proof of a part 
of the Ricci O'Neill theorem of \cite[Theorem 2]{Lott (2003)}.

\begin{corollary} \label{Oneill}
Let $p \: : \: M \rightarrow B$ be a Riemannian submersion of compact
connected manifolds, with fibers $Z_b$. Choose $N\geq \dim(M)$ and
$\Psi_M \in C^\infty(M)$ with 
$\int_M e^{-\Psi_M} \: \dvol_M \: = \: 1$; if $N=\dim(M)$ then we assume that
$\Psi_M$ is constant. Define $\Psi_B \in C^\infty(B)$
by $p_* \left( e^{-\Psi_M} \: \dvol_M \right) \: = \: 
e^{-\Psi_B} \: \dvol_B$. Suppose that the fiber parallel
transport of the Riemannian submersion preserves the fiberwise measures
$e^{-\Psi_M} \big|_Z \: \dvol_Z$ up to multiplicative constants. (That is,
if $\gamma \: : \: [0, 1] \rightarrow B$ is a smooth path in $B$, let
$P_\gamma \: : \: Z_{\gamma(0)} \rightarrow Z_{\gamma(1)}$ denote the
fiber transport diffeomorphism. Then we assume that there is a constant 
$C_\gamma > 0$ so that 
\begin{equation}
\left. P_\gamma^* \left(
e^{-\Psi_M} \big|_{Z_{\gamma(1)}} \: \dvol_{Z_{\gamma(1)}} \right) \: = \:
C_\gamma \: e^{-\Psi_M} \big|_{Z_{\gamma(0)}} \: \dvol_{Z_{\gamma(0)}}. \qquad \right)
\end{equation}
With these assumptions, \\
a. If $\Ric_N(M) \ge 0$ then $\Ric_N(B) \ge 0$. \\
b. For any $K \in \R$, if $\Ric_\infty(M) \ge K g_M$ then
$\Ric_\infty(B) \ge K g_B$.
\end{corollary}

\begin{proof}
Put $\nu_M \: = \: e^{- \Psi_M} \: \dvol_M$ and 
$\nu_B \: = \: e^{- \Psi_B} \: \dvol_B$. We can decompose $\nu_M$
with respect to $p$ as $\sigma(b) \: \nu_B(b)$, with $\sigma(b) \in
P^{ac}_2 \left( Z_b \right)$. 
{From} the assumptions, the family $\{\sigma(b)\}_{b \in B}$ of
vertical densities is invariant under fiber parallel transport.

To prove part (a), let $\{\mu_t\}_{t \in [0,1]}$ be a 
Wasserstein geodesic in $P^{ac}_2(B)$. Define 
$\{\mu^\prime_t\}_{t \in [0,1]}$ in $P^{ac}_2(M)$ by
$\mu^\prime_t \: \equiv \: \sigma(b) \: \mu_t(b)$. By construction,
the corresponding densities satisfy $\rho_t^\prime \: = \: p^* \rho_t$.
Thus $H_{N, \nu_M}(\mu^\prime_t) \: =\: H_{N, \nu_B}(\mu_t)$.
Furthermore, $\{\mu^\prime_t\}_{t \in [0,1]}$ is a Wasserstein geodesic; 
if $(F_t)_{t \in [0,1]}$ is an optimal 
Monge transport from $\mu_0$ to $\mu_1$ then its horizontal lift is
an optimal Monge transport from $\mu^\prime_0$ to $\mu^\prime_1$, with 
generating
function $\phi_M \: = \: p^* \phi_B$. {From}
Theorem~\ref{Riccimfolds}(a) and Remark \ref{acequiv}, 
$H_{N, \nu_M}$ is a.c. displacement convex on $P^{ac}_2(M)$.
In particular, (\ref{defldc}) is satisfied along
$\{\mu^\prime_t\}_{t \in [0,1]}$ with $U_\nu \: = \:  
H_{N,\nu_M}$ and $\lambda \: = \: 0$.
Then the same equation is satisfied along
$\{\mu_t\}_{t \in [0,1]}$ with $U_\nu \: = \:  
H_{N,\nu_B}$ and $\lambda \: = \: 0$. Thus $H_{N, \nu_B}$ is
a.c. displacement convex on $P^{ac}_2(B)$.
Theorem~\ref{Riccimfolds}(a) now implies that $\Ric_N(B) \: \ge \:0$.

The proof of part (b) is similar.
\end{proof}

\begin{remark} In fact, for any $N \in [1, \infty]$ and any $K \in \R$,
if $\Ric_N(M) \ge K g_M$ then $\Ric_N(B) \ge K g_B$. This was proven in
\cite[Theorem 2]{Lott (2003)} in the cases $N = \infty$ and $N = \dim(M)$ by
explicit tensor calculations. (The paper \cite{Lott (2003)} writes
$\Ric_q$ for what we write as $\Ric_N$, where $q \: = \: N-n$.)
The same method of proof works for all $N$.
\end{remark}

\begin{remark}
Suppose that $M$ is a compact connected
Riemannian manifold on which a compact Lie group 
$G$ acts isometrically, with all orbits of the same orbit type.  Put $B \: = \:
M/G$.
If $\Psi_M \in C^\infty(M)$ is a $G$-invariant function that
satisfies
$\int_M e^{-\Psi_M} \: \dvol_M \: = \: 1$, and
$\left( M, g_M, e^{-\Psi_M} \: \dvol_M \right)$ has $\Ric_N(M) \: \ge \:0$,
then Corollaries \ref{Gquotient} and \ref{Oneill} overlap in saying that
$\left( B, g_B,  e^{-\Psi_B} \: \dvol_B \right)$ has $\Ric_N(B) \: \ge \:0$.
There is a similar statement when
$\Ric_\infty(M) \: \ge \: K \: g_M$.
\end{remark}

\begin{remark} There is an obvious analogy between the 
Ricci O'Neill theorem and O'Neill's theorem that sectional curvature 
is nondecreasing under pushforward by a Riemannian submersion.
There is also a ``synthetic'' proof of O'Neill's theorem, obtained
by horizontally lifting a geodesic hinge from $B$ and using
triangle comparison results, along with the fact that $p$ is
distance nonincreasing.
\end{remark}

\appendix

\section{The Wasserstein space as an Alexandrov space} \label{secgeom}

This section is concerned with the geometry of the Wasserstein space
$P_2(M)$ of a Riemannian manifold $M$.
Otto introduced a formal infinite-dimensional Riemannian metric on
$P_2(M)$ and showed that $P_2(\R^n)$ formally has nonnegative sectional
curvature \cite{Otto (2001)}.  We make such results rigorous by
looking at $P_2(M)$ as an Alexandrov space.

We first give a general lower bound on Wasserstein distances in terms
of Lipschitz functions. We show that if $M$ is a compact Riemannian
manifold with nonnegative
sectional curvature then $P_2(M)$ has nonnegative Alexandrov curvature.
Using the above-mentioned lower bound on Wasserstein distances, we
compute the tangent cones of $P_2(M)$ at the absolutely continuous
measures.

\subsection{Lipschitz functions and optimal transport}
In general, one can estimate Wasserstein distances from above
by choosing particular transference plans. The next lemma provides
a way to estimate these distances
from below by using Lipschitz functions.

\begin{lemma} \label{lemlip}
If $X$ is a compact length space and
$\{ \mu_t \}_{t \in [0,1]}$ is a Wasserstein geodesic then
for all $f \in \Lip(X)$,
\begin{equation}
\left | \int_X f d\mu_1 \: - \: \int_X f d\mu_0 \right |^2 \: \le \:
W_2(\mu_0, \mu_1)^2 \: 
\int_0^1 \left (\int_X |\nabla f|^2 \: d\mu_t\right) \: dt.
\end{equation}
\end{lemma}

\begin{proof}
By Proposition \ref{everygeod}, the Wasserstein geodesic arises as the
displacement interpolation associated to some optimal dynamical
transference plan $\Pi$.
We have
\begin{align}
\int_X f \: d\mu_1 \: - \: \int_X f \: d\mu_0 \: & = \:
\int_X f \: d(e_1)_* \Pi \: - \: \int_X f \: d(e_0)_* \Pi \: = \:
\int_\Gamma ((e_1)^* f  \: - \: (e_0)^* f) \: d \Pi \\
& = \:
\int_\Gamma (f(\gamma(1))  \: - \: f(\gamma(0))) \: d \Pi(\gamma)\notag.
\end{align}
As $f \circ \gamma \in \Lip([0,1])$, we have
\begin{equation}
|f(\gamma(1))  \: - \: f(\gamma(0))| \: =\:
\left| \int_0^1 \frac{df(\gamma(t))}{dt} \: dt \right| \: \le \:
\int_0^1 \left| \frac{df(\gamma(t))}{dt} \right| \: dt \: \le \:
\int_0^1 |\nabla f|(\gamma(t)) \: L(\gamma) \:  dt.
\end{equation}
Then
\begin{equation}
\left|
\int_X f \: d\mu_1 \: - \: \int_X f \: d\mu_0 \right| \: \le \:
\int_\Gamma \int_0^1 |\nabla f|(\gamma(t)) \: L(\gamma) \: dt \: d\Pi(\gamma).
\end{equation}
{From} the Cauchy--Schwarz inequality,
\begin{align}
\left| \int_X f \: d\mu_1 \: - \: \int_X f \: d\mu_0 \right|^2 \: & \le \:
\int_\Gamma \int_0^1 |\nabla f|^2(\gamma(t)) \:  dt \: d\Pi(\gamma)
\int_\Gamma \int_0^1 L(\gamma)^2 \: dt \: d\Pi(\gamma) \\
& = \: \int_0^1 \int_\Gamma |\nabla f|^2(\gamma(t)) \:  dt \: d\Pi(\gamma)
\int_\Gamma L(\gamma)^2 \: d\Pi(\gamma). \notag
\end{align}
We have $W_2(\mu_0, \mu_1)^2 \: = \: \int_\Gamma L(\gamma)^2 \: d\Pi(\gamma)$.
To conclude the proof, we note that
\begin{equation}
\int_\Gamma |\nabla f|^2(\gamma(t)) \: d\Pi(\gamma) \: = \:
\int_\Gamma |\nabla f|^2(e_t(\gamma)) \: d\Pi(\gamma) \: = \:
\int_X |\nabla f|^2 \: d(e_t)_* \Pi \: = \:
\int_X |\nabla f|^2 \: d\mu_t.
\end{equation}
\end{proof}

\subsection{The case of nonnegatively curved manifolds} \label{sectAlex}

\begin{theorem} \label{thmalex}
A smooth compact connected manifold $M$ 
has nonnegative sectional curvature if and only if $P_2(M)$ has
nonnegative Alexandrov curvature.
\end{theorem}

\begin{proof}
Suppose that $M$ has nonnegative sectional curvature.  We first show
that $P^{ac}_2(M)$ has
nonnegative Alexandrov curvature. Let 
$\mu_0$, $\mu_1$, $\mu_2$ and $\mu_3$ be points in $P^{ac}_2(M)$ with
$\mu_i \neq \mu_0$ for $1 \le i \le 3$. For $1 \le i \neq j \le 3$,
let $\widetilde{\angle} \mu_i \mu_0 \mu_j$ denote the comparison
angle at $\mu_0$ of the triangle formed by
$\mu_i$, $\mu_0$ and $\mu_j$
\cite[Definition 3.6.25]{Burago-Burago-Ivanov (2001)}.
For $1 \le i \le 3$, let $c_i \: : \: [0,1] \rightarrow P^{ac}_2(M)$ 
be a Wasserstein geodesic from
$\mu_0$ to $\mu_i$. Let $\phi_i$ be the corresponding
$\frac{d^2}{2}$-concave function on $M$ that generates $c_i$.
That is, $c_i(t) \: = \: (F_{i,t})_* \mu_0$,
with $F_{i,t}(m) \: = \: \exp_m(- \: t \: \nabla \phi_i(m))$.
Then $W_2(\mu_0, \mu_i)^2 \: = \: 
\int_M d(m, F_{i,t}(m))^2 \: d\mu_0(m) \: = \: 
\int_M |\nabla \phi_i|^2 \: d\mu_0$.

Let us define a particular transference plan from $\mu_i$ to
$\mu_j$ as a Monge transport $F_{j,1} \circ F_{i,1}^{-1}$. 
It gives an upper bound on $W_2(\mu_i, \mu_j)^2$ by
\begin{equation}
W_2(\mu_i, \mu_j)^2 \: \le \: 
\int_M d(F_{i,1}(m), F_{j,1}(m))^2 \: d\mu_0(m).
\end{equation}
For almost all $m$,
the nonnegative curvature of $M$, applied to the hinge at $m$ formed
by the geodesic segments
$t \rightarrow F_{i,t}(m)$ and $t \rightarrow F_{j,t}(m)$, implies
\begin{equation} \label{ineqgeodseg}
d(F_{i,1}(m), F_{j,1}(m))^2 \: \le \:
|\nabla \phi_i(m)|^2 \: + \: |\nabla \phi_j(m)|^2 \: - \:
2 \: \langle \nabla \phi_i(m), \nabla \phi_j(m)  \rangle.
\end{equation}
Integrating~\eqref{ineqgeodseg} with respect to $\mu_0$ yields
\begin{equation}
W_2(\mu_i, \mu_j)^2 \: \le \: W_2(\mu_0, \mu_i)^2 \: + \: 
W_2(\mu_0, \mu_j)^2 \: - \: 2 \:
\int_M  \: \langle \nabla \phi_i(m), \nabla \phi_j(m)  \rangle \:
d\mu_0(m).
\end{equation}
Thus $\widetilde{\angle} \mu_i \mu_0 \mu_j \: \le \: \theta_{ij}$, where
$\theta_{ij} \in [0, \pi]$ is defined by
\begin{equation} \label{defangle}
\cos \theta_{ij} \: = \: 
\frac{\int_M \langle \nabla \phi_i, \nabla \phi_j \rangle \: 
d\mu_0}{\sqrt{\int_M | \nabla \phi_i |^2 \: 
d\mu_0 \: \int_M | \nabla \phi_j |^2 \: 
d\mu_0}}.
\end{equation}
It follows from the geometry of an inner product space that
\begin{equation}
\theta_{12} \: + \: \theta_{23} \: + \: \theta_{31} \: 
\le \: 2\pi.
\end{equation}
Thus
\begin{equation} \label{sumangles}
\widetilde{\angle} \mu_1 \mu_0 \mu_2 \: + \: 
\widetilde{\angle} \mu_2 \mu_0 \mu_3 \: + \: 
\widetilde{\angle} \mu_3 \mu_0 \mu_1 \: \le \: 2\pi,
\end{equation}
which implies
that $P^{ac}_2(M)$ has nonnegative Alexandrov curvature
\cite[Proposition 10.1.1]{Burago-Burago-Ivanov (2001)}.

As $P_2(M)$ is the completion of $P^{ac}_2(M)$, the fact that
(\ref{sumangles}) can be written solely in terms of distances implies that it
also holds for $P_2(M)$. Thus $P_2(M)$ has nonnegative Alexandrov curvature.

Now suppose that $P_2(M)$ has nonnegative Alexandrov curvature. As
the embedding $M \rightarrow P_2(M)$ by delta functions defines 
a totally geodesic subspace of $P_2(M)$,
it follows that $M$ has nonnegative Alexandrov curvature.  Thus
$M$ has nonnegative sectional curvature.
\end{proof}

\begin{remark} \label{existgradflow}
The fact that $P_2(M)$ has nonnegative Alexandrov curvature
ensures the existence and uniqueness of the gradient flow of a 
$\lambda$-concave function on $P_2(M)$ \cite[Appendix]{Perelman-Petrunin}.
(The conventions of \cite{Perelman-Petrunin} are such that
the function increases along the flowlines of its gradient flow; 
some other authors have the convention
that a function decreases along the flowlines of its gradient flow,
and hence consider $\lambda$-convex functions.) Other approaches to
geometrizing $P_2(M)$, with a view toward defining gradient flows,
appear in \cite{Ambrosio-Gigli-Savare 
(2004),Ambrosio-Gigli-Savare,Carrillo-McCann-Villani (2004)}.
\end{remark}

Now suppose that $M$ has nonnegative sectional curvature. Let
$c_0, c_1 \: : \: [0, 1] \rightarrow P_2(M)$
be nontrivial Wasserstein geodesics, with
$c_0(0) \: = \: c_1(0) \: = \: \mu$.
Theorem~\ref{thmalex} implies that 
the comparison
angle $\widetilde{\angle}c_0(s_0) \mu c_1(s_1)$ is monotonically
nonincreasing as $s_0$ and $s_1$ increase, separately in $s_0$ and $s_1$
\cite[Definition 4.3.1, Theorem 4.3.5 and 
Theorem 10.1.1]{Burago-Burago-Ivanov (2001)}.
Then
there is a well-defined angle $\angle (c_0, c_1)$ that
$c_0$ and $c_1$ form at $\mu$, in the sense of
\cite[Definition 3.6.26]{Burago-Burago-Ivanov (2001)}, given by
\begin{equation}
\angle (c_0, c_1) \: = \: \lim_{s_0, s_1 \rightarrow 0^+}
\widetilde{\angle}c_0(s_0) \mu c_1(s_1).
\end{equation}

\begin{proposition} \label{anglecomp}
Let
$c_0, c_1 \: : \: [0, 1] \rightarrow P_2^{ac}(M)$
be nontrivial Wasserstein geodesics, with
$c_0(0) \: = \: c_1(0) \: = \: \mu$.
If $\phi_0$ and $\phi_1$ are the $\frac{d^2}{2}$-concave functions 
that generate $c_0$ and $c_1$, respectively, then
\begin{equation}\label{anglec0c1}
\cos \angle (c_0, c_1) \: = \: 
\frac{\int_M \langle \nabla \phi_0, \nabla \phi_1 \rangle \: 
d\mu}{\sqrt{\int_M | \nabla \phi_0 |^2 \: 
d\mu \: \int_M | \nabla \phi_1 |^2 \: 
d\mu}},
\end{equation}
\end{proposition}

\begin{proof}
Applying~\eqref{defangle} (and the sentence preceding it) 
to the triangle $\triangle c_0(s_0) \mu c_1(s_1)$
gives
\begin{equation}
\cos \widetilde{\angle}c_0(s_0) \mu c_1(s_1) \: \ge \: 
\frac{\int_M \langle \nabla \phi_0, \nabla \phi_1 \rangle \: 
d\mu}{\sqrt{\int_M | \nabla \phi_0 |^2 \: 
d\mu \: \int_M | \nabla \phi_1 |^2 \: 
d\mu}},
\end{equation}
and so
\begin{equation}
\cos \angle (c_0, c_1) \: \ge \: 
\frac{\int_M \langle \nabla \phi_0, \nabla \phi_1 \rangle \: 
d\mu}{\sqrt{\int_M | \nabla \phi_0 |^2 \: 
d\mu \: \int_M | \nabla \phi_1 |^2 \: 
d\mu}}.
\end{equation}
{From} the monotonicity of $\widetilde{\angle}c_0(s_0) \mu c_1(s_1)$,
it suffices to show that
\begin{equation}
\lim_{s \rightarrow 0} \:
\cos \widetilde{\angle} c_0(s) \mu c_1(s) \: \le \: 
\frac{\int_M \langle \nabla \phi_0, \nabla \phi_1 \rangle \: 
d\mu}{\sqrt{\int_M | \nabla \phi_0 |^2 \: 
d\mu \: \int_M | \nabla \phi_1 |^2 \: 
d\mu}}.
\end{equation}
This amounts to showing a lower bound on $W_2(c_0(s), c_1(s))$. 

Let $\{\mu_{t,s}\}_{t \in [0,1]}$ be a Wasserstein geodesic from
$c_0(s)$ to $c_1(s)$.
{From} Lemma~\ref{lemlip}, for any $f \in C^1(M)$,
\begin{equation}
\left( \int_M f dc_1(s) \: - \: \int_M f dc_0(s) \right)^2 \: \le \:
W_2(c_0(s), c_1(s))^2 \: \int_0^1 \int_M |\nabla f|^2 \: d\mu_{t,s} \: dt.
\end{equation}
In terms of the Monge transport maps $F_{0,t}$ and $F_{1,t}$,
\begin{equation}
\int_M f \: dc_1(s) \: - \: \int_M f \: dc_0(s) \: = \:
\int_M f \: d(F_{1,s})_* \mu \: - \: \int_M f \: d(F_{0,s})_* \mu \: = \:
\int_M ((F_{1,s})^* f \: - \: (F_{0,s})^* f) \: d \mu.
\end{equation}
Thus
\begin{equation} \label{eqnref}
\left( \frac{\int_M ((F_{1,s})^* f \: - 
\: (F_{0,s})^* f) \: d \mu}{s} \right)^2 
\: \le \:
\frac{W_2(c_0(s), c_1(s))^2}{s^2} \: 
\int_0^1 \int_M |\nabla f|^2 \: d\mu_{t,s} \: dt.
\end{equation}

Since $\{\mu_{t,s}\}_{t \in [0,1]}$ is minimizing between its
endpoints, we must have
\begin{equation}
W_2(\mu, \mu_{t,s}) \: \le \: W_2(\mu, c_0(s)) \: + \:
W_2(\mu, c_1(s))
\end{equation}
for all $t \in [0,1]$. 
(Otherwise the length of $\{\mu_{t,s}\}_{t \in [0,1]}$
would have to be greater than $W_2(\mu, c_0(s)) \: + \:
W_2(\mu, c_1(s))$. Then 
there would be a
path from $c_0(s)$ to $c_1(s)$ that is shorter 
than $\{\mu_{t,s}\}_{t \in [0,1]}$,
 obtained by going from $c_0(s)$ to $\mu$ along $c_0$ and then
from $\mu$ to $c_1(s)$ along $c_1$.)
If $\pi_{t,s}$ is an optimal transference plan between
$\mu$ and $\mu_{t,s}$ then
$\int_0^1 \pi_{t,s} \: dt$ is a transference plan between
$\mu$ and $\int_0^1 \mu_{t,s} \: dt$, showing that
\begin{equation}
W_2 \left( \mu, \int_0^1 \mu_{t,s} \: dt \right)^2 \: \le \:
\int_0^1 W_2 \left( \mu, \mu_{t,s} \right)^2 \: dt \: \le \:
\left( W_2(\mu, c_0(s)) \: + \:
W_2(\mu, c_1(s)) \right)^2.
\end{equation}

Thus
$\lim_{s \rightarrow 0} \int_0^1 \mu_{t,s} \: dt \: = \: \mu$
in the weak-$*$ topology. As $|\nabla f|^2 \in C(M)$,
taking $s \rightarrow 0$ in (\ref{eqnref}) gives
\begin{equation} \label{triangineqlip}
\left( \int_M \langle \nabla f, \nabla \phi_0 \: - \: \nabla \phi_1 
\rangle \: d\mu \right)^2 
\: \le \: \left( \lim_{s \rightarrow 0} 
\frac{W_2(c_0(s), c_1(s))^2}{s^2} \right) \: \int_M |\nabla f|^2 \: d\mu.
\end{equation}

We claim that in fact~\eqref{triangineqlip} holds for any 
$f \in \Lip(M)$. To see this, let $e^{v \triangle}$ be the heat operator
on $M$. Given $f \in \Lip(M)$, for any $v > 0$, $e^{v \triangle} f \in
C^1(M)$. It follows from spectral theory that
$\lim_{v \rightarrow 0} \nabla e^{v \triangle} f \: = \: \nabla f$ in
the Hilbert space of square-integrable vector fields on $M$. Then
$\lim_{v \rightarrow 0} \left| \nabla e^{v \triangle} f \right|^2
\: = \: \left| \nabla f \right|^2$ in $L^1(M, \dvol_M)$. There is
a uniform bound on the $L^\infty$-norm of
$\nabla e^{v \triangle} f \: \equiv \: 
e^{v \triangle_1} df $ for $v \in [0, 1]$. Writing
$\mu \: = \rho \: \dvol_M$ with $\rho \in L^1(M, \dvol_M)$,
for any $N \in \Z^+$ we have
\begin{align}
& \left| \int_M \left| \nabla e^{v \triangle}
 f \right|^2 \: \rho \: \dvol_M \: - \: \int_M
|\nabla f|^2 \: \rho \: \dvol_M \right| \: = \\
& \left| \int_{\rho^{-1}([0,N])} \left( \left| \nabla e^{v \triangle}
 f \right|^2  \: - \: 
|\nabla f|^2  \right) \: \rho \: \dvol_M \: + \:
\int_{\rho^{-1}((N, \infty))} \left( \left| \nabla e^{v \triangle}
 f \right|^2  \: - \: 
|\nabla f|^2  \right) \: \rho \: \dvol_M \right| \: \le \notag \\
& N \: \parallel \left| \nabla e^{v \triangle}
 f \right|^2  \: - \: 
|\nabla f|^2 \parallel_1 \: + \:
\left( \parallel \nabla e^{v \triangle}
 f \parallel_\infty^2 \: + \:
\parallel \nabla 
 f \parallel_\infty^2 \right) \:
\int_{\rho^{-1}((N, \infty))} \rho \: \dvol_M. \notag
\end{align}
For any $\epsilon > 0$, by taking $N$ large we can make
$\left( \parallel \nabla e^{v \triangle}
 f \parallel_\infty^2 \: + \:
\parallel \nabla 
 f \parallel_\infty^2 \right) \:
\int_{\rho^{-1}((N, \infty))} \rho \: \dvol_M$ less than $\epsilon$.
Then by taking $v$ small, we can make
$N \: \parallel \left| \nabla e^{v \triangle}
 f \right|^2  \: - \: 
|\nabla f|^2 \parallel_1$ less than $\epsilon$.
It follows that 
\begin{equation}
\lim_{v \rightarrow 0} \int_M \left| \nabla e^{v \triangle}
 f \right|^2 \: \rho \: \dvol_M \: = \: \int_M
|\nabla f|^2 \: \rho \: \dvol_M.
\end{equation}
By a similar argument,
\begin{equation} \label{similarspectr}
\lim_{v \rightarrow 0} \int_M \left\langle \nabla e^{v \triangle}
f, \nabla \phi_0 \: - \: \nabla \phi_1 \right\rangle \: \rho \: \dvol_M 
\: = \: \int_M
\left\langle 
\nabla f, \nabla \phi_0 \: - \: \nabla \phi_1 \right\rangle \: \rho \: 
\dvol_M.
\end{equation}
Thus (\ref{triangineqlip}) holds for $f$.

In particular, taking $f \: = \: \phi_0 \: - \: \phi_1$ 
in~\eqref{similarspectr} gives
\begin{equation}
\lim_{s \rightarrow 0} 
\frac{W_2(c_0(s), c_1(s))^2}{s^2} \: \ge \: 
\int_M |\nabla \phi_0 \: - \: \nabla \phi_1 |^2
\: d\mu,
\end{equation}
or
\begin{align}
\lim_{s \rightarrow 0} 
\frac{W_2(c_0(s), c_1(s))^2}{s^2} \: \ge \:  
&\frac{W_2(\mu, c_0(s))^2}{s^2} \: + \: \frac{W_2(\mu, c_1(s))^2}{s^2}
\: - \: \\
& 2 \:
\frac{W_2(\mu, c_0(s))}{s} \: \frac{W_2(\mu, c_1(s))}{s} \: 
\frac{\int_M \langle \nabla \phi_0, \nabla \phi_1 \rangle \: 
d\mu}{\sqrt{\int_M | \nabla \phi_0 |^2 \: 
d\mu \: \int_M | \nabla \phi_1 |^2 \: 
d\mu}}. \notag
\end{align}
Equation~\eqref{anglec0c1} follows.
\end{proof}

\subsection{Application to the geometric description of $P_2(M)$}

Let us recall some facts about a finite-dimensional Alexandrov space $Y$ 
with curvature bounded below
\cite{Burago-Burago-Ivanov (2001),Burago-Gromov-Perelman (1992)}.
Let $n$ be the dimension of $Y$. 
A point $y \in Y$ is a {\em regular point} if its tangent cone is 
isometric to $\R^n$. The complement of the regular points is the set $S$ of
{\em singular points}. The regular points $Y - S$ form a dense totally
convex subset of $Y$, but need not be open or closed in $Y$;
see \cite[p. 632-633]{Otsu-Shioya (1994)} for simple but relevant
examples. The existence of a Riemannian metric on $Y$ was studied 
in~\cite{Otsu (1997),Otsu-Shioya (1994),Perelman}. We recall the
results of~\cite{Perelman}: There is a dense open totally
convex subset $Y^\delta$ of $Y$, containing $Y-S$, 
which is a topological manifold with DC (=difference of concave)
transition maps; there is a Riemannian metric $g$ on $Y^\delta$ which in 
local charts is bounded, measurable and of bounded variation, with the 
restriction of $g$ to $Y-S$ being continuous; the Christoffel symbols exist
as measures in local charts on $Y^\delta$; the lengths of curves in $Y-S$ can
be computed using $g$. 

There is an evident analogy between $Y - S \subset Y$ and
$P^{ac}_2(M) \subset P_2(M)$. The arguments of the above papers do not 
directly extend to infinite-di\-men\-sio\-nal Alexandrov spaces. Nevertheless,
in order to make a zeroth order approximation to a Riemannian geometry
on $P_2(M)$, it makes sense to look at the tangent cones.
We recall that for a finite-dimensional Riemannian manifold, the
tangent cone at a point $p$ is isometric to $T_pM$ equipped with the
inner product coming from the Riemannian metric at $p$.

\begin{proposition} \label{formalmetric}
Let $M$ be a smooth compact connected Riemannian 
manifold. If $M$ has nonnegative sectional curvature then for each
absolutely continuous measure $\mu \in P_2^{ac}(M)$, the tangent 
cone of $P_2(M)$ at $\mu$ is an inner product space.
\end{proposition}

\begin{proof}
Given $\mu \in P^{ac}_2(M)$,
we consider the space $\Sigma_{\mu}^\prime$ of equivalence classes
of geodesic segments emanating from $\mu$, 
with the equivalence relation identifying
two segments if they form a zero angle at $\mu$
\cite[Section 9.1.8]{Burago-Burago-Ivanov (2001)} 
(which in the case of curvature bounded below means that one
segment is contained in the other).
The metric on $\Sigma_{\mu}^\prime$ is the angle.
By definition, the space of directions 
$\Sigma_{\mu}$ is the metric completion of
$\Sigma^\prime_{\mu}$. The tangent cone $K_\mu$ is the
union of $\Sigma_{\mu} \times \R^+$ and a ``vertex'' point,
with the metric described in 
\cite[Section 10.9]{Burago-Burago-Ivanov (2001)}.

We first note that any Wasserstein geodesic $\{\mu_t\}_{t \in [0,1]}$
emanating from $\mu$ is of the form
$\mu_t \: = \: (F_t)_* \mu$, with $F_t$ as in Subsection \ref{recall}
\cite[Theorem 2.47]{Villani (2003)}. It follows that we can apply the angle
calculation in Proposition \ref{anglecomp} to all
Wasserstein geodesics emanating from $\mu$. (One can also use the fact that
any such Wasserstein geodesic $\{\mu_t\}_{t \in [0,1]}$ has
$\mu_t \in P_2^{ac}(M)$ for $t \in
[0, 1)$ \cite[Lemma 22]{Bernard-Buffoni (2004)}; see
\cite[Proposition 5.9(iii)]{Villani (2003)} for the $\R^n$-case.)

Now to identify $\Sigma_{\mu}^\prime$, 
consider the space ${\cal S}$ of Lipschitz functions on $M$
that are $rd^2$-concave for some $r > 0$. In terms of the function $\phi$,
we can identify the geodesic
segments from $\mu$ with ${\cal S}^\prime = {\cal S}/\R$, where $\R$
acts additively on ${\cal S}$.
There is an action of $\R^+$ on ${\cal S}^\prime$ by multiplication.  
As the angle between geodesic segments is given by~\eqref{anglec0c1}, we can
identify $\Sigma^\prime_{\mu}$ with the corresponding quotient of 
the space of Lipschitz functions on $M$ that are $rd^2$-concave for some 
$r > 0$.

We can approximate a Lipschitz function on $M$ with respect to
the quadratic form $Q(\phi) \: = \: \int_M |\nabla \phi|^2 \: d\mu$ by
functions that are $rd^2$-concave for various $r > 0$, for example
by flowing the Lipschitz function for a short time under the heat equation 
on $M$.  Hence when considering
the metric completion of $\Sigma^\prime_{\mu}$, it doesn't matter whether 
we start with Lipschitz 
functions on $M$ that are $rd^2$-concave for some $r > 0$ or arbitrary
Lipschitz functions on $M$. It follows that
$K_{\mu}$ is the inner product space constructed by starting with
$\Lip(M)$, quotienting by the kernel of $Q$ and
taking the metric completion with respect to $Q$.
\end{proof}

The tangent cone constructed in 
Proposition \ref{formalmetric} agrees with
the formal infinite-dimensional Riemannian metric on $P_2(M)$
considered by Otto \cite{Otto (2001)}. Proposition \ref{formalmetric}
can be considered as a way of making this formal Riemannian metric
rigorous, and Theorem \ref{thmalex} as a rigorous version of Otto's formal 
argument that his Riemannian metric on $P_2(\R^n)$ has
nonnegative sectional curvature. 

\section{Some properties of the functionals $U_\nu$} \label{applsc}

The goal of this section is to gather several results about functionals 
of the form $U_\nu$, as defined in
Definition \ref{deffunctional}. (We will
generalize slightly to allow $X$ to be
a compact Hausdorff space, but the same definition makes sense.)

We show that $U_\nu(\mu)$ is lower semicontinuous in $\mu$ and
$\nu$. Such a lower semicontinuity in $\mu$ is well known in the setting of 
the weak topology on $L^p$ functions,
but we need to prove it in the weak-$*$ topology on Borel measures.
To do so, we derive a Legendre-type formula for $U_\nu(\mu)$;
this Legendre formula is also well-known in certain cases,
e.g. $U(r)=r\log r$, but it is not so easy to find a precise
reference for general nonlinearities.

We will also show that $U_\nu(\mu)$ is nonincreasing under pushforward.

For notation, if $U$ is a convex function then
$U^*$ is its Legendre transform and $U'$ is its right-derivative.

\subsection{The functional $U_\nu$ via Legendre transform}

We start by recalling, without proof, a consequence of Lusin's theorem.

\begin{theorem} \label{lusin}
Let $X$ be a compact Hausdorff space. Let $\mu$ be a Borel
probability measure
on $X$. Then for all $f \in L^\infty(X)$ there is
a sequence $\{f_k\}_{k=1}^\infty$ of continuous functions on $X$ such that\\
(i) $\inf f \leq \inf f_k \leq \sup f_k \leq \sup f$ and \\
(ii) $\lim_{k \rightarrow \infty} f_k(x) \: = \:  f(x)$ for 
$\mu$-almost all $x \in X$.
\end{theorem}

We now prove a useful Legendre-type representation formula.

\begin{theorem} \label{thmdual}
Let $X$ be a compact Hausdorff space. Let $U:[0, \infty)\to\R$ be a continuous
convex function with $U(0) = 0$. 
Given $\mu, \nu \in P(X)$, we have
\begin{align} \label{fmldualL8}
U_\nu(\mu) & = 
\sup \left \{ \int_X\varphi\,d\mu - \int_X U^*(\varphi)\,d\nu \: : \: 
\quad \varphi\in L^\infty(X), \> \varphi\leq U'(\infty) \right \} \\ 
& =  \sup _{M \in \Z^+} \sup
\left \{ \int_X\varphi\,d\mu - 
\int_X U^*(\varphi)\,d\nu \: : \:  \quad \varphi\in C(X), \> 
\varphi \: \leq \: U'(M) \right \}. \notag
\end{align}
\end{theorem}

\begin{remark}
The reason to add the condition $\varphi \: \leq \: U'(M)$ is to ensure that
$U^*(\phi)$ is continuous on $X$. This will be used in the proof of
Theorem \ref{thmlsc}(i).
\end{remark}
\begin{proof}[Proof of Theorem~\ref{thmdual}]
As an initial reduction, if $U^\prime(0) \:  = \: U^\prime(\infty)$ then $U$ is linear
and the result of the theorem is easy to check. If
$U^\prime(0) \: < \: U^\prime(\infty)$, choose $c \in (U^\prime(0), U^\prime(\infty))$.
Replacing $U(r)$ by $U(r) \: - \: cr$, we can reduce to the case when
$U^\prime(0) < 0$ and $U^\prime(\infty) > 0$.

Let $\mu=\rho\nu+\mu_s$ be the Lebesgue decomposition of $\mu$ 
with respect to $\nu$.
Let $S$ be a measurable set such that $\mu_s(S)=\mu_s(X)$ and $\nu(S)=0$. 
Without loss
of generality, we may
assume that $\rho<\infty$ everywhere on $X-S$, and we set
$\rho=\infty$ on $S$.

We will prove that
\begeq\label{ineq1dual}
U_\nu(\mu) \geq
\sup \left \{ \int_X\varphi\,d\mu - \int_X U^*(\varphi)\,d\nu \: : \: 
\quad \varphi\in L^\infty(X), \> \varphi\leq U'(\infty) \right \}
\endeq
and
 \begeq\label{ineq2dual} U_\nu(\mu) \leq
\sup_{M \in \Z^+}\: \sup \left \{ \int_X\varphi\,d\mu - 
\int_X U^*(\varphi)\,d\nu \: : \:  \quad \varphi\in C(X), \> 
U' \left ( \frac1{M} \right ) \leq
\varphi\leq U'(M) \right \}.
\endeq 
As the right-hand-side of (\ref{ineq2dual}) is clearly less than or equal to the right-hand-side of
(\ref{ineq1dual}), this will imply that
\begin{align}
U_\nu(\mu) & = 
\sup \left \{ \int_X\varphi\,d\mu - \int_X U^*(\varphi)\,d\nu \: : \: 
\quad \varphi\in L^\infty(X), \> \varphi\leq U'(\infty) \right \} \\ 
& =  \sup _{M \in \Z^+} \sup
\left \{ \int_X\varphi\,d\mu - 
\int_X U^*(\varphi)\,d\nu \: : \:  \quad \varphi\in C(X), \> 
U' \left ( \frac1{M} \right ) \leq
\varphi \: \leq \: U'(M) \right \}, \notag
\end{align}
which in turn implies the theorem.

The proof of~\eqref{ineq1dual} is obtained by a direct argument : for any
$\varphi\in L^\infty(X)$ with $\varphi \: \le \: U^\prime(\infty)$, 
we will show that
\begeq \label{eqntoshow}
U_\nu(\mu) \geq
 \int_X\varphi\,d\mu - \int_X U^*(\varphi)\,d\nu.
\endeq
We may assume that $U(\rho)\in L^1(X, \nu)$, as 
otherwise there is nothing to prove.
For all $x\in X$, we have
\begeq \label{Uvarphirho} 
U(\rho(x)) \: - \: \varphi(x)\rho(x) \: \ge \:  - \: U^*(\varphi(x)). \endeq
Also, $\varphi\rho \in L^1(X, \nu)$. 
Integrating (\ref{Uvarphirho}) with respect to $\nu$ gives
\begin{equation} \label{intUvarphirho}
\int_X U(\rho) \: d\nu \: 
- \: \int_X \varphi \rho \: d\nu \: \ge \:  - \int_X  U^*(\varphi) \: 
d\nu.
\end{equation}
On the other hand, since $\varphi\leq U'(\infty)$, we also have
\begeq\label{otherUvarphirho} U'(\infty)\mu_s(X) \geq 
\int_X \varphi\,d\mu_s.\endeq
Adding (\ref{intUvarphirho}) and~\eqref{otherUvarphirho} 
gives (\ref{eqntoshow}).
\med

To prove~\eqref{ineq2dual}, it suffices to show the existence
of a sequence $\{\varphi_M\}_{M=1}^\infty$ in $C(X)$ such that
\begeq \begin{cases} \label{toineq2dual}
U' \left (\frac1{M} \right ) \leq \varphi_M \leq U'(M) \text{ and } \\
U_\nu(\mu) \leq \liminf_{M\to\infty}
\left( \int_X \varphi_M\,d\mu - \int_X U^*(\varphi_M)\,d\nu \right)
.\end{cases} \endeq
\med

For $M\geq 1$, we define
\begin{equation} 
\rho_M = \max \left ( \frac1{M}, \min (\rho,M) \right ).
\end{equation}
It is clear that\\
\quad (i) $M^{-1}\leq \rho_M\leq M$;\\
\quad (ii) for all $x\in X$, $\lim_{M \rightarrow \infty} \rho_M(x) = \rho(x)$;\\
\quad (iii) if $0<\rho(x)<\infty$
then $\rho_M(x)=\rho(x)$ for $M$ large enough.

Choose $\epsilon > 0$ so that $U$ is nonincreasing
on $[0, \epsilon)$ and nondecreasing on $[\epsilon, \infty)$. 
Monotone convergence implies that 
\begin{equation}
\int_{\rho^{-1}[\epsilon, \infty)} U(\rho)\,d\nu = \lim_{M\to\infty} 
\int_{\rho^{-1}[\epsilon, \infty)} U(\rho_M)\,d\nu
\end{equation}
and
\begin{equation}
\int_{\rho^{-1}[0, \epsilon)} U(\rho)\,d\nu = \lim_{M\to\infty} 
\int_{\rho^{-1}[0, \epsilon)} U(\rho_M)\,d\nu.
\end{equation}
Hence
\begeq \label{monconvrhoM} 
\int_X U(\rho)\,d\nu = \lim_{M\to\infty} \int_X U(\rho_M)\,d\nu.
\endeq
\med

Define now a function $\ov{\varphi}_M \: : \: X \rightarrow \R$ by 
\begin{equation} 
\ov{\varphi}_M = U' (\rho_M). 
\end{equation}
Since $U'$ is nondecreasing, we have
\begin{equation} 
U' \left ( \frac1{M} \right )\leq \ov{\varphi}_M \leq U'(M). 
\end{equation}

We also have the pointwise equality
\begin{equation} 
U(\rho_M) = \ov{\varphi}_M \rho_M - U^*(\ov{\varphi}_M).
\end{equation}
All of the functions appearing in this identity are integrable with respect to $\nu$, so
\begeq\label{intXUrhoM} 
\int_X U(\rho_M)\,d\nu = \int_X \ov{\varphi}_M\rho_M\,d\nu -
\int_X U^*(\ov{\varphi}_M)\,d\nu. 
\endeq

Our first goal is to prove that
\begeq\label{goalUU} 
U_\nu(\mu) \leq \liminf_{M\to\infty}
\left( \int_X \ov{\varphi}_M\,d\mu - \int_X U^*(\ov{\varphi}_M)\,d\nu \right).
\endeq
If this is true then we have shown that
the sequence $\{\ov{\varphi}_M\}_{M=1}^\infty$ satisfies all of the
properties required in~\eqref{toineq2dual}, except maybe continuity.
We split~\eqref{goalUU} into two parts :
\begeq\label{UUpart1}
U'(\infty) \mu_s(X) = \lim_{M\to\infty} \int_X \ov{\varphi}_M\,d\mu_s \endeq
and
\begeq\label{UUpart2}
\int_X U(\rho)\,d\nu \le \liminf_{M\to\infty} \left [\int_X \ov{\varphi}_M\,
\rho \, d\nu
- \int_X U^*(\ov{\varphi}_M)\,d\nu \right ]. \endeq

To prove~\eqref{UUpart1}, we write
\begin{equation} 
\int_X \ov{\varphi}_M\,d\mu_s = \int_S \ov{\varphi}_M \,d\mu_s
= U'(M) \mu_s(S) = U'(M)\mu_s(X) 
\stackrel{M \rightarrow \infty}{\longrightarrow} U'(\infty)\mu_s(X). 
\end{equation}

To prove~\eqref{UUpart2}, we note that for large $M$,
 \begin{equation} 
U'(\rho_M)(\rho_M-\rho) \leq 0.
\end{equation}
Indeed, if $1/M\leq\rho\leq M$ then $\rho_M=\rho$; 
if $\rho>M$ then $\rho>\rho_M$ and $U'(\rho_M)\geq 0$;
while if $\rho<1/M$ then $\rho<\rho_M$ and $U'(\rho_M)\leq 0$. 
Thus
\begin{equation} 
\int_X \ov{\varphi}_M\,\rho_M\,d\nu \leq \int_X \ov{\varphi}_M\,\rho\,d\nu.
\end{equation}
Combining this with~\eqref{monconvrhoM} and~\eqref{intXUrhoM}, we find
\begin{align}
\int_X U(\rho)\,d\nu & = \lim_{M\to\infty} \int_X U(\rho_M)\,d\nu \\
& = \lim_{M\to\infty} \left [ \int_X \ov{\varphi}_M\,\rho_M\,d\nu -
\int_X U^*(\ov{\varphi}_M)\,d\nu \right ]  \notag\\
& \leq \liminf_{M\to\infty} \left [ \int_X \ov{\varphi}_M\rho\,d\nu - 
\int_X U^*(\ov{\varphi}_M)\,d\nu \right ]. \notag
\end{align}
\med

This proves (\ref{goalUU}).
To conclude the proof of the theorem, 
it suffices to show that for any $M \in \Z^+$ there is
a {\em continuous} function $\varphi_M$ such that $U'(1/M)\leq \varphi_M
\leq U'(M)$
and
\begin{equation} 
\left | \left ( \int_X 
\varphi_M\,d\mu - \int_X U^*(\varphi_M)\,d\nu \right ) - 
\left ( \int_X \ov{\varphi}_M\,d\mu - 
\int_X U^* (\ov{\varphi}_M)\,d\nu \right ) \right | \leq
\frac1{M}. 
\end{equation}

Fix $M$. By Theorem~\ref{lusin}, 
there is a sequence $\{\psi_k\}_{k=1}^\infty$ of continuous
functions such that $U'(1/M)\leq \inf \ov{\varphi}_M \leq \inf\psi_k \leq
\sup \psi_k\leq \sup\ov{\varphi}_M \leq U'(M)$ (in particular 
$\{\psi_k\}_{k=1}^\infty$ is
uniformly bounded) and $\lim_{k \rightarrow \infty}
\psi_k(x) \: = \: 
\ov{\varphi}_M(x)$ $(\mu+\nu)$-almost everywhere.

At this point we note that
\begin{equation} 
U^*(p) = \sup_{r\geq 0} [pr- U(r)] 
\end{equation}
is bounded below by $-U(0) \: = \: 0$ and is a nondecreasing function of $p$.
Also,
\begin{equation} 
0 \leq U^*(U'(1/M)) \leq U^*(U'(M)) = MU'(M)-U(M) <\infty.
\end{equation}
Thus $U^*$ is bounded on $[U'(1/M), U'(M)]$. Since it is also lower 
semicontinuous and convex, 
it follows that it is actually continuous on that interval. So 
$\{ U^*(\psi_k) \}_{k=1}^\infty$ converges
$\nu$-almost everywhere to $U^*(\ov{\varphi}_M)$. By dominated convergence,
\begin{equation}
\lim_{k \rightarrow \infty} \int_X U^*(\psi_k)\,d\nu \: = \:
\int_X U^*(\ov{\varphi}_M)\,d\nu.
\end{equation}
Also by dominated convergence,
\begin{equation} 
\lim_{k \rightarrow \infty}
\int_X \psi_k\,d\mu \: = \: \int_X \ov{\varphi}_M\,d\mu.
\end{equation}
We can conclude by choosing $\varphi_M = \psi_k$ for some large $k$.
\end{proof}

\subsection{Lower semicontinuity and contraction}

The following theorem is an easy consequence of the duality formulas
established above.

\begin{theorem} \label{thmlsc}
Let $X$ and $U$ satisfy the assumptions of
Theorem~\ref{thmdual}. Then

\quad (i) $U_\nu(\mu)$ is a lower semicontinuous function of $(\mu,\nu)
\in P(X) \times P(X)$.
That is, if $\{\mu_k\}_{k=1}^\infty$ and 
$\{\nu_k\}_{k=1}^\infty$ are sequences in $P(X)$
with $\lim_{k \rightarrow \infty} \mu_k \: = \: \mu$ and
$\lim_{k \rightarrow \infty} \nu_k \: = \: \nu$ in the weak-$*$ topology then
\begin{equation} 
U_\nu(\mu) \leq \liminf_{k\to\infty} U_{\nu_k}(\mu_k).
\end{equation}

\quad (ii) $U_\nu(\mu)$ is nonincreasing under pushforward.
That is, if $Y$ is a compact Hausdorff space and
$f \: : \: X \rightarrow Y$ is a Borel map then
\begin{equation} 
U_{f_* \nu} (f_*\mu) \leq U_\nu(\mu).
\end{equation}
\end{theorem}

\begin{proof}
Using the second representation formula in Theorem~\ref{thmdual},
we can write 
\begin{equation} 
U_\nu(\mu) = \sup_{(L_1,L_2)\in {\cal L}} [L_1(\mu) +L_2(\nu)],
\end{equation}
where ${\cal L}$ is a certain subset of $C(X) \oplus C(X)$, and
$L_1$ and $L_2$ define continuous linear functionals on the space of
measures $C(X)^*$. As the supremum of a set of lower semicontinuous functions
(in particular linear functions) is lower semicontinuous,
it follows that $U_\nu$ is lower semicontinuous in $(\mu,\nu)$.
\sm

To prove part (ii), we use the first representation formula in 
Theorem~\ref{thmdual}
to obtain
\begin{align} 
U_{f_*\nu}(f_*\mu) & = \sup \left \{ \int_Y \varphi\,d(f_*\mu) - 
\int_Y U^*(\varphi)\,d(f_*\nu);
\quad \varphi\in L^\infty(Y), \> \varphi\leq U'(\infty) \right \} \\
& = \sup \left \{ \int_X (\varphi\circ f)\,d\mu - 
\int_X U^*(\varphi\circ f)\,d\nu;
\quad \varphi\in L^\infty(Y), \> \varphi\leq U'(\infty) \right \}. \notag
\end{align}
If $\varphi \in L^\infty(Y)$ and $\varphi\leq U'(\infty)$ then
$\varphi\circ f \in L^\infty(X)$ and
$\varphi\circ f \leq U'(\infty)$.
So the above supremum is bounded above by
\begin{equation} 
U_\nu(\mu) \: = \: \sup \left \{ \int_X \psi\,d\mu - \int_X U^*(\psi)\,d\nu;
\quad \psi\in L^\infty(X), \> \psi\leq U'(\infty) \right \}.
\end{equation}
\end{proof}

\section{Approximation in $P_2(X)$} \label{secmoll}

In this section we show how to effectively approximate a measure
$\mu \in P_2(X,\nu)$ by measures $\{\mu_k\}_{k=1}^\infty$ whose densities, with
respect to $\nu$, are continuous. The approximation will be such that
$\lim_{k \rightarrow \infty} \mu_k \: = \: \mu$ in the weak-$*$ topology
and $\lim_{k \rightarrow \infty} U_\nu(\mu_k) \: = \: U_\nu(\mu)$.

We first construct a mollification operator on measures, in terms of a
partition of unity for $X$.  We then use finer and finer partitions of unity
to construct the sequence $\{\mu_k\}_{k=1}^\infty$.

\subsection{Mollifiers}

Let $(X,d,\nu)$ be a compact 
metric space equipped with a reference Borel probability
measure $\nu$.
Let $K:X\times X\to[0, \infty)$ be a  symmetric continuous kernel
satisfying
\begin{equation} 
\forall x \in \supp(\nu), \quad \int_X K(x,y)\,d\nu(y)=1.
\end{equation}
For $\rho\in L^1(X,\nu)$, define $K\rho \in C(X)$ by
\begin{equation} 
(K\rho) (x) = \int_X K(x,y)\,\rho(y)\,d\nu(y).
\end{equation}
Note that
$\int_X K\rho \: d\nu \: = \: \int_X \rho \: d\nu$ and the map
$\rho \rightarrow 1_{\supp(\nu)} \: K\rho$ 
is a bounded operator on $L^1(X,\nu)$ with norm~1.
For $\mu \in P_2(X, \nu)$, define $K\mu \in P_2(X, \nu)$ by
saying that for $f \in C(\supp(\nu))$,
\begin{equation} \label{Kdef}
\int_{\supp(\nu)}  f \: d(K\mu) \: = \: \int_X (Kf) \: d\mu.
\end{equation}
More explicitly,
\begin{equation}
K\mu \: = \left( \: \int_X K(\cdot,y) \: d\mu(y) \right) \: \nu.
\end{equation}
In particular, $K\mu \in P^{ac}_2(X, \nu)$ is the product of a continuous
function on $X$ with $\nu$. The notation is consistent, in the sense that if
$\rho \in L^1(X, \nu)$ then $K(\rho \: d\nu) \: = \: K(\rho) \: d\nu$.
Moreover, taking $f=1$ in (\ref{Kdef}), 
it follows that if $\mu$ is a probability
measure then $K\mu$ is a probability
measure.
\med

\begin{theorem} \label{approxthm}
Let $(X,d)$ be a compact metric space equipped with a reference Borel
probability
measure $\nu$. Then there is a
sequence $\{K_I\}_{I=1}^\infty$ of continuous 
nonnegative kernels with
the following properties : \\
(i) $\forall x,y \in X,\quad K_I(x,y) = K_I(y,x)$. \\
(ii) $\forall x \in \supp(\nu), \quad \int_X K_I(x,y)\,d\nu(y)=1$. \\
(iii) There is a sequence $\{\epsilon_I\}_{I=1}^\infty$ converging
to zero so that $K_I(x,y) =0$ whenever $d(x,y)\geq \epsilon_I$. \\
(iv) For all $\mu\in P_2(X, \nu)$,
$\lim_{I \rightarrow \infty} K_I\mu = \mu$
in the weak-$*$ topology.
\end{theorem}
\begin{proof}
Let ${\cal C} \: \equiv \: \{U_j\}$ be a finite open cover of $X$.
Let $\{ \phi_j \}$ be a subordinate partition of unity. Put 
\begin{equation}
K_{\cal C}(x, y) \: = \: \sum_{j \: : \: \int_X \phi_j \: d\nu \: > \: 0} \frac{\phi_j(x) \: 
\phi_j(y)}{\int_X \phi_j \: d\nu}.
\end{equation}
Then $K_{\cal C}(x,y)=K_{\cal C}(y,x)$. If $x \in \supp(\nu)$ and
$\phi_j(x) > 0$ then $\int_X \phi_j \: d\nu \: >\: 0$, so
$\int_X K_{\cal C}(x, y) \: d\nu(y) \: = \: 1$.
Properties (i) and (ii) are ensured if our sequence
is made of such kernels. Moreover, (iii) will be satisfied if
each $\phi_j$ has support in a small ball of radius $\epsilon_I/2$. 

Given $\delta > 0$, let ${\cal C}_\delta$ denote a finite
cover of $X$ by $\delta$-balls.
Given $f \in C(\supp(\nu))$ and 
$\epsilon > 0$, suppose that $\delta >0$ is such that 
$|f(x^\prime) \: - \: f(x)| \: \le \:
\epsilon$ whenever $x^\prime, x \in \supp(\nu)$ satisfy
$d(x^\prime, x) \: \le \: 2 \delta$. With
such a cover ${\cal C}_\delta$, if $x \in \supp(\nu)$ has
$\phi_j(x) > 0$ then
\begin{equation}
\left| \frac{\int_X f \phi_j \: d\nu}{\int_X \phi_j \: d\nu} \: - \: 
f(x) \right| 
\: \le \: \epsilon.
\end{equation}
As
\begin{equation}
(K_{{\cal C}_\delta} f)(x) \: = \: \sum_{j \: : \: \int_X \phi_j \: d\nu \: > \: 0} \phi_j(x) \: 
\frac{\int_X f \phi_j \: d\nu}{\int_X \phi_j \: d\nu}
 \: = \: \sum_{j \: : \: \phi_j(x) \: > \: 0} \phi_j(x) \: 
\frac{\int_X f \phi_j \: d\nu}{\int_X \phi_j \: d\nu}
\end{equation}
and
\begin{equation}
f(x) \: = \: \sum_j \phi_j(x) \: f(x) \: = \: \sum_{j \: : \: \phi_j(x) \: > \: 0} \phi_j(x) \: f(x),
\end{equation}
it follows that
\begin{equation}
\left| (K_{{\cal C}_\delta} f)(x) \: - \: f(x) \right| \: \le \: 
\epsilon.
\end{equation}
Thus $\lim_{I \rightarrow \infty} \: K_{{\cal C}_{1/I}} f \: = \: f$ in
$C(\supp(\nu))$.

Now put $K_I \: = \: K_{{\cal C}_{1/I}}$. For any 
$f\in C(\supp(\nu))$, we have
\begin{equation}
\lim_{I \rightarrow \infty} \int_X f \: d K_{I} \mu \: = \: 
\lim_{I \rightarrow \infty} \int_X (K_I f) \: d\mu \: = \:
\int_X f \: d\mu.
\end{equation}
This proves (iv).
\end{proof}

\subsection{Approximation by continuous densities}

\begin{theorem} \label{thmapproxU}
Let $U$ be a continuous convex function on $[0, \infty)$
with $U(0) \: = \: 0$.
Let $(X,d)$ be a compact metric space equipped with a Borel probability
measure $\nu$.
Let $\mu \in P_2(X,\nu)$ satisfy
$U_\nu(\mu)<\infty$. Then
there is a sequence $\{f_k\}_{k=1}^\infty$ in $C(X)$
such that $\lim_{k \rightarrow \infty} f_k\nu \: = \: \mu$ in the
weak-$*$ topology and 
$\lim_{k \rightarrow \infty} U_\nu(f_k\nu) \: = \: U_\nu(\mu)$.
\end{theorem}

\begin{proof}
We introduce an sequence of mollifying kernels $K_I$ as above.
We will prove that $\{K_I \mu\}_{I=1}^\infty$
does the job. Each $K_I \mu$ is the product of a continuous function
on $X$ with 
$\nu$.  Theorem \ref{approxthm}(iv) gives that 
$\lim_{I \rightarrow \infty} K_I \mu \: = \: \mu$.
By Theorem \ref{thmlsc}(i), 
$U_\nu(\mu) \: \le \: \liminf_{I \rightarrow
\infty} U_\nu(K_I \mu)$.
Hence it suffices to show that
$U_\nu(K_I \mu) \: \le \: U_\nu(\mu)$ for all $I$.

Before giving the general proof, it might be enlightening 
to first consider two ``extreme'' cases. \\ \\
{\bf First particular case: $\mu$ is absolutely continuous.} 
Assume that $\mu=\rho\nu$.
We write $K_I \mu(x)$ for the density of $K_I \mu$ with respect to
$\nu$.
By Jensen's inequality,
for $x \in \supp(\nu)$,
\begin{equation} 
U(K_I\mu(x)) = U \left ( \int_X K_I(x,y)\,\rho(y)\,d\nu(y) \right )
\leq \int_X K_I(x,y) \: U(\rho(y))\,d\nu(y).
\end{equation}
It follows that
\begin{equation} 
U_\nu(K_I\mu) = \int_X U(K_I\mu(x)) \, d\nu(x) \leq \int_X U(\rho)\,d\nu.
\end{equation}
\noindent
{\bf Second particular case: $\mu$ is completely singular.}
With the notation used before, $\mu=\mu_s$.
We write $K_I \mu_s(x)$ for the density of 
$K_I \mu_s$ with respect to
$\nu$.
Then
\begin{equation} 
\int_X U((K_I\mu_s)(x))\,d\nu(x) \leq 
U'(\infty) \int_X (K_I\mu_s)(x)\,d\nu(x) = U'(\infty) \mu_s(X).
\end{equation}
\noindent
{\bf General case:} Now we introduce the Lebesgue decomposition
$\mu=\rho\nu +\mu_s$. In view of the first particular case above,
we may assume that $\mu_s(X)>0$. 
If $U'(\infty)=\infty$ then $U_\nu(\mu)=\infty$, so we can restrict to the
case $U'(\infty)<\infty$. 

We write $K_I \mu(x)$ for the density of $K_I \mu$ with respect to
$\nu$, and similarly for $K_I \mu_s(x)$.
For all $\theta\in (0,1)$, there
is a pointwise inequality
\begin{equation} \label{comb0}
U(K_I\mu) \: = \: U( K_I\rho + K_I \mu_s)
\: \leq \: \theta \: U \left ( \frac{K_I \rho}{\theta} \right)
\: + \: (1-\theta) \: U \left ( \frac{K_I\mu_s}{1-\theta} \right ). 
\end{equation}
Then
\begin{align} \label{comb1}
U_\nu(K_I \mu) \: = \: \int_X U(K_I\mu(x)) \: d\nu(x) \: & \leq \: \theta \:
\int_X U \left ( \frac{K_I \rho}{\theta} \right ) \: d\nu
\: + \: U^\prime(\infty) \: \int_X K_I\mu_s(x) \: d\nu(x) \\
& = \: \theta \:
\int_X U \left ( \frac{K_I \rho}{\theta} \right ) \: d\nu
\: + \: U^\prime(\infty) \: \mu_s(X). \notag
\end{align}

As $K_I \rho \in C(X)$, 
\begin{equation} \label{comb2}
\lim_{\theta\to 1^-} \: \left( \theta \:
\int_X U \left ( \frac{K_I \rho}{\theta} \right ) \: d\nu \right)
\: = \:
\int_X U(K_I \rho) \: d\nu.
\end{equation}
As in the proof of the first particular case,
\begin{equation} \label{comb3}
\int_X U(K_I \rho) \: d\nu \: \le \: \int_X U(\rho) \: d\nu.
\end{equation}
Combining (\ref{comb1}), (\ref{comb2}) and (\ref{comb3}) gives
$U_\nu(K_I \mu) \: \le \: U_\nu(\mu)$.
\end{proof}

\section{Hessian calculations} \label{appformaldc}

How can one check, in practice, the displacement convexity
of a given functional on $P_2(X)$, say in the case when $X$ is 
a smooth compact connected Riemannian manifold $M$?
We provide below some explicit computations to that effect,
following~\cite{Otto-Villani (2000)}. The computations are
purely formal and we do not rigorously justify them, ignoring
all regularity issues. 
Although formal, these computations motivate the definition of
nonnegative $N$-Ricci curvature in terms of displacement
convexity.

Denote by $\dvol_M$ the Riemannian density on $M$, and
introduce a reference measure
\begin{equation}
d\nu \: = \: e^{-\Psi} \: \dvol_M,
\end{equation}
with $\Psi \in C^\infty(M)$ satisfying $\int_M e^{-\Psi} \: \dvol_M \: = \: 1$.

The direction vector along a curve $\{\mu_t\}$ in the space of
probability measures $P(M)$
can be ``represented'' as
\begin{equation} \label{appeqn1}
\frac{\partial \mu}{\partial t} \: = \: - \: \nabla \cdot (\mu \nabla \Phi),
\end{equation}
where $\Phi \equiv \Phi(t)$ is a function on $M$ that is defined
up to constants.
The meaning of (\ref{appeqn1}) is that
\begin{equation}
\frac{d}{dt} \: \int_M \xi \: d\mu \: = \: 
\int_M \nabla \xi \cdot \nabla \Phi \: d\mu 
\end{equation}
for all $\xi \in C^\infty(M)$.
Thus we can parametrize the tangent space $T_\mu P(M)$ by the function $\Phi$.
Otto's formal inner product on $T_\mu P(M)$ is given by the
quadratic form $\int_M \langle \nabla \Phi, \nabla \Phi \rangle \:
d\mu$.

With this Riemannian metric, 
the geodesic equation in $P(M)$ becomes
\begin{equation}
\partial_t \Phi \: + \: \frac12 \: |\nabla \Phi|^2 \: = \: 0.
\end{equation}
It has the solution 
\begin{equation}
\Phi(t)(x) \: = \: \inf_{y \in X} \left[ \Phi(0)(y) \: + \: \frac{d(x,y)^2}{2t} \right].
\end{equation}
The corresponding length
metric on $P(M)$ is formally the Wasserstein metric.

Let $U$ be a continuous convex function on $[0, \infty)$
that is $C^2$-regular on $(0, \infty)$. 
Put
\begin{equation} 
E(\mu) = \int_M U \left( \frac{d\mu}{d\nu} \right) \: d\nu.
\end{equation}
Recall that
\begin{equation}
\p(r) \: = \: r U^{\prime}(r) \: - \: U(r), \qquad
\pp(r) \: = r\p^{\prime}(r) - \p(r).
\end{equation}
Then along a curve $\{\mu_t\}$ in $P(M)$, the derivative of
$E(\mu_t)$ is given by  
\begin{align} \label{EE}
\frac{dE}{dt} \: & = \: \int_M U^\prime \left( \frac{d\mu}{d\nu} \right) \: 
\partial_t \frac{d\mu}{d\nu} \: d\nu \\
& = \: \int_M \nabla U^\prime \left( \frac{d\mu}{d\nu} \right) \cdot 
\nabla \Phi \: d\mu \notag \\
& = \: \int_M \nabla \Phi \cdot \left( \frac{d\mu}{d\nu} \: 
U^{\prime \prime} \left( \frac{d\mu}{d\nu} \right) \: \nabla \frac{d\mu}{d\nu}
\right) \: d\nu \notag \\
& = \: \int_M \nabla \Phi \cdot \nabla \p \left( \frac{d\mu}{d\nu} \right) \: 
d\nu \notag \\
& = \: \int_M (- \: \triangle \Phi \: + \: \nabla \Psi \cdot \nabla \Phi) \:
\p \left( \frac{d\mu}{d\nu} \right) \: 
d\nu. \notag
\end{align}
Parametrizing $T_\mu P(M)$ by $\{\Phi\}$, equation
(\ref{EE}) shows 
in particular that 
$\grad E$ is represented by the function
$U^\prime \left( \frac{d\mu}{d\nu} \right)$.

To compute the Hessian of $E$, we assume that $\{\mu_t\}$ is a geodesic
curve in $P_2(M)$. Then
\begin{align} \label{mess2}
\frac{d^2 E}{dt^2} \: = \: 
\int_M & \left( \triangle \left( \frac12 |\nabla \Phi|^2
\right) \: - \: 
\nabla \Psi \cdot \nabla \left( \frac12 |\nabla \Phi|^2 \right) \right)
\: \p \left( \frac{d\mu}{d\nu} \right) \: d\nu \\
& + \: \nabla \left( \left( - \: \triangle \Phi \: + \: \nabla \Psi \cdot
\nabla \Phi \right) \: \p^\prime \left( \frac{d\mu}{d\nu} \right) \right)
\cdot \nabla \Phi \: d\mu. \notag
\end{align}
Now
\begin{align} \label{mess3}
& \int_M \nabla \left( \left( - \: \triangle \Phi \: + \: \nabla \Psi \cdot
\nabla \Phi \right) \: \p^\prime \left( \frac{d\mu}{d\nu} \right) \right)
\cdot \nabla \Phi \: d\mu \: = \: \\
& \int_M \nabla \left( \left( - \: \triangle \Phi \: + \: \nabla \Psi \cdot
\nabla \Phi \right) \: \p^\prime \left( \frac{d\mu}{d\nu} \right) \right)
\cdot \nabla \Phi \: \frac{d\mu}{d\nu} \: d\nu \: = \: \notag \\
& \int_M \left( - \: \triangle \Phi \: + \: \nabla \Psi \cdot
\nabla \Phi \right) \: \p^\prime \left( \frac{d\mu}{d\nu} \right)
\: \left( - \: \nabla \cdot \left( \nabla \Phi \: \frac{d\mu}{d\nu} \right) \:
+ \: \nabla \Psi \cdot \nabla \Phi \: \frac{d\mu}{d\nu} \right) 
\: d\nu \: = \: \notag \\
& \int_M \left( - \: \triangle \Phi \: + \: \nabla \Psi \cdot
\nabla \Phi \right)^2 \: \p^\prime \left( \frac{d\mu}{d\nu} \right)
\: \frac{d\mu}{d\nu} \: d\nu \: - \: \notag \\
& \qquad \int_M \left( - \: \triangle \Phi \: + \: \nabla \Psi \cdot
\nabla \Phi \right) \: 
\nabla \Phi \cdot \nabla \p \left( \frac{d\mu}{d\nu} \right)
\: d\nu \: = \notag \\
& \int_M \left( - \: \triangle \Phi \: + \: \nabla \Psi \cdot
\nabla \Phi \right)^2 \: \p^\prime \left( \frac{d\mu}{d\nu} \right)
\: \frac{d\mu}{d\nu} \: d\nu \: + \: \notag \\
& \qquad \int_M \left[ \nabla \left( - \: \triangle \Phi \: + \: \nabla \Psi \cdot
\nabla \Phi \right) \cdot \nabla \Phi \: + \:
\left( - \: \triangle \Phi \: + \: \nabla \Psi \cdot
\nabla \Phi \right) \cdot \triangle \Phi \: - \: \right. \notag \\
& \qquad \: \: \: \: \: \: \: \:
\left. \nabla \Psi \cdot  \left( - \: \triangle \Phi \: + \: \nabla \Psi \cdot
\nabla \Phi \right) \nabla \Phi \right]
\p \left( \frac{d\mu}{d\nu} \right)
\: d\nu \: = \notag \\
& \int_M \left( - \: \triangle \Phi \: + \: \nabla \Psi \cdot
\nabla \Phi \right)^2 \: \p^\prime \left( \frac{d\mu}{d\nu} \right)
\: \frac{d\mu}{d\nu} \: d\nu \: + \: \notag \\
& \qquad \int_M \left[ \nabla \left( - \: \triangle \Phi \: + \: \nabla \Psi \cdot
\nabla \Phi \right) \cdot \nabla \Phi \: - \:
\left( - \: \triangle \Phi \: + \: \nabla \Psi \cdot
\nabla \Phi \right)^2 \:  \right]
\p \left( \frac{d\mu}{d\nu} \right)
\: d\nu. \notag
\end{align} 
 
Combining (\ref{mess2}) and (\ref{mess3}) gives
\begin{align}
\frac{d^2E}{dt^2} \: = \: & 
\int_M \biggl[ \triangle \left( \frac12 |\nabla \Phi|^2
\right) \: - \: 
\nabla \Psi \cdot \nabla \left( \frac12 |\nabla \Phi|^2 \right) 
\: + \:  \\
& \: \: \: \: \: \: \: \:
\nabla \left( - \: \triangle \Phi \: + \: \nabla \Psi \cdot
\nabla \Phi \right) \cdot \nabla \Phi \: - \: \notag \\
& \: \: \: \: \: \: \: \:
 \left( - \: \triangle \Phi \: + \: \nabla \Psi \cdot
\nabla \Phi \right)^2
\biggr]
\: \p \left( \frac{d\mu}{d\nu} \right) \: d\nu \: + \notag \\
& \int_M \left( - \: \triangle \Phi \: + \: \nabla \Psi \cdot
\nabla \Phi \right)^2 \: \p^\prime \left( \frac{d\mu}{d\nu} \right)
\: \frac{d\mu}{d\nu} \: d\nu \notag \\
= \: & \int_M \left[ |\Hess \Phi|^2 \: + \:
\nabla \Phi \cdot (\Ric \: + \: \Hess \Psi) \nabla \Phi
\right]
\: \p \left( \frac{d\mu}{d\nu} \right) \: d\nu \: + \notag \\
& \qquad \int_M \left( - \: \triangle \Phi \: + \: \nabla \Psi \cdot
\nabla \Phi \right)^2 \: \left( - \: \p \left( \frac{d\mu}{d\nu} \right)
\: + \: \p^\prime \left( \frac{d\mu}{d\nu} \right)
\: \frac{d\mu}{d\nu} \right) \: d\nu \notag \\
= \: & \int_M \left[ |\Hess \Phi|^2 \: + \:
\nabla \Phi \cdot (\Ric \: + \: \Hess \Psi) \nabla \Phi
\right]
\: \p \left( \frac{d\mu}{d\nu} \right) \: d\nu \: + \notag \\
& \qquad \int_M \left( - \: \triangle \Phi \: + \: \nabla \Psi \cdot
\nabla \Phi \right)^2 \: \pp \left( \frac{d\mu}{d\nu} \right)
\: d\nu. \notag
\end{align}

In particular, if $U\in \DC_N$ with $n < N$ then from Lemma \ref{p2lem}(b),
\begin{align}
\frac{d^2E}{dt^2} \: \ge \: & \int_M \left[ |\Hess \Phi|^2 \: + \:
\nabla \Phi \cdot (\Ric \: + \: \Hess \Psi) \nabla \Phi
\: - \: \frac{1}{N} \: \left( - \: \triangle \Phi \: + \: \nabla \Psi \cdot
\nabla \Phi \right)^2\right] \: \p \left( \frac{d\mu}{d\nu} \right) \: d\nu \\
\ge \: & \int_M \left[ \frac{1}{n} \: (\triangle \Phi)^2 \: + \:
\nabla \Phi \cdot (\Ric \: + \: \Hess \Psi) \nabla \Phi
\: - \: \frac{1}{N} \: \left( - \: \triangle \Phi \: + \: \nabla \Psi \cdot
\nabla \Phi \right)^2\right] \: 
\p \left( \frac{d\mu}{d\nu} \right) \: d\nu \notag \\
\ge \: & \int_M \left[
\nabla \Phi \cdot (\Ric \: + \: \Hess \Psi) \nabla \Phi
\: - \: \frac{1}{N-n} \: \left(\nabla \Psi \cdot
\nabla \Phi \right)^2\right] \: 
\p \left( \frac{d\mu}{d\nu} \right) \: d\nu \notag \\
= \: & \int_M \left[
\nabla \Phi \cdot \left( \Ric \: + \: 
\Hess \Psi \: - \: \frac{1}{N-n} \: \nabla \Psi
\otimes \nabla \Psi \right) \nabla \Phi
\right] \: \p \left( \frac{d\mu}{d\nu} \right) \: d\nu \notag \\
= \: & \int_M \left[
\nabla \Phi \cdot \Ric_N (\nabla \Phi)
\right] \: \p \left( \frac{d\mu}{d\nu} \right) \: d\nu. \notag
\end{align}

The same conclusion applies if $N=n$ and $\nabla\Psi=0$,
in which case $\Ric_N=\Ric$.

If $\Ric_N \: \ge \: K \: g_M$ then 
\begin{equation} \label{conveqn}
\frac{d^2E}{dt^2} \: \ge \: K \: \int_M 
|\nabla \Phi|^2 \: 
\frac{\p \left( \frac{d\mu}{d\nu} 
\right)}{\frac{d\mu}{d\nu}} \: d\mu.
\end{equation}

If $K = 0$ then (\ref{conveqn}) gives $\frac{d^2E}{dt^2} \: \ge \: 0$.
That is, $E$ is formally convex on $P_2(M)$, no matter what the value of 
$N$ is.

If $N = \infty$ then (\ref{conveqn}) gives
\begin{equation} 
\frac{d^2E}{dt^2} \: \ge \: \left( \inf_{r > 0} K \: \frac{p(r)}{r} \right)
\: \int_M 
|\nabla \Phi|^2 \: d\mu \: = \: 
\lambda(U) \: \int_M 
|\nabla \Phi|^2 \: d\mu.
\end{equation}
As ${\int_M |\nabla \Phi|^2 \: d\mu}$ is the square of
the speed of the geodesic, it follows that $E$ is formally
$\lambda(U)$-convex on $P_2(M)$.

\section{The noncompact case} \label{noncompact}

In the preceding part of the paper we worked with compact spaces $X$.
We now discuss how to adapt our arguments to (possibly noncompact)
pointed metric spaces. To avoid expanding the size of this section
too much, we sometimes simplify the proofs by slightly restricting
the generality of the discussion, and we give details mainly
for the case of nonnegative $N$-Ricci curvature with $N < \infty$.

\subsection{Pointed spaces}

In this section
we will always assume our metric spaces have
distinguished basepoints.
In other words, the objects under study will
be complete pointed metric spaces, which are pairs $(X,\bp)$ where $X$
is a complete metric space and $\bp\in X$. A map $f$ between pointed spaces
$(X_1,\bp_1)$ and $(X_2,\bp_2)$ is said to be a pointed map
if $f(\bp_1)=f(\bp_2)$. In this setting, the analog of
Gromov--Hausdorff convergence is the following:

\begin{definition} \label{pointedGH}
Let $\{(X_i, \bp_i)\}_{i=1}^\infty$ be a sequence of
complete pointed metric spaces. It converges to a
complete pointed metric space $(X, \bp)$ in
the {\em pointed Gromov--Hausdorff topology}, by means of pointed
approximations
$f_i \: : \: X_i \rightarrow X$, if for every $R >0$ there is
a sequence $\{ \epsilon_{R, i} \}_{i=1}^\infty$ of positive numbers
converging to zero so that \\
1. For all $x_i, y_i \in B_{R}(\bp_i)$, we have
$| d_X(f_i(x_i), f_i(y_i)) \: - \: d_{X_i}(x_i, y_i)| \: \le \: \epsilon_{R,i}$.\\
2. For all $x \in B_{R}(\bp)$, there is some $x_i \in 
B_{R}(\bp_i)$ so that $d_X(f_i(x_i), x) \: \le \: \epsilon_{R,i}$.
\end{definition}

A more usual definition would involve approximations defined just on
balls in $X_i$, instead of all of $X_i$. It is convenient for us to
assume that $f_i$ is defined on all of $X_i$, as will be seen when
defining maps between Wasserstein spaces.  The notions of convergence
are equivalent.

Next, a pointed metric-measure space is a complete pointed metric space
$(X, \bp)$ equipped with a nonnegative nonzero Radon measure $\nu$.
We do not assume that $\nu$ has finite mass. In this context,
a pointed map $f \: : \: X_1 \rightarrow X_2$ will be assumed to be Borel,
with the preimage of a compact set being precompact. (This ensures
that the pushforward of a Radon measure is a Radon measure.)

\begin{definition} \label{pointedMGH}
Let $\{(X_i, \bp_i, \nu_i)\}_{i=1}^\infty$ be a sequence of complete
pointed locally compact metric-measure spaces. 
It is said to converge to a complete pointed locally compact
metric-measure space $(X, \bp,\nu)$ in the 
{\em pointed measured Gromov--Hausdorff topology} if 
$\{(X_i, \bp_i)\}_{i=1}^\infty$ converges to $(X, \bp)$ in
the pointed Gromov--Hausdorff topology by means of pointed
approximations $f_i \: : \: X_i \rightarrow X$ 
which additionally satisfy
$\lim_{i \rightarrow \infty} (f_i)_* \nu_i \: = \: \nu$ in
the weak-$*$ topology on $C_c(X)^*$.
\end{definition}

The pointed measured Gromov--Hausdorff topology was used, for
example, in \cite{Cheeger-Colding (2000)}. In what follows we will 
examine its compatibility with the Wasserstein space.

\subsection{Wasserstein space}

If $X$ is a complete pointed metric space, possibly noncompact, 
let $P_2(X)$ be
the space of probability measures $\mu$ on $X$ with a finite
second moment, namely
\begin{equation}
P_2(X) \: = \: \left\{ \mu \in P(X) \: : \: \int_X d(\bp,x)^2 \: d\mu(x) 
\: < \: \infty \right\}.
\end{equation}
One can still introduce the Wasserstein distance
$W_2$ by the same formula as in~\eqref{var}.
Then $W_2$ is a well-defined metric on $P_2(X)$ 
\cite[Theorem~7.3]{Villani (2003)}.
The metric space $(P_2(X),W_2)$ will be called the $2$-Wasserstein space
associated to $X$. It does not depend on the choice of 
basepoint $\bp\in X$.

We will assume that $X$ is a complete separable metric space, 
in which case $P_2(X)$ 
is also a complete separable metric space. Put
\begin{equation}
\bigl(1 + d(\bp,\cdot)^2\bigr) \: C_b(X) \: = \:
\left\{ f \in C(X) \: : \: \sup_{x \in X} \frac{|f(x)|}{1+d(\bp,x)^2} \: < \:
\infty \right\}.
\end{equation}
Then $(1 + d(\bp,\cdot)^2) \: C_b(X)$ is a Banach space with norm
\begin{equation}
\parallel f \parallel \: = \: \sup_{x \in X} \frac{|f(x)|}{1+d(\bp,x)^2},
\end{equation}
and the underlying topological vector space is independent of
the choice of basepoint $\bp$. The dual space
$\left( (1 + d(\bp,\cdot)^2) \: C_b(X) \right)^*$
contains $P_2(X)$ as a subset. As such, $P_2(X)$ inherits a subspace
topology from the weak-$*$ topology on
$\left( (1 + d(\bp,\cdot)^2) \: C_b(X) \right)^*$, which turns out
to coincide with the topology on $P_2(X)$ arising from the metric $W_2$
\cite[Theorem 7.12]{Villani (2003)}.
(If $X$ is noncompact then $P_2(X)$ is not a closed subset
of $\left( (1 + d(\bp,\cdot)^2) \: C_b(X) \right)^*$.)
A subset $S \subset P_2(X)$ is relatively
compact if and only if it satisfies the ``tightness'' condition that
for every $\epsilon > 0$, there is some
$R > 0$ so that for all $\mu \in S$, 
$\int_{X-B_R(\bp)} d(\bp, x)^2 \: d\mu(x) \: < \: \epsilon$
\cite[Theorem 7.12(ii)]{Villani (2003)}. Applying this to a ball
in $P_2(X)$ around $\delta_\bp$, it follows that
if $X$ is noncompact then $P_2(X)$ is not locally compact,
while if $X$ is compact then $P_2(X)$ is compact.

If $X$ is a complete locally compact
length space then for all $k > 0$, $\Lip_k([0,1], X)$ is locally compact, with
the set of geodesic paths between
any two given points in $X$ forming a compact subset.
Then the proof of
Proposition~\ref{Wlengthspace} carries through to show that $P_2(X)$ is
a length space. Finally, if $X$ is pointed then there is a distinguished 
basepoint in $P_2(X)$, namely the Dirac mass $\delta_\bp$. 

The next proposition is an extension of Proposition \ref{everygeod}.

\begin{proposition} \label{everygeod2} Let $(X,\bp)$ be a complete pointed
locally compact length space
and let $\{\mu_t\}_{t \in [0,1]}$ be a geodesic path in $P_2(X)$.
Then there exists some optimal dynamical transference plan $\Pi$
such that $\{\mu_t\}_{t \in [0,1]}$ 
is the displacement interpolation associated to $\Pi$.
\end{proposition}
\begin{proof}
We can go through the proof of Proposition \ref{everygeod}, constructing
the probability measures $R^{(i)}$ with support on the locally
compact space $\Gamma$. For each $i$, we have
$(e_0)_* R^{(i)} \: = \: \mu_0$ and $(e_1)_* R^{(i)} \: = \: \mu_1$.
In order
to construct a weak-$*$ accumulation point 
$R^{(\infty)}$, as a probability measure
on $\Gamma$, it suffices to show
that for each $\epsilon > 0$ there is a compact set $K \subset
\Gamma$ so that for all $i$, $R^{(i)}(\Gamma - K) \: < \: 
\epsilon$.

Let
$E \: : \: \Gamma \rightarrow X \times X$ be the endpoints map. 
Given $r > 0$, put
$K \: = \: E^{-1}(\overline{B_r(\bp)} \times \overline{B_r(\bp)})$,
a compact subset of $\Gamma$. 
As 
\begin{equation}
\Gamma - K \: = \: E^{-1}
\left( \bigl((X - \overline{B_r(\bp)}) \times X\bigr) \cup
\bigl(X \times (X - \overline{B_r(\bp)})\bigr) \right),
\end{equation}
we have
\begin{align}
R^{(i)}(\Gamma - K) \: & \le \:
(E_* R^{(i)}) ((X - \overline{B_r(\bp)}) \times X) \: + \:
(E_* R^{(i)}) (X \times (X - \overline{B_r(\bp)})) \\
& = \:
\mu_0 (X - \overline{B_r(\bp)}) \: + \:
\mu_1 (X - \overline{B_r(\bp)}). \notag
\end{align}
Taking $r$ sufficiently large, we can ensure that
$\mu_0(X - \overline{B_r(\bp)}) \: + \: 
\mu_1(X - \overline{B_r(\bp)}) \: < \:
\epsilon$. 
\end{proof}

Using Proposition \ref{everygeod2}, we show that geodesics with
endpoints in a given compact subset of $P_2(X)$ will all lie in a compact
subset of $P_2(X)$.

\begin{proposition} \label{contained}
For any compact set $K \subset P_2(X)$, there is a compact set
$K^{\prime} \subset P_2(X)$ with the property 
that for any $\mu_0, \mu_1 \in K$,
if $\{\mu_t\}_{t \in [0,1]}$ is a Wasserstein
geodesic from $\mu_0$ to $\mu_1$ then $\mu_t \in K^\prime$ for all
$t \in [0,1]$. 
\end{proposition}
\begin{proof}
Given $\mu_0, \mu_1 \in K$, let 
$\{\mu_t\}_{t \in [0,1]}$ be a Wasserstein geodesic from
$\mu_0$ to $\mu_1$.
Then
\begin{align} \label{breakupp}
\int_{X - B_R(\bp)} d(\bp, x)^2 \:
\: d\mu_t(x) \: & = \:
\int_\Gamma d(\bp, \gamma(t))^2 \: 1_{\gamma(t) \notin B_R(\bp)}
\: d\Pi(\gamma) \\
& = \:
\int_\Gamma d(\bp, \gamma(t))^2 \: 1_{\gamma(t) \notin B_R(\bp)} \:
1_{\max(d(\bp, \gamma(0)),d(\bp, \gamma(1)))
\: \ge \: \frac{R}{2}}
\: d\Pi(\gamma), \notag
\end{align}
where $\Pi$ comes from Proposition \ref{everygeod2}.
We break up the integral in the last term of (\ref{breakupp}) into two
pieces according to whether $d(\bp, \gamma(0)) \: \le \: d(\bp, \gamma(1))$
or $d(\bp, \gamma(1)) \: < \: d(\bp, \gamma(0))$. If
$d(\bp, \gamma(0)) \: \le \: d(\bp, \gamma(1))$ then
\begin{align}
d(\bp, \gamma(t)) \: & \le \: d(\bp, \gamma(0)) \: + \: d(\gamma(0), \gamma(t)) \:
\le \: d(\bp, \gamma(0)) \: + \: d(\gamma(0), \gamma(1)) \\
& \le \: 2 \: d(\bp, \gamma(0)) \: + \: d(\bp, \gamma(1)) \: \le \:
3 \: d(\bp, \gamma(1)). \notag
\end{align}
Then the contribution to the last term of (\ref{breakup}), when
$d(\bp, \gamma(0)) \: \le \: d(\bp, \gamma(1))$, is bounded above by
\begin{equation}
9 \: \int_\Gamma d(\bp, \gamma(1))^2 \: 
1_{d(\bp, \gamma(1))
\: \ge \: \frac{R}{2}}
\: d\Pi(\gamma) \: = \: 9 \: 
\int_{X - B_{R/2}(\bp)} d(\bp, x)^2 \:  d\mu_1(x).
\end{equation}
Adding the contribution when 
$d(\bp, \gamma(1)) \: < \: d(\bp, \gamma(0))$ gives
\begin{equation} \label{contribution}
\int_{X - B_R(\bp)} d(\bp, x)^2 \:
\: d\mu_t(x) \: \le \: 9 \: \int_{X - B_{R/2}(\bp)} d(\bp, x)^2 \:  d\mu_0(x)
\: + \: 9 \: \int_{X - B_{R/2}(\bp)} d(\bp, x)^2 \:  d\mu_1(x).
\end{equation}
For any $\epsilon > 0$ we can choose $R > 0$ so that the
right-hand-side of (\ref{contribution}) is bounded above by $\epsilon$,
uniformly in $\mu_0, \mu_1 \in K$. This proves the
proposition.
\end{proof}

\begin{remark}
Although we will consider optimal transport between elements of $P_2(X)$,
there are also interesting issues concerning
``optimal transport'' in a generalized sense, with possibly infinite cost, on
the whole of $P(X)$. For example, one has McCann's theorem
about existence of ``generalized optimal transport'' 
between arbitrary probability measures on $\R^n$ 
\cite[Theorem~2.32]{Villani (2003)}. 
\end{remark}

\subsection{Functionals}

In the nonpointed part of the paper we dealt with a compact measured 
length space $(X, d, \nu)$, with the background measure $\nu$ lying in 
$P_2(X)$. When generalizing 
to the case when $X$ is pointed and possibly noncompact,
one's first inclination might be to again have 
$\nu$ lie in $P_2(X)$.
This is indeed the appropriate choice for some purposes, such as to
extend the functional analytic results of Sections \ref{HWI inequalities} and 
\ref{secineq}. However, requiring $\nu$ to lie in $P_2(X)$ would rule out
such basic cases as $X \: = \: \R^N$ with the Lebesgue measure.
Additionally, it would preclude the tangent cone construction
for a compact space with $N$-Ricci curvature bounded below.
Because of this, in what follows 
we will allow $\nu$ to have infinite mass, at the
price of some additional complications

Let $\nu$ be a nonnegative nonzero Radon measure on $X$.
Let $U$ be a continuous
convex function on $[0, \infty)$ with $U(0) \: = \: 0$.
One would like to define the functional $U_\nu$ as in 
Definition \ref{deffunctional}, but this requires some care.
Even if  
we use (\ref{defUnu}) to define $U_\nu(\mu)$ for $\mu \: = \: \rho \nu$ and
$\rho \in C_c(X)$, in general there is no way to extend $U_\nu$
to a lower semicontinuous function on $P_2(X)$. A way to circumvent
this difficulty is to impose a growth assumption on $\nu$. 

\begin{definition}
For $k>0$, we define $M_{-k}(X)$ to be the space of nonnegative
Radon measures $\nu$ on $X$ such that
\begin{equation}
\int_X (1 + d(\bp, x)^2)^{- \: \frac{k}{2}} \: d\nu(x) \: < \: \infty.
\end{equation}
Equivalently, $\nu$ is a nonnegative Radon measure on $X$ that
lies in the dual space of $(1 + d(\bp,\cdot)^2)^{- \: \frac{k}{2}} \: C_b(X)$.
We further define $M_{-\infty}(X)$ by the condition
$\int_X e^{- \: c \: d(x, *)^2} \: d\nu(x) \: < \: \infty$,
where $c$ is a fixed positive constant.
\end{definition}

\begin{proposition} Let $X$ be an arbitrary metric space.
Given $N\in [1,\infty]$, suppose that $U\in \DC_N$ and $\nu\in M_{-2(N-1)}(X)$.
Then $U_\nu$ is a well-defined functional on $P_2(X)$, 
with values in $\R\cup\{\infty\}$.
\end{proposition}

\begin{proof}
Suppose first that $N < \infty$.
From the definition of $\DC_N$, there is
a constant $A > 0$ so that
\begin{equation}
\lambda^N \:  U(\lambda^{-N}) \: \ge \:
- \: A\lambda \: - \: A,
\end{equation}
so
\begin{equation} \label{Ubound}
U(\rho) \: \ge \: - \: A \: 
\left( \rho \: + \: \rho^{1 - \frac{1}{N}}
\right).
\end{equation}

Of course, $\rho$ lies in $L^1(X,\nu)$. To prove that
$U_\nu(\mu)$ is well-defined, it suffices to show that
$\rho^{1-1/N}$ also lies in $L^1(X,\nu)$. For that we use
H\"older's inequality:
\begin{align}
\int_X \rho(x)^{1 - \frac{1}{N}} \: d\nu(x) 
& = \int_X \bigl((1+d(\bp,x)^2) \rho(x)\bigr)^{1 - \frac{1}{N}} \:
(1+d(\bp,x)^2)^{- 1 + \frac{1}{N}} \: d\nu(x) \\
& \leq \left( \int_X (1+d(\bp,x)^2) \rho(x) \: d\nu(x) 
\right)^{1 - \frac{1}{N}} 
\left( \int_X (1+d(\bp,x)^2)^{- (N-1)} \: d\nu(x) 
\right)^{\frac{1}{N}}. \notag
\end{align}

Now suppose that $N = \infty$. From Lemma \ref{increasing}, 
if $U$ is nonlinear then there are constants $a,b > 0$ so that
\begin{equation}
U(\rho) \: \ge \: a \:  \rho \: \log \rho \: - \: b \:  \rho.
\end{equation}
Thus it is sufficient to show that $(\rho\log\rho)_-\in L^1(X,\nu)$.
Applying Jensen's inequality with the probability measure
$\frac{e^{-c\, d(\bp,\cdot)^2} \: d\nu}
{\int_X e^{-c\, d(\bp,\cdot)^2} \: d\nu}$
gives
\begin{align}
& \int_X \rho(x) \: \log(\rho(x)) \: d\nu(x) \: = \\
& 
\int_X \rho(x) \: e^{c d(\bp, x)^2} \: \log \left( \rho(x) \: e^{c d(\bp, x)^2}
\right) \: e^{- c d(\bp, x)^2} \: d\nu(x) \: - 
c \int_X d(\bp, x)^2 \: \rho(x) \: d\nu(x) \: = \notag \\
& \left( \int_X e^{- c d(\bp, x)^2} \: d\nu(x) \right)
\left( 
\int_X \rho(x) \: e^{c d(\bp, x)^2} \: \log \left( \rho(x) \: e^{c d(\bp, x)^2}
\right) \: \frac{e^{- c d(\bp, x)^2} \: d\nu(x)}{\int_X e^{- c d(\bp, x)^2} \: 
d\nu(x)} \right) \: - \notag \\
& \: \: \: \: \: \: c \:
\int_X d(\bp, x)^2 \: \rho(x) \: d\nu(x) \: \ge \: \notag \\
& \left( \int_X e^{- c d(\bp, x)^2} \: d\nu(x) \right) 
\left( \frac{\int_X \rho \: d\nu}{\int_X e^{- c d(\bp, x)^2} \: 
d\nu(x)} \right) 
\log \left( \frac{\int_X \rho \: d\nu}{\int_X e^{- c d(\bp, x)^2} \: 
d\nu(x)} \right) 
\: - \notag \\
& \: \: \: \: \: \:  c \: \int_X d(\bp, x)^2 \: \rho(x) \: d\nu(x). \notag
\end{align}
This concludes the argument.
\end{proof}

\subsection{Approximation arguments} \label{approxarg}

Now we check that the technical results in 
Appendices~\ref{applsc} and~\ref{secmoll} go through to the pointed
locally compact case. 

There is an obvious way to generalize the conclusion 
of Theorem~\ref{thmdual} by introducing the quantity
\begin{equation} \label{exp}
\sup_{\varphi\in (1 + d(\bp,\cdot)^2) \: L^\infty(X), \> 
\varphi\leq U'(\infty)} 
\left( \int_X \varphi\,d\mu - \int_X U^*(\varphi) \,d\nu.
\right )
\end{equation}
We claim that this quantity is not $-\infty$ as long as $\nu$ lies
in $M_{-2(N-1)}(X)$. To prove this, it
suffices to exhibit a $\varphi$ such that $\int\varphi\,d\mu>-\infty$
and $\int U^*(\varphi) \,d\nu <\infty$. It turns out that
$\varphi(x)=-c\, d(\bp,x)^2$ will do the job, taking into account
the following lemma:

\begin{lemma} \label{Ustar}
If $U \in \DC_N$ with $N < \infty$ then as $p \rightarrow -\infty$, 
\begin{equation} \label{negp}
U^*(p) \: =  \: O \left( (-p)^{1-N} \right).
\end{equation}
If $U \in \DC_\infty$ then as $p \rightarrow -\infty$, 
\begin{equation}
U^*(p) \: = \: O \left( e^{p} \right).
\end{equation}
\end{lemma}

\begin{proof}
Suppose first that $U \in \DC_N$ with 
$N < \infty$.
Then for $p$ sufficiently negative, using (\ref{Ubound}) we have
\begin{equation}
0 \: \le \: U^*(p) \: = \: 
\sup_{r \ge 0} \left[ pr \: - \: U(r)
\right] \: \le \: 
\sup_{r \ge 0} \left[ pr \: + 2A \: r^{1 - \frac{1}{N}}
\right] \: = \:  \const \: (-p)^{1-N}. 
\end{equation}
The case $N = \infty$ is similar.
\end{proof}

The analog of Theorem \ref{thmdual} becomes
\begin{align}
U_\nu(\mu) & = 
\sup_{\varphi\in (1 + d(\bp, \cdot)^2) \:
L^\infty(X), \> \varphi\leq U'(\infty)} 
\left( \int_X \varphi\,d\mu - \int_X U^*(\varphi) \,d\nu
\right ) \\ 
& =  \sup_{M \in \Z^+}
\sup_{\varphi\in (1+d(\bp, \cdot)^2)C_b(X), \> \varphi\leq U'(M)} 
\left( \int_X \varphi\,d\mu - 
\int_X U^*(\varphi) \,d\nu \right). \notag
\end{align}
The proof is similar to the proof of Theorem \ref{thmdual}, with the following
modifications.
Given $R > 0$ and $M \in \Z^+$ large numbers, we
construct $\overline{\phi}_{R,M}$ on $B_R(\bp)$ as in the proof of
Theorem \ref{thmdual}. We extend $\overline{\phi}_{R,M}$ to $X$ by
setting it equal to $- \: d(\bp, \cdot)^2$ on $X - B_R(\bp)$.  Then
after passing to an appropriate subsequence of $\{\overline{\phi}_{R,M}\}$, 
the analog of (\ref{goalUU}) holds. As in the proof of Lusin's theorem,
we can find a sequence $\{\psi_k\}$ of uniformly bounded continuous functions 
that converges pointwise to
$\frac{\overline{\phi}_{R,M}}{1+d(\bp,\cdot)^2}$. Considering the function
$(1 + d(\bp, \cdot)^2) \: \psi_k$ for large $k$ proves the theorem.

Then Theorem \ref{thmlsc}(i) extends, where $\mu$ lies in 
$P_2(X)$ and $\nu$ is a measure on $X$ that lies in the dual space of
$(1 + d(\bp,\cdot)^2)^{-(N-1)} \: C_b(X)$, which we endow with the
weak-$*$ topology.

Theorem~\ref{thmlsc}(ii) is a bit more subtle because we have to be careful
about how the pushforward map $f_*$ acts on the measures at spatial infinity.
The discussion is easier when $N<\infty$, so from now on we restrict 
to this case. Then the statement in Theorem~\ref{thmlsc}(ii) 
goes through as soon as we impose that the map $f$ is a pointed
Borel map satisfying
\begin{equation} \label{condgrowthf}
A^{-1} \: d_X(\bp_X, x) - A \: \le \:
d_Y(\bp_Y, f(x)) \: \le \: A \: d_X(\bp_X, x) + A
\end{equation} 
for some constant $A>0$. This condition ensures that $f_*$ maps 
$P_2(X)$ to $P_2(Y)$, and maps measures on $X$ that lie in the dual space of
$(1 + d(\bp_X, \cdot)^2)^{-(N-1)} \: C_b(X)$ to measures on $Y$ that lie 
in the dual space of $(1 + d(\bp_Y,\cdot)^2)^{-(N-1)} \: C_b(Y)$.
\med

Finally, we wish to extend Theorem~\ref{thmapproxU} to the pointed locally
compact setting, replacing
$C(X)$ by $C_c(X)$. 
Given $\delta > 0$, 
let $\{x_j\}$ be a maximal $\frac{\delta}{2}$-separated net in $X$.
Then ${\mathcal C} \: = \: \{B_{\delta}(x_j)\}$ is an open cover of $X$. 
It is locally finite, as $X$ is
locally compact. If $\{\phi_j\}$ is a subordinate partition of unity then we 
define the operator $K_{\mathcal C}$ as in Appendix~\ref{secmoll}. 
Given $\mu \in P_2(X, \nu)$,
we claim that $K_{\mathcal C} \mu \in P_2(X, \nu)$. To see this, we write
\begin{align}
\int_X (1+ d(\bp, x)^2) \: dK_{\mathcal C}\mu(x) \: & = \:
\int_X (1+ d(\bp, x)^2) \:
\sum_{j \: : \: \int_X \phi_j \: d\nu > 0}  
\phi_j(x) \: \frac{\int_X \phi_j \: d\mu}{\int_X \phi_j \: d\nu} \: d\nu(x) \\
& = \: \int_X 
\sum_{j \: : \: \int_X \phi_j \: d\nu > 0} 
\frac{\int_X (1+d(\bp, \cdot)^2) \: \phi_j \: d\nu}{\int_X \phi_j \: d\nu}
\: \phi_j(y) \: d\mu(y). \notag
\end{align}
From the choice of $\{x_j\}$ and $\{\phi_j\}$, 
there is a constant $C < \infty$ so that
\begin{equation}
\frac{\int_X (1+d(\bp, \cdot)^2) \: \phi_j \: d\nu}{\int_X \phi_j \: d\nu} \le \: C \:
(1+d(\bp, x_j)^2)
\end{equation}
for all $j$ with $\int_X \phi_j \: d\nu > 0$. There is another constant $C^\prime < \infty$ so that
\begin{equation}
(1+d(\bp, x_j)^2) \: \phi_j(y) \: \le \: C^\prime \: (1+d(\bp, y)^2) \: \phi_j(y)
\end{equation}
for all $j$ and all $y \in X$, from which the claim follows.

Next, we claim that $U_\nu(K_{\mathcal C} \mu) \: \le \: U_\nu(\mu)$.  We use the fact that
$K_{\mathcal C}\mu$ is the product of a continuous function on $X$ with $\nu$. As in
(\ref{comb0})-(\ref{comb2}), for each $R > 0$ we have
\begin{equation}
\int_{B_R(\bp)} U(K_{\mathcal C}\mu(x)) \: d\nu(x) \: \le \: 
\int_{B_R(\bp)} U(K_{\mathcal C}\rho) \: d\nu \: + \: U^\prime(\infty) \: \int_{B_R(\bp)}
K_{\mathcal C}\mu_s(x) \: d\nu(x).
\end{equation}
Taking $R \rightarrow \infty$ and applying the arguments of the
particular cases in the proof of Theorem \ref{thmapproxU} gives the claim.

For $R>0$ and $\epsilon > 0$, let $\phi_{R, \epsilon} \: : \: X \rightarrow [0,1]$
be a continuous function which is one on $B_R(\bp)$ and vanishes outside of
$B_{R+\epsilon}(\bp)$. We can find sequences $\{\delta_k\}$,
$\{{\mathcal C}_k\}$, $\{R_k\}$ and
$\{\epsilon_k\}$ so that 
$\lim_{k \rightarrow \infty} \frac{\phi_{R_k, \epsilon_k } 
K_{{\mathcal C}_k} \mu}{\int_X \phi_{R_k, \epsilon_k } 
d(K_{{\mathcal C}_k} \mu)} \: = \: 
\mu$ in $P_2(X, \nu)$ and
\begin{equation}
\limsup_{k \rightarrow \infty} \: U_\nu \left( 
\frac{\phi_{R_k, \epsilon_k } 
K_{{\mathcal C}_k} \mu}{\int_X \phi_{R_k, \epsilon_k } 
d(K_{{\mathcal C}_k} \mu)}\right)
\: \le \: \limsup_{k \rightarrow \infty} \: U_\nu(K_{{\mathcal C}_k} \mu).
\end{equation}
From the previously-shown lower semicontinuity of $U_\nu$,
we know that
\begin{equation}
U_\nu(\mu) \: \le \: \liminf_{k \rightarrow \infty} \:
U_\nu \left( \frac{\phi_{R_k, \epsilon_k } 
K_{{\mathcal C}_k} \mu}{\int_X \phi_{R_k, \epsilon_k } 
d(K_{{\mathcal C}_k} \mu)}
\right).
\end{equation}
We have already show that
$U_\nu(K_{{\mathcal C}_k} \mu) \: \le \: U_\nu(\mu)$. Hence
$U_\nu(\mu) \: = \: \lim_{k \rightarrow \infty} 
U_\nu \left( \frac{\phi_{R_k, \epsilon_k } 
K_{{\mathcal C}_k} \mu}{\int_X \phi_{R_k, \epsilon_k } 
d(K_{{\mathcal C}_k} \mu)}
\right)$. This proves the desired
extension of Theorem \ref{thmapproxU}.

\subsection{Stability of $N$-Ricci curvature bounds}
\label{ppmghc}

We now define the notion of a complete pointed measured locally compact
length space $(X, \bp, \nu)$ having nonnegative $N$-Ricci
curvature as in Definition \ref{nonnegdef}, provided that
$\nu\in M_{-2(N-1)}(X)$. Note that this notion is independent of 
the choice of basepoint.

Most of the geometric inequalities discussed in 
Sections~\ref{HWI inequalities}, \ref{secricci} 
and~\ref{secineq} have evident extensions
to the pointed case. When discussing HWI, log Sobolev, Talagrand and
Poincar\'e inequalities we assume that $\nu\in P_2(X)$.
If $X$ is a smooth Riemannian manifold such that the
reference measure $\nu$ lies in $M_{-2(N-1)}(X)$ then there is an
analog of Theorem~\ref{Riccimfolds}, expressing the condition of
having nonnegative $N$-Ricci curvature in terms of classical tensors.

\begin{remark}
For $n > 2$, if $X=\R^n$ is endowed with the
Lebesgue measure $\nu$ then
$\nu \in M_{-2(n-1)}(\R^n)$ and 
$(X, \nu)$ will
have nonnegative $n$-Ricci curvature.
In case of $X \: = \: \R^2$, endowed with the Lebesgue measure $\nu$,
it is not true that $\nu \in M_{-2(n-1)}(\R^n)$. The borderline case
$n=2$ merits further study; see also Corollary \ref{tangentcone}.
\end{remark}

The issue of showing that Ricci bounds are preserved by
pointed measured Gromov--Hausdorff convergence is more involved than
in the nonpointed case.
The following definition seems to be useful.

\begin{definition} \label{properpointedGH}
A sequence $\{(X_i, \bp_i)\}_{i=1}^\infty$ of pointed metric spaces 
converges to the pointed metric space $(X, \bp)$ in
the {\em proper} pointed Gromov--Hausdorff topology if \\
1. It converges
in the pointed Gromov--Hausdorff topology, by means of pointed 
approximations
$f_i \: : \:  X_i \rightarrow X$, \\
2. There is a function $\widehat{R} \: : (0, \infty) \rightarrow (0,\infty)$ with
$\widehat{R}(R) \: > \: R$ for all $R$, \\
3. There are nondecreasing functions
$G_i \: :  \: (0, \infty) \rightarrow (0, \infty)$, each increasing linearly at
infinity,  and \\
4. There is a constant $A > 0$ \\
such that \\
1. For all $R >0$, we have $R \:< \: \liminf_{i \rightarrow \infty} G_i(\widehat{R}(R))$, and \\
2. For all $x_i \in X_i$,
\begin{equation} \label{Gidi}
G_i(d_i(\bp_i, x_i)) \: \le \: d(\bp, f_i(x_i)) \: \le \:
A \: d_i(\bp_i, x_i) \: + \: A.
\end{equation}
\end{definition}

Here are the main motivations for the definition.
The condition that $d(\bp, f_i(x_i)) \: \le \:
A \: d_i(\bp_i, x_i) \: + \: A$ ensures that
$(f_i)_*$ maps $P_2(X_i)$ to $P_2(X)$.
Condition 3. and~\eqref{Gidi}
ensure that $(f_i)_*$ maps a measure on $X_i$ lying in 
the dual space of
$(1 + d_i(\bp_i, \cdot)^2)^{-(N-1)} \: C_b(X_i)$ to a
measure on $X$ lying in the dual space of
$(1 + d(\bp, \cdot)^2)^{-(N-1)} \: C_b(X)$.
The condition
$G_i(d_i(\bp_i, x_i)) \: \le \: d(\bp, f_i(x_i))$ implies that
$f_i$ is metrically proper; for example,
it cannot map an unbounded sequence in $X_i$ to a bounded
sequence in $X$.
The conditions
$R \: < \: \liminf_{i \rightarrow \infty} G_i(\widehat{R}(R))$ and
$G_i(d_i(\bp_i, x_i)) \: \le \: d(\bp, f_i(x_i))$ imply that for any $R > 0$, we have
$f_i^{-1}(B_R(\bp)) \subset B_{\widehat{R}(R)}(\bp_i)$
for sufficiently large $i$. It then follows that in fact,
$f_i^{-1}(B_R(\bp)) \subset B_{R+\epsilon_{R,i}}(\bp_i)$ for large $i$.

\begin{definition}
A sequence of pointed metric-measure spaces 
$\{(X_i, \bp_i, \nu_i)\}_{i=1}^\infty$
converges to $(X, \bp, \nu)$ in the proper
pointed $M_{-k}$-measured Gromov--Hausdorff topology
if $\lim_{i \rightarrow \infty} (X_i, \bp_i) \: = \: (X, \bp)$ 
in the proper pointed Gromov--Hausdorff topology, by means of
pointed Borel 
approximations $f_i \: : \: X_i \rightarrow X$ as above, and
in addition $\lim_{i \rightarrow \infty} (f_i)_* \nu_i \: = \:
\nu$ in the weak-$*$ topology on the dual space of
$(1 + d(\bp,\cdot)^2)^{- \: \frac{k}{2}} \: C_b(X)$.
\end{definition}

Now we can prove the stability of Ricci curvature bounds with
respect to the proper pointed measured Gromov--Hausdorff topology.
Again, for simplicity 
we restrict to the case of nonnegative $N$-Ricci curvature
with $N < \infty$.

\begin{theorem} \label{properlimits}
Let $\{(X_i, \bp_i, \nu_i)\}_{i=1}^\infty$ be a sequence of
complete pointed measured locally compact 
length spaces with $\lim_{i \rightarrow \infty} 
(X_i, \bp_i,\nu_i) \: = \: (X, \bp, \nu_\infty)$ in the
proper pointed $M_{-2(N-1)}$-measured Gromov--Hausdorff topology.
If each $(X_i, \nu_i)$ has nonnegative $N$-Ricci curvature then
$(X, \nu_\infty)$ has nonnegative $N$-Ricci curvature.
\end{theorem}
\begin{proof}
Given
$\mu_0, \mu_1 \in P_2(X, \nu_\infty)$, we wish to show that there is a
geodesic joining them along which (\ref{defldc}) holds for $U_{\nu_\infty}$,
with $\lambda = 0$.

We first show that the claim is true if $\mu_0 \: =
\rho_0 \: \nu_\infty$ and
$\mu_1 \: = \: \rho_1 \: \nu_\infty$, with $\rho_0$ and $\rho_1$
being compactly-supported continuous functions on $X$. 
We will follow the proof of Theorem~\ref{thmstabGH}. This involved constructing
a limiting geodesic using the Arzel\`a--Ascoli
theorem, which in turn used the compactness of $P_2(X)$.
If $X$ is noncompact then $P_2(X)$ is not locally compact. Nevertheless,
we will show that the needed arguments can be carried out in a
compact subset of $P_2(X)$.

By assumption, there
is some $R > 0$ so that $\rho_0$ and $\rho_1$ have support in
$B_R(\bp)$. Put
$\mu_{i,0} \: = \: (f_i^* \rho_0) \: \nu_i$ and
$\mu_{i,1} \: = \: (f_i^* \rho_1) \: \nu_i$. From the definition
of proper pointed Gromov--Hausdorff convergence, for large $i$ we know that
$f_i^* \rho_0$ and $f_i^* \rho_1$ have support in
$B_{R + \epsilon_{R,i}}(\bp_i)$.
Choose Wasserstein geodesics $c_i$ as in the proof of Theorem \ref{thmstabGH}.
If $\gamma$ is a geodesic joining two points of $B_{R + \epsilon_{R,i}}(\bp_i)$
then $\gamma([0,1]) \subset B_{2R + 2\epsilon_{R,i}}(\bp_i)$, 
so Proposition \ref{everygeod2} implies that each $c_i(t)$ has support in
$B_{2R + 2\epsilon_{R,i}}(\bp_i)$. Then
$(f_i)_*(c_i(t))$ has support in $B_{2R + 3\epsilon_{R,i}}(\bp)$.

Hence for large $i$, each measure $(f_i)_*(c_i(t))$ has support
in $B_{2R + 1}(\bp)$. As the elements of
$P_2(X)$ with support in $B_{2R + 1}(\bp)$ form a relatively compact
subset of $P_2(X)$, we can now carry out the arguments of the proof of
Theorem~\ref{thmstabGH}.

This proves the theorem when $\mu_0, \mu_1 \in P_2(X, \nu_\infty)$ 
have compactly-supported
continuous densities. To handle the general case, we will use the
arguments of Proposition \ref{lemsuffdc}. Again, the main issue is
to show that one can carry out the arguments in a compact subset of
$P_2(X)$.

Let $r_0 > 0$ be such that
$\mu_0(B_{r_0}(\bp)) > 0$ and $\mu_1(B_{r_0}(\bp)) > 0$. 
For $r > r_0$, put $\mu_{0,r} \: = \: 
\frac{1_{B_r(\bp)}}{\mu_0(B_r(\bp))} \: \mu_0$ and $\mu_{1,r} \: = \: 
\frac{1_{B_r(\bp)}}{\mu_1(B_r(\bp))} \: \mu_1$. 
Let
$\{\mu_{\delta,0,r}\}$ and $\{\mu_{\delta,1,r}\}$ be
mollifications of $\mu_{0,r}$ and $\mu_{1,r}$, respectively,
using a maximal $\delta$-separated net as discussed in
Section \ref{approxarg}.
Then
\begin{align}
\int_{X - B_R(\bp)} d(\bp, x)^2 \: d\mu_{\delta,0,r}(x) \: & = \:
\frac{1}{\mu_0(B_r(\bp))}
 \int_{X - B_R(\bp)} d(\bp, x)^2 \: dK_{\mathcal C} (1_{B_r(\bp)} \mu_0)(x) \\
& \le \:
\frac{1}{\mu_0(B_r(\bp))} \int_{X - B_{R-10\delta}(\bp)} 
\left( d(\bp, x) \: + \:
10\delta
\right)^2 \: 1_{B_r(\bp)}(x) \:  d\mu_0(x) \notag \\
& \le \:
\frac{1}{\mu_0(B_{r_0}(\bp))} \int_{X - B_{R-10\delta}(\bp)} 
\left( d(\bp, x) \: + \:
10\delta
\right)^2 \: d\mu_0(x). \notag
\end{align}
For small $\delta$, we obtain
\begin{equation} \label{rhs1}
\int_{X - B_R(\bp)} d(\bp, x)^2 \: d\mu_{\delta,0,r}(x) \: \le \:
\frac{2}{\mu_0(B_{r_0}(\bp))} \int_{X - B_{R/2}(\bp)} 
d(\bp, x)^2 \: d\mu_0(x).
\end{equation}
Similarly,
\begin{equation} \label{rhs2}
\int_{X - B_R(\bp)} d(\bp, x)^2 \: d\mu_{\delta,1,r}(x) \: \le \:
\frac{2}{\mu_1(B_{r_0}(\bp))} \int_{X - B_{R/2}(\bp)} 
d(\bp, x)^2 \: d\mu_1(x).
\end{equation}
As the right-hand-sides of (\ref{rhs1}) and (\ref{rhs2}) can be
made arbitrarily small by taking $R$ sufficiently large, it 
follows that $\bigcup_{r > r_0} \bigcup_{i=1}^\infty
\{ \mu_{i^{-1},0,r}, \mu_{i^{-1},1,r} \}$ is relatively compact in $P_2(X)$.
With an appropriate choice of $i_j$ for $j$ large, we have
$\lim_{j \rightarrow \infty} \mu_{i_j^{-1},0,j} \: = \: \mu_0$ and
$\lim_{j \rightarrow \infty} \mu_{i_j^{-1},1,j} \: = \: \mu_1$. 
Using Proposition \ref{contained},
the argument in the proof of Proposition \ref{lemsuffdc} can now be
applied to show that there is a geodesic from $\mu_0$ to $\mu_1$ 
along which (\ref{defldc}) holds for $U_{\nu_\infty}$,
with $\lambda = 0$.
\end{proof}

\subsection{Tangent Cones}

We now give an application of Theorem \ref{properlimits} 
that just involves the {\em pointed measured Gromov--Hausdorff topology}
that was introduced in Definition \ref{pointedMGH}.

\begin{corollary} \label{tangentcone}
Let $\{(X_i, \bp_i, \nu_i)\}_{i=1}^\infty$ be a sequence of
complete pointed measured locally compact
length spaces.
Suppose that $\lim_{i \rightarrow \infty} 
(X_i, \bp_i,\nu_i) \: = \: (X, \bp, \nu)$ in the
pointed measured Gromov--Hausdorff topology of
Definition \ref{pointedMGH}. If $N \in (2, \infty)$ and 
each $(X_i, \nu_i)$ has nonnegative $N$-Ricci curvature then
$(X, \nu)$ has nonnegative $N$-Ricci curvature.
\end{corollary}

\begin{proof}[Proof of Corollary~\ref{tangentcone}]
If $X$ is compact then the result follows from
Theorem \ref{stabRicci}, so we will assume that $X$ is noncompact.
Let $\{ f_i \}_{i=1}^\infty$ be a sequence of approximations
as in Definition \ref{pointedGH}. Given $R_i >0$, let 
$\widehat{f}_i \: : \: X_i
\rightarrow X$ be an arbitrary Borel map such that
$\widehat{f}_i(x_i) \: = \: f_i(x_i)$
if $d_i(\bp_i, x_i) \: < \: R_i$ and
$d(\bp, \widehat{f}_i(x_i)) \: = \: d(\bp_i, x_i)$ if
$d_i(\bp_i, x_i) \: \ge \: R_i$. (For example, if
$\gamma \: : \: [0, \infty) \rightarrow X$ is a normalized ray
with $\gamma(0) \: = \: \bp$ then we can put
$\widehat{f}_i(x_i) \: = \: \gamma(d(\bp_i, x_i))$ when
$d_i(\bp_i, x_i) \: \ge \: R_i$.) After passing to a subsequence
of $\{f_i\}_{i=1}^\infty$ (which we relabel as $\{f_i\}_{i=1}^\infty$)
and replacing $f_i$ by
$\widehat{f}_i$ for an appropriate choice of $R_i$, we can
assume that $\lim_{i \rightarrow \infty} 
(X_i, \bp_i,\nu_i) \: = \: (X, \bp, \nu)$ in the
{\em proper} pointed measured Gromov--Hausdorff topology,
with $\widehat{R}(R) \: = \: 3R$ and $G_i(r) \: = \: \frac{r}{2}$.

As each $(X_i, \nu_i)$ has nonnegative $N$-Ricci curvature, and
$\lim_{i \rightarrow \infty} (\widehat{f}_i)_* \nu_i \: = \: \nu$ in
$C_c(X)^*$, the Bishop--Gromov inequality of Theorem \ref{thmBG}
(as extended to the noncompact case) implies that there are 
constants $C, r_0 > 0$ so that for all $i$, whenever
$r \ge r_0$ we have $\nu_i(B_r(\bp_i)) \: \le \: C \: r^N$.
As $N > 2$, it follows from dominated convergence that 
$\lim_{i \rightarrow \infty} (\widehat{f}_i)_* \nu_i \: = \: \nu$ in
the weak-$*$ topology on 
the dual space of $(1 + d(\bp,\cdot)^2)^{-(N-1)} \: C_b(X)$.
The claim now
follows from Theorem \ref{properlimits}.
\end{proof}

\begin{example} \label{tcone}
We apply Corollary \ref{tangentcone} to tangent cones.
Suppose that $(X, d, \nu)$ is a complete measured locally compact
length space with nonnegative
$N$-Ricci curvature for some $N \in (2, \infty)$. Suppose that
$\supp(\nu) \: = \: X$.
For $i \: \ge \: 1$, put $(X_i, d_i) \: = \: (X, i \cdot d)$.
Given $\bp \in X$, let $\bp_i$ be the same point in $X_i$. 
Using Theorem \ref{thmBG} (as extended to the noncompact case),
after passing to a subsequence we may assume that
$\{(X_i, \bp_i)\}_{i=1}^\infty$ converges in the pointed
Gromov--Hausdorff topology to a tangent cone $(T_\bp X, o)$; 
see Corollary \ref{MGHcompactness}.
Let $\nu_i$ be the pushforward from $X$ to $X_i$, via the identity
map, of the measure
$\frac{\nu}{\nu({{B_{i^{-1}}(\bp))}}}$. After
passing to a further subsequence, we can assume that
$\{(X_i, \bp_i, \nu_i)\}_{i=1}^\infty$ converges in the
pointed measured Gromov--Hausdorff topology to a pointed measured length
space $(T_\bp X, o, \nu_\infty)$, where $\nu_\infty$ is a nonnegative
Radon measure on $T_\bp X$ that is normalized so that
$\nu_\infty({B_1(o)}) \: = \: 1$; see 
\cite[Section 1]{Cheeger-Colding (1997)}.
From Corollary \ref{tangentcone}, 
$(T_\bp X, \nu_\infty)$ has nonnegative $N$-Ricci curvature.
\end{example}

\section{Bibliographic notes on optimal transport} \label{bibnotes}
The following notes are by no means exhaustive, but may provide some
entry points to the literature.

Wasserstein was one of many authors who discovered, rediscovered or 
studied optimal transportation metrics \cite{Wasserstein (1969)}.
He was interested in the case 
when the cost coincides with the distance. 
Tanaka~\cite{Tanaka (1973)} may have been the first to
take advantage of geometric properties of $W_2$, in
his study of the Boltzmann equation. Accordingly, other names
could be used for $W_2$, 
such as Monge--Kantorovich distance or Tanaka distance.
The terminology ``Wasserstein distance'' was used by Otto and
coworkers, and naturally gave rise to the term ``Wasserstein space''.
Otto studied this metric space from a geometric point of view and 
showed that $P_2(\R^n)$ can be equipped with a formal
infinite-dimensional Riemannian metric, thereby
allowing insightful computations \cite{Otto (2001)}. 
He also showed that his Riemannian metric formally has nonnegative
sectional curvature.
Otto's motivation came from partial differential equations, and in
particular from earlier work by 
Brenier in fluid mechanics \cite{Brenier (1999)}. 
The formal gradient flow of a ``free energy'' functional on the
Wasserstein space was considered by
Jordan, Kinderlehrer and Otto~\cite{Jordan-Kinderlehrer-Otto (1998)}
and Otto \cite{Otto (2001)}.

The notion of displacement convexity was introduced by 
McCann~\cite{McCann (1997)}, and later refined to the notion of
$\lambda$-uniform displacement convexity. 
A formal differential calculus on $P_2(M)$, when $M$ is a smooth
Riemannian manifold, was described by Otto and 
Villani~\cite{Otto-Villani (2000)}.
It was ``shown'' that the entropy 
functional $\int_M \rho\log\rho \: \dvol_M$ 
is displacement convex on a manifold
with nonnegative Ricci curvature. Appendix~\ref{appformaldc} of the present 
paper follows up on the calculations in~\cite{Otto-Villani (2000)}.

Simultaneously, a rigorous theory of optimal
transport on manifolds was initiated by McCann~\cite{McCann (2001)} 
and further developed by 
Cordero-Erausquin, McCann and Schmuckenschl\"ager
\cite{Cordero-Erausquin-McCann-Schmuckenschlager (2001)}.
In particular, these authors prove
the implication $(1) \Rightarrow (4)$ of 
Theorem \ref{Riccimfolds}(a) of the present paper when
$\Psi$ is constant and $N = n$. The paper
\cite{Cordero-Erausquin-McCann-Schmuckenschlager (2001)}
was extended by von Renesse and Sturm~\cite{Sturm-von Renesse}, 
whose paper contains a
proof of the implications $(1) \Leftrightarrow (5)$ of
Theorem \ref{Riccimfolds}(b) of the present paper when
$\Psi$ is constant, and also indicates that
the condition (5) may make sense for some metric-measure spaces.
In a more recent contribution, which was done independently of the 
present paper,
Cordero-Erausquin, McCann and 
Schmuckenschl\"ager~\cite{Cordero-Erausquin-McCann-Schmuckenschlager (new)}
prove the implication $(1) \Rightarrow (5)$ of Theorem \ref{Riccimfolds}(b)
for general $\Psi$.

Connections between optimal transport and the theory of log
Sobolev inequalities and Poincar\'e inequalities were established
by Otto and Villani~\cite{Otto-Villani (2000)} and 
developed by many authors.
This was the starting point for 
Section~\ref{secineq} of the present paper. 
More information can be found in~\cite{Villani (2003)}.

A proof of a weak Bonnet--Myers theorem, based only on Riemannian growth
control and concentration estimates, was given by Ledoux~\cite{Ledoux (1999)}
as a special case of a more general result about the control of the
diameter of manifolds satisfying a log Sobolev inequality.
The simplified proof used in the present paper, based on a
Talagrand inequality, is taken from~\cite{Otto-Villani (2000)}.

\bibliographystyle{acm}

\begin{thebibliography}{9}

\bibitem{Alberti-Ambrosio (1999)}
G. Alberti and L. Ambrosio,
``A geometrical approach to monotone functions in $\R^n$'',
Math. Z. 230, p. 259-316 (1999)

\bibitem{Ambrosio-Gigli-Savare (2004)}
L. Ambrosio, N. Gigli and G. Savar\`e,
\underline{Gradient flows in metric spaces and in the space of
probability measures}, Lectures in Mathematics ETH Z\"urich,
Birkha\"user, Basel (2005)

\bibitem{Ambrosio-Gigli-Savare}
L. Ambrosio, N. Gigli and G. Savar\`e,
``Gradient flows with metric and differentiable structures, and 
applications to the Wasserstein space'',
Atti Accad. Naz. Lincei Cl. Sci. Fis. Mat. Natur. Rend. Lincei (9) Mat. Appl. 15, p.
327--343 (2004)

\bibitem{Toulouse (2000)}
C. An\'e, S. Blach\`ere, D.~Chafa\"{\i}, P.~Foug\`eres, I.~Gentil,
F.~Malrieu, C.~Roberto and G.~Scheffer,
\underline{Sur les in\'egalit\'es de Sobolev logarithmiques},
Panoramas et Synth\`eses 10, Soci\'et\'e Math\'ematique de France
(2000)

\bibitem{Bakry (1994)} D. Bakry,
``L'hypercontractivit\'e et son utilisation en th\'eorie des semigroupes'',
in \underline{Lectures on probability theory (Saint-Flour, 1992)},
Lecture Notes in Math. 1581,
Springer, Berlin, p. 1-114 (1994)

\bibitem{Bakry-Emery (1985)} D. Bakry and M. \'Emery,
``Diffusions hypercontractives'', in
\underline{S\'eminaire de probabilit\'es XIX}, 
Lecture Notes in Math. 1123, Springer, Berlin, p. 177-206 (1985)

\bibitem{Bernard-Buffoni (2004)} P. Bernard and B. Buffoni,
``Optimal mass transportation and Mather theory'', 
to appear, J. of European Math. Soc.

\bibitem{Bobkov-Gentil-Ledoux (2001)}
S. Bobkov, I. Gentil and M. Ledoux,
``Hypercontractivity of Hamilton-Jacobi equations'',
J. Math. Pures Appl. 80, p. 669-696 (2001)

\bibitem{Bobkov-Gotze (1999)}
S.G. Bobkov and F. G\"otze,
``Exponential integrability and transportation cost related to
logarithmic Sobolev inequalities'',
J. Funct. Anal. 163, p. 1-28 (1999)

\bibitem{Brenier (1999)} Y. Brenier,
``Minimal geodesics on groups of volume-preserving maps and generalized 
solutions of the Euler equations'',
Comm. Pure Appl. Math. 52, p. 411-452 (1999)

\bibitem{Burago-Burago-Ivanov (2001)}
D. Burago, Y. Burago and S. Ivanov,
\underline{A course in metric geometry}, 
Graduate Studies in Mathematics 33, 
American Mathematical Society, Providence (2001)

\bibitem{Burago-Gromov-Perelman (1992)} Y. Burago, M. Gromov and
G. Perelman, 
``A. D. Aleksandrov spaces with curvatures bounded below'',
Russian Math. Surveys 47, p. 1-58 (1992)

\bibitem{Buttazzo-Giaquinta-Hildebrandt (1998)}
G. Buttazzo, M. Giaquinta and S. Hildebrandt, 
\underline{One-dimensional variational problems},
Oxford Lecture Series in Mathematics and its Applications 15, 
Oxford University Press, New York (1998) 

\bibitem{Carrillo-McCann-Villani (2004)} J.A. Carrillo, R.J. McCann and 
C. Villani,
``Contractions in the $2$-Wasserstein length space and thermalization of
granular media'', Arch. Rational Mech. Anal. 179, p. 217-263 (2006)

\bibitem{Cheeger-Colding (1996)} J. Cheeger and T. Colding,
`` Lower bounds on Ricci curvature and the almost rigidity of warped 
products'', Ann. of Math. 144, p. 189-237 (1996)

\bibitem{Cheeger-Colding (1997)} J. Cheeger and T. Colding,
``On the structure of spaces with Ricci curvature bounded below I'', 
J. Diff. Geom. 46, p. 37-74 (1997)

\bibitem{Cheeger-Colding II (2000)} J. Cheeger and T. Colding,
``On the structure of spaces with Ricci curvature bounded below II'', 
J. Diff. Geom. 54, p. 13-35 (2000)

\bibitem{Cheeger-Colding (2000)} J. Cheeger and T. Colding,
``On the structure of spaces with Ricci curvature bounded below III'', 
J. Diff. Geom. 54, p. 37-74 (2000)

\bibitem{Cordero-Erausquin-McCann-Schmuckenschlager (2001)}
D. Cordero-Erausquin, R. McCann and M. Schmuckenschl\"ager,
``A Riemannian interpolation inequality \`a la Borell, Brascamp
and Lieb'', Invent. Math. 146, p. 219-257 (2001)

\bibitem{Cordero-Erausquin-McCann-Schmuckenschlager (new)}
D. Cordero-Erausquin, R. McCann and M. Schmuckenschl\"ager,
``Pr\'ekopa--Leindler type inequalities on Riemannian manifolds,
Jacobi fields, and optimal transport'', to appear in
Ann. Fac. Sci. Toulouse Math., available at
http://perso-math.univ-mlv.fr/users/cordero.dario/

\bibitem{Fukaya (1987)} K. Fukaya,
``Collapsing of Riemannian manifolds and eigenvalues of Laplace operator'',
Invent. Math. 87 , p. 517-547 (1987)

\bibitem{Gromov (1981)} M. Gromov, 
``Groups of polynomial growth and expanding maps'', Publ. Math. de l'IHES 53,
p. 53-78 (1981)

\bibitem{Gromov (1999)} M. Gromov,
\underline{Metric structures for Riemannian and non-Riemannian spaces},
Progress in Mathematics 152, Birkh\"auser, Boston (1999)

\bibitem{Grove-Petersen (1991)} K. Grove and P. Petersen,
``Manifolds near the boundary of existence'', 
J. Diff. Geom. 33, p. 379-394 (1991)

\bibitem{Jordan-Kinderlehrer-Otto (1998)}
R. Jordan, D. Kinderlehrer and F. Otto,
``The variational formulation of the Fokker-Planck equation'',
SIAM J. Math. Anal. 29, p. 1-17 (1998)

\bibitem{Ledoux (1999)}
M. Ledoux, ``Concentration of measure and logarithmic Sobolev inequalities'',
S\'eminaire de Probabilit\'es XXXIII, 
Lecture Notes in Mathematics 1709, p. 120-216
Springer-Verlag, Berlin (1999)

\bibitem{Ledoux (2001)}
M. Ledoux, \underline{The concentration of measure phenomenon},
American Mathematical Society, Providence (2001)

\bibitem{Lichnerowicz (1958)} 
A. Lichnerowicz,
\underline{G\'eom\'etrie des groupes de transformations},
Travaux et Recherches Math\'ematiques III, Dunod, Paris (1958)

\bibitem{Liese-Vajda (1987)} F. Liese and I. Vajda,
\underline{Convex statistical distances},
Teubner-Texte zur Mathematik 95,
BSB B. G. Teubner Verlagsgesellschaft, Leipzig (1987)

\bibitem{Lott (2003)} J. Lott,
``Some geometric properties of the Bakry-\'Emery-Ricci tensor'',
Comment. Math. Helv. 78, p. 865-883 (2003)

\bibitem{McCann (1997)} R.J. McCann,
``A convexity principle for interacting gases'',
Adv. Math. 128, p. 153-179 (1997)

\bibitem{McCann (2001)} R.J. McCann,
``Polar factorization of maps on Riemannian manifolds'',
Geom. Funct. Anal. 11, p. 589-608 (2001)

\bibitem{Otsu (1997)} Y. Otsu,
``Differential geometric aspects of Alexandrov spaces'',
in
\underline{Comparison geometry}, MSRI Publ. 30,
Cambridge Univ. Press, Cambridge, p. 135-148 (1997) 

\bibitem{Otsu-Shioya (1994)} Y. Otsu and T. Shioya,
``The Riemannian structure of Alexandrov spaces'',
J. Diff. Geom. 39, p. 629-658 (1994)

\bibitem{Otto (2001)} F. Otto,
``The geometry of dissipative evolution equations: the porous medium 
equation'', Comm. Partial Differential Equations 26, p. 101-174 (2001)

\bibitem{Otto-Villani (2000)} F. Otto and C. Villani,
``Generalization of an inequality by Talagrand, and links with the
logarithmic Sobolev inequality'',
J. Funct. Anal. 173, p. 361-400 (2000)

\bibitem{Perelman} G. Perelman,
``DC Structure on Alexandrov space with curvature bounded below'',
unpublished preprint

\bibitem{Perelman-Petrunin} G. Perelman and A. Petrunin,
``Quasigeodesics and gradient curves in Alexandrov spaces'',
unpublished preprint

\bibitem{Qian (1997)} Z. Qian,
``Estimates for weighted volumes and applications'',
Quart. J. Math. Oxford 48, p. 235-242 (1997)

\bibitem{Sturm-von Renesse} 
M.-K. von Renesse and K.-T. Sturm,
``Transport inequalities, gradient estimates and Ricci curvature'',
Comm. Pure Appl. Math. 68, p. 923-940 (2005)

\bibitem{Sturm1} K.-T. Sturm, ``Convex functionals of probability
measures and nonlinear diffusions on manifolds'', 
J. Math. Pures. Appl. 84, p. 149-168 (2005)

\bibitem{Sturm2} K.-T. Sturm, ``On the geometry of metric measure
spaces'', to appear, Acta Math.

\bibitem{Tanaka (1973)}
H. Tanaka,
``An inequality for a functional of probability
distributions and its application to Kac's
one-dimensional model of a Maxwellian gas'',
Z. Wahrscheinlichkeitstheorie und Verw. Gebiete 27,
p. 47-52 (1973)

\bibitem{Villani (2003)}
C. Villani,
\underline{Topics in optimal transportation},
Graduate Studies in Mathematics 58, 
American Mathematical Society, Providence (2003)

\bibitem{Wasserstein (1969)} L. Wasserstein, ``Markov processes over
denumerable products of spaces describing large systems of automata'',
Problems of Information Transmission 5, p. 47-82 (1969)

\bibitem{Zimmer (1984)}
R. Zimmer,
\underline{Ergodic theory and semisimple groups}, 
Monographs in Mathematics 81, 
Birkh\"auser, Basel (1984)
\end{thebibliography}

\end{document}